\documentclass{amsart}
\usepackage{amssymb,latexsym}

\newcommand{\light}[1]{{#1}} 
\newcommand{\bulletcolor}[1]{{#1}} 

\numberwithin{equation}{section}

\newtheorem{theorem}{Theorem}[section]
\newtheorem{proposition}[theorem]{Proposition}
\newtheorem{corollary}[theorem]{Corollary}
\newtheorem{definition}[theorem]{Definition}
\newtheorem{conjecture}[theorem]{Conjecture}

\newtheorem{lemma}[theorem]{Lemma}

\newtheorem{example}[theorem]{Example}

\def\ZZ{\mathbb{Z}}
\def\CC{\mathbb{C}}
\def\RR{\mathbb{R}}

\def\gg{\mathfrak{g}}
\def\nn{\mathfrak{n}}
\def\hh{\mathfrak{h}}

\def\A{\mathcal{A}}
\def\M{\mathcal{M}}
\def\B{\mathcal{B}}
\def\C{\mathcal{C}}
\def\D{\mathcal{D}}

\def\Arr{\mathrm{Arr}}
\def\ii{\mathbf{i}}
\def\jj{\mathbf{j}}

\def\l{\ell}

\def\<{\langle}
\def\>{\rangle}

\def\circi{{\bigcirc\!\!\!\! i\,}}
\def\circj{{\bigcirc\!\!\!\! j\,}}

\def\wnot{w_\mathrm{o}}

\newcommand{\circled}[1]{{\,\bigcirc\hspace{-.123in} #1\,}}

\newcommand{\beal}{\begin{eqnarray}\begin{array}{l} }
\newcommand{\bear}{\begin{eqnarray}\begin{array}{r} }
\newcommand{\beac}{\begin{eqnarray}\begin{array}{c} }
\newcommand{\bealn}{\[\begin{array}{l} }
\newcommand{\bearn}{\[\begin{array}{r} }
\newcommand{\beacn}{\[\begin{array}{c} }
\newcommand{\eea}{\end{array}\end{eqnarray}}
\newcommand{\eean}{\end{array}\]}

\newcommand{\beq}{\begin{equation} }
\newcommand{\eeq}{\end{equation} }

\newcommand{\mat}[4]{\left(\!\!\begin{array}{cc}
#1 & #2 \\
#3 & #4 \\
\end{array}\!\!\right)}

\newcommand{\doublesubscript}[3]{
\displaystyle\mathop{\displaystyle #1_{#2}}_{#3}
}

\begin{document}

\title{{Double Bruhat cells and total positivity}}

\author{Sergey Fomin}
\address{Department of Mathematics, Massachusetts Institute of
  Technology, Cambridge, Massachusetts 02139}
\email{fomin@math.mit.edu}

\author{Andrei Zelevinsky}
\address{Department of Mathematics, 
Northeastern University, Boston, Massachusetts 02115} 
\email{andrei@neu.edu}

\thanks{The authors were supported in part 
by NSF grants \#DMS-9400914, \#DMS-9625511, and \#DMS-9700927,
and by MSRI (NSF grant \#DMS-9022140).}


\date{February 10, 1998}


\maketitle

\tableofcontents

\section{Introduction} 
\label{sec:introduction}

The classical theory of total positivity studies matrices 
whose all minors are nonnegative.
Recently, G.~Lusztig~\cite{lusztig-reductive} 
extended this classical subject 
(pioneered in the 1930s by Gantmacher, Krein, and Schoenberg)
by introducing the totally nonnegative variety $G_{\geq 0}$ in an arbitrary 
reductive group~$G$.
Lusztig's study, motivated by surprising connections he discovered between
total positivity and his theory of canonical bases for quantum groups, 
was mainly focused on the structure of the intersection
$G_{\geq 0} \cap N$, where $N$ is the unipotent radical of a Borel subgroup
$B$ in~$G$.
This structure was also the main object of study in~\cite{BFZ} 
(for type $A$) and \cite{BZ} (for an arbitrary semisimple group).
In this paper we extend the results of \cite{BFZ,BZ} to the whole variety
$G_{\geq 0}$.

We will try to make the point that the natural framework
for the study of $G_{\geq 0}$ is provided by the decomposition of
$G$ into the disjoint union  
of \emph{double Bruhat cells} $G^{u,v} = B u B  \cap B_- v B_-$;
here $B$ and $B_-$ are two opposite Borel subgroups in~$G$,
and $u$ and $v$ belong to the Weyl group $W$ of~$G$.
We believe these double cells to be a very interesting object of study
in its own right.
The term ``cells'' might be misleading: in fact,
the topology of $G^{u,v}$ is in general quite nontrivial.
(In some special cases, the ``real part'' of $G^{u,v}$ was studied
in~\cite{rietsch,SSV}.
V.~Deodhar~\cite{deodhar} studied the intersections
$B u B  \cap B_- v B$ whose properties are very different 
from those of~$G^{u,v}$.)

We study a family of birational parametrizations of $G^{u,v}$, 
one for each reduced expression $\ii$ of the element
$(u,v)$ in the Coxeter group $W \times W$.
Every such parametrization can be thought of as a system of local 
coordinates in $G^{u,v}$.
We call these coordinates the \emph{factorization parameters}
associated to $\ii$. 
They are obtained by expressing a generic element $x \in G^{u,v}$
as an element of the maximal torus $H = B \cap B_-$ multiplied by
the product of elements of various one-parameter subgroups in $G$
associated with 
simple roots and their negatives; the reduced expression 
$\ii$ prescribes the order of factors in this product. 
The main technical result of this paper
(Theorem~\ref{th:t-through-x}) is an explicit formula 
for these factorization parameters as rational functions on
the double Bruhat cell $G^{u,v}$.

Theorem~\ref{th:t-through-x} is formulated in terms of 
a special family of regular functions $\Delta_{\gamma, \delta}$ 
on the group $G$.
These functions are suitably normalized matrix coefficients corresponding 
to pairs of extremal weights $(\gamma, \delta)$ in some fundamental
representation of~$G$.
Again, we believe these functions to be a very interesting object
that deserves a study of its own.
For the type $A$, they specialize to the minors of a matrix,
and their properties are of course developed in great detail.
It would be very interesting to extend the main body of the classical
theory of determinantal identities to the family of functions 
$\Delta_{\gamma, \delta}$.
In this paper, we make the first steps in this direction
(see especially Theorems~\ref{th:minors-Plucker}
and \ref{th:minors-Dodgson} below).

Returning to total positivity, our explicit formulas for factorization 
parameters allow us to obtain a family of \emph{total positivity criteria},
each of which efficiently tests whether a given element~$x$ from an 
arbitrary double Bruhat cell $G^{u,v}$ is totally nonnegative.
More specifically, each of our criteria consists in
verifying whether $x$ satisfies a system of 
inequalities of the form $\Delta_{\gamma, \delta}(x) > 0$,
the number of these inequalities being equal to the dimension 
of $G^{u,v}$. 

As in \cite{BFZ,BZ}, the main algebraic relations involving  
factorization parameters and generalized minors 
$\Delta_{\gamma, \delta}$ can be written in a ``subtraction-free'' form,
and thus the theory can be developed over an arbitrary \emph{semifield}.
The readers familiar with \cite{BFZ,BZ} will have no trouble
extending the corresponding results there 
(cf.~\cite[Section~2]{BFZ}) to the more general context of this paper.
We do not pursue this path here since at the moment we have not 
developed applications of this more general setup.

In Section~\ref{sec:main} we give precise formulations of 
our main results.
Their proofs are given in Sections~\ref{sec:prerequisites}
and~\ref{sec:proofs-general}.

The last Section~\ref{sec:gln} contains applications of our theory
to the case of the \emph{general linear group}.
The case $G=GL_n$ is treated separately for a number of reasons.
First, $GL_n$ is not a semisimple group
(although everything reduces easily to $SL_n$).
Furthermore, the questions that we consider become some  
very natural linear-algebraic questions whose understanding does not 
require any Lie-theoretic background.
For instance, the factorization parameters become the parameters 
in factorizations of a square matrix into the smallest possible
number of elementary Jacobi matrices. 
To our surprise, our main results seem to be new even in this
case. In Section~\ref{sec:gln}, we tried to present them in an
elementary form, making this section as self-contained as possible.
Last but not least, our results 
in the $GL_n$ case have a particularly transparent formulation 
in the language of pseudoline arrangements. 

\medskip

\textsc{Acknowledgments.} 
Part of this paper was written when the authors were participating
in the special program ``Combinatorics'' at MSRI in Berkeley
in Spring 1997. 
The second author (A.Z.) gratefully acknowledges the hospitality 
of his colleagues in Buenos Aires (Alicia Dickenstein and Fernando 
Cukierman) and Strasbourg (Peter Littelmann and Olivier Mathieu),
where he worked on parts of this paper; these visits were 
supported by the University of Buenos Aires and CNRS, France.
A large part of our computations in Section~\ref{sec:gln} 
were performed with \texttt{Maple}.

\section{Main results}
\label{sec:main}

\subsection{Semisimple groups}
\label{sec:lie}

We begin by introducing general terminology and notation 
(mostly standard) for semisimple Lie groups and algebras
(cf., e.g., \cite{springer}).
Let $\gg$ be a \emph{semisimple complex Lie algebra} of rank~$r$ 
with the \emph{Cartan decomposition} $\gg = \nn_- \oplus \hh \oplus \nn$. 
Let $e_i, h_i, f_i\,$, for $i=1, \ldots, r$, be the standard 
generators of~$\gg$, and let $A = (a_{ij})$ be the \emph{Cartan matrix}.
Thus $a_{ij} = \alpha_j (h_i)$, where 
$\alpha_1, \ldots, \alpha_r \in \hh^*$ are the \emph{simple roots} of~$\gg$.
Let $G$ be a simply connected complex Lie group with 
the Lie algebra $\gg$. 
Let $N_-$, $H$ and $N$ be closed subgroups of $G$ with Lie algebras 
$\nn_-$, $\hh$ and $\nn$, respectively.
Thus $H$ is a \emph{maximal torus},  
and $N$ and $N_-$ are two opposite \emph{maximal unipotent subgroups} of~$G$.  
Let $B_- = HN_-$ and $B = HN$ be the corresponding pair of opposite 
\emph{Borel subgroups}. 
For $i=1,\dots,r$ and $t \in \CC$, we write
\beal\label{eq:x,y}
x_i (t) = \exp(t e_i)\ , \quad
x_{\overline i} (t) = \exp(t f_i) \ ,
\eea
so that $t \mapsto x_i (t)$ (resp.\ $t \mapsto x_{\overline i} (t)$) 
is a one-parameter subgroup in $N$ (resp.\ in $N_-$).
We prefer the notation~$x_{\overline i} (t)$ to 
the usual~$y_i(t)$, for 
the reasons that will become clear later.
It will be convenient to denote $[1,r] = \{1, \ldots, r\}$
 and $[\overline 1, \overline r] = \{\overline 1, \ldots, \overline r\}$.

The \emph{weight lattice} $P$ is the set of all weights $\gamma \in \hh^*$ 
such that $\gamma (h_i) \in \ZZ$ for all $i$. 
The group $P$ has a $\ZZ$-basis formed by the \emph{fundamental weights}  
$\omega_1, \ldots, \omega_r$ defined by $\omega_i (h_j) = \delta_{ij}$.  
Every weight $\gamma \in P$ gives rise to a multiplicative character 
$a \mapsto a^\gamma$ of the maximal torus $H$; this character 
is given by $\exp (h)^\gamma = e^{\gamma (h)}  \,\, (h \in \hh)$. 
 
The \emph{Weyl group} $W$ of $G$ is defined by $W = {\rm Norm}_G (H)/H$. 
The action of~$W$ on~$H$ by conjugation gives rise to the action 
of $W$ on the weight lattice $P$ given by
\beal
\label{eq:W-action}
a^{w (\gamma)} = (w^{-1} a w)^\gamma  \quad
(w \in W, \, a \in H, \, \gamma \in P) \ .
\eea
As usual, we identify $W$ with the corresponding group of 
linear transformations of~$\hh^*$.  
The group $W$ is a \emph{Coxeter group} generated by \emph{simple reflections} 
$s_1, \ldots, s_r$ given by 
$s_i (\gamma) = \gamma - \gamma (h_i) \alpha_i\,$, for 
$\gamma \in \hh^*$.  

A \emph{reduced word} for $w \in W$ is a sequence of indices 
$\ii = (i_1, \ldots, i_m)$ of shortest possible length~$m$
such that $w = s_{i_1} \cdots s_{i_m}\,$. 
The number~$m$ is denoted by $\l(w)$ and  
is called the \emph{length} of~$w$.
The set of reduced words for~$w$ will be denoted by~$R(w)$. 
The Weyl group $W$ has the unique element $\wnot$ of maximal length,
and $\l (\wnot) = \l (w) + \l (w^{-1} \wnot)$
for any $w \in W$.

\subsection{Factorization problem}
\label{sec:factorization-problem}

Recall that the group $G$ has two \emph{Bruhat decompositions}, 
with respect to opposite Borel subgroups $B$ and $B_-$:
$$G = \bigcup_{u \in W} B u B = \bigcup_{v \in W} B_- v B_-  \ . $$
The \emph{double Bruhat cells}~$G^{u,v}$ are defined by 
$$G^{u,v} = B u B  \cap B_- v B_- \ ;$$
thus $G$ is the disjoint union of all $G^{u,v}$ for $(u,v)\in
W\times W$.

\begin{theorem}
\label{th:Guv-affine}
The variety $G^{u,v}$ is biregularly isomorphic to a Zariski
open subset of an affine space of dimension $r+\l(u)+\l(v)$. 
\end{theorem}

We will study a family of birational parametrizations of $G^{u,v}$. 
To describe these parametrizations, we will need the following
combinatorial notion.

A \emph{double reduced word} for the elements $u,v\in W$
is a reduced word for an element $(u,v)$ 
of the Coxeter group $W \times W$.
To avoid confusion, we will use the indices 
$\overline 1, \overline 2, \ldots, \overline r$ for the simple reflections 
in the first copy of $W$, and $1, 2, \ldots, r$ for the second copy.
A double reduced word for $(u,v)$ is nothing but a shuffle of a reduced
word for $u$ written in the alphabet 
$[\overline 1, \overline r]$ and a reduced word for $v$ written in 
the alphabet $[1,r]$.
We denote the set of double reduced words for $(u,v)$ by $R(u,v)$.

For any sequence $\ii= (i_1, \ldots, i_m)$ of indices 
from the alphabet $[1,r] \cup [\overline 1, \overline r]$,
consider the map $x_\ii: H \times \CC^m \to G$ defined by
\begin{equation}
\label{eq:productmap}
x_\ii (a; t_1, \ldots, t_m) = a\, x_{i_1} (t_1) \cdots x_{i_m} (t_m) \ ,
\end{equation} 
where we use the notation of~(\ref{eq:x,y}). 

Let $\CC_{\neq 0}$ denote the set of nonzero complex numbers. 

\begin{theorem}
\label{th:double}
For any $u,v\in W$ and $\ii= (i_1, \ldots, i_m)\in R(u,v)$, 
the map $x_\ii$ restricts to a biregular isomorphism between 
$H \times\CC_{\neq 0}^m$ and a Zariski open subset of the 
double Bruhat cell~$G^{u,v}$.
\end{theorem}

Thus $x_\ii$ gives rise to a birational isomorphism between
$H \times\CC^m$ and $G^{u,v}$.
We remark that this property holds if and only if  
$\ii$ is a double reduced word for~$(u,v)$. 

Theorem~\ref{th:double} tells that for a generic 
element $x\in G^{u,v}$ 
and any $\ii= (i_1, \ldots, i_m) \in R(u,v)$,
there are uniquely defined \emph{factorization parameters}
$a,t_1,\dots,t_m$ such that 
\beal\label{eq:x=axxx}
x=a\, x_{i_1} (t_1) \cdots x_{i_m} (t_m) \ .
\eea
One of our main results is the solution of the following
\emph{factorization problem}:
find explicit formulas for the inverse birational isomorphism 
$x_\ii^{-1}$ between $G^{u,v}$ and $H \times \CC^m$.
In other words, we express the factorization parameters
in terms of the element~$x$ and the double reduced
word~$\ii$ underlying the factorization. 
Our solution of the factorization problem generalizes 
Theorems~1.4 and~6.2 in~\cite{BZ}
(the case of $x\in N$),
which in turn generalize Theorems~1.4 and~5.4.2 in~\cite{BFZ}
(same, for type~$A$).

\subsection{Total positivity}
\label{sec:TPgeneral}

We will apply our solution of the factorization problem to the
study of total positivity. 
Following G.~Lusztig~\cite{lusztig-reductive}, 
let us define \emph{totally nonnegative elements} in~$G$. 
Let $H_{>0}$ be the subgroup of~$H$ consisting of all $a\in H$ such that
$a^\gamma \in \RR_{>0}$ for any weight $\gamma \in P$.  
(We denote by $\RR_{>0}$ the set of positive reals.) 
The set~$G_{\geq 0}$ of totally nonnegative elements is, by definition, 
the multiplicative semigroup in $G$ generated by 
$H_{>0}$ and the elements $x_i (t)$ and $x_{\overline i} (t)$,  
for $i\in [1,r]$ and $t\in\RR_{>0}$. 
It is easy to see that a totally nonnegative element $x\in G$  
can be represented as $x=x_\ii (a; t_1, \ldots, t_m)$,
for some sequence $\ii$, with all the $t_k$
positive and $a\in H_{>0}\,$. 
For the type~$A_r$ (i.e., for $G=SL_n(\CC)$, $n=r+1$), 
a theorem of C.~Loewner~\cite{loewner}, based on a result 
by A.~Whitney~\cite{whitney}, tells that the above definition of
total nonnegativity coincides with the usual
one~\cite{ando,karlin}: a matrix (with
determinant~1) is totally nonnegative if and only if all its
minors are nonnegative. 

The set $G_{\geq 0}$ is the disjoint union of the subsets 
$G^{u,v}_{> 0}$ obtained by intersecting it with double Bruhat
cells: 
$$G^{u,v}_{> 0} = G_{\geq 0} \cap G^{u,v} \ .$$
We call the~$G^{u,v}_{> 0}$ \emph{totally positive varieties}; 
they will be one of the main objects of study in this paper. 
The terminology is justified by the following observation made 
by Lusztig~\cite{lusztig-reductive}: in the special case
$G=SL_n(\CC)$ and $u=v=\wnot\,$, the variety $G^{u,v}_{> 0}$
is the set of all $n\!\times\! n$-matrices (with determinant~1)
which are (strictly) totally positive in the usual sense,
i.e., all their minors are positive. 

\begin{theorem}
\label{th:LusztigTP}
For any $u,v \in W$ and any double reduced word
$\ii \in R(u,v)$, the map $x_\ii$  restricts to a bijection
$H_{>0} \times \RR_{>0}^m \to G^{u,v}_{> 0}$. 
\end{theorem} 

Informally speaking, Theorem~\ref{th:LusztigTP} asserts that 
an element $x\in G^{u,v}$ 
is totally nonnegative if and only if for some 
(equivalently, any) 
double reduced word $\ii \in R(u,v)$,
the factorization parameters 
$a,t_1,\dots,t_m$ appearing in (\ref{eq:x=axxx}) 
are well-defined and positive.  
Thus the solution of the factorization problem will lead to a family of
total positivity criteria---one for each double reduced word. 

\subsection{Generalized minors}
\label{sec:minors}

The main ingredients of our answer to the factorization problem
are similar to those in \cite{BZ}: 
a family of regular functions
on $G$ generalizing minors of a square matrix,
and a biregular ``twist'' $G^{u,v} \to G^{u^{-1}, v^{-1}}$. 

We will denote by $G_0=N_-HN$ 
the set of elements $x\in G$ that admit Gaussian decomposition;
for the latter, we will use the notation
\beal
\label{eq:Gaussian-decomp}
x = [x]_- [x]_0 [x]_+ \ .
\eea

Following \cite{BZ}, for any fundamental weight $\omega_i\,$,
let $\Delta^{\omega_i}$  denote 
the regular function on $G$ whose restriction to the open set $G_0$
is given by
\begin{equation}
\label{eq:Delta-principal}
\Delta^{\omega_i} (x) = [x]_0^{\omega_i} \ .
\end{equation}
For the type $A_r\,$, the $\Delta^{\omega_i}(x)$ are the
\emph{principal minors} of a matrix~$x$. We will use the same
terminology in the general case as well. 

To define the analogues of arbitrary minors, we will need
two special representatives 
$\overline w,\overline{\overline {w}}\in G$ for any element $w \in W$. 
For a simple reflection~$s_i\,$, set 
$$
\overline {s_i} = \varphi_i \mat{0}{-1}{1}{0} \ ,\quad 
\overline{\overline {s_i}} = \varphi_i \mat{0}{1}{-1}{0} \ ,
$$
where $\varphi_i: SL_2 \to G$ is the group homomorphism
given by
\beal
\label{eq:phi_i}
\varphi_i \mat{1}{t}{0}{1} = x_i (t)\ , \quad 
\varphi_i \mat{1}{0}{t}{1} = x_{\overline i} (t) \ .
\eea
Alternatively, we could define 
\beal
\label{eq:s_i-via-x_i}
\overline {s_i}=x_i(-1)x_{\overline i}(1) x_i(-1)\ ,\quad
\overline {\overline {s_i}}=x_i(1)x_{\overline i}(-1) x_i(1)\ . 
\eea
It is known (and easy to check) that the families  
$\{\overline {s_i}\}$ and $\{\overline{\overline {s_i}}\}$
satisfy the braid relations in~$W$.
It follows that the representatives $\overline w$ and 
$\overline{\overline{w}}$ can be uniquely and unambigiously defined for any 
$w \in W$ by the condition that
\beal
\label{eq:braid} 
\overline {w' w''} = \overline {w'} \cdot \overline {w''} \ , \,\,
\overline {\overline {w' w''}}  = \overline {\overline {w'}} \cdot 
\overline {\overline {w''}}
\eea
whenever $\l (w' w'') = \l (w') + \l (w'')$.

\begin{definition}
\label{def:general minors}
{\rm
For $u,v \in W$, define a regular function 
$\Delta_{u \omega_i, v \omega_i}$ on $G$ by setting
\begin{equation}
\label{eq:Delta-general}
\Delta_{u \omega_i, v \omega_i} (x) = 
\Delta^{\omega_i} \left(\overline {\overline {u^{-1}}}  
   x \overline v\right)\ . 
\end{equation}
}\end{definition}

One has to check that this is well-defined, i.e., 
the right-hand side of (\ref{eq:Delta-general}) only depends on  
the weights $u \omega_i$ and $v \omega_i$, not on the particular choice 
of $u$ and~$v$. 
This is done in Section~\ref{sec:det-identities}
(cf.\ Proposition~\ref{prop:principal-minors-s_j}).

For the type $A_r$, the functions $\Delta_{u \omega_i, v \omega_i}(x)$
are the minors of a matrix $x$. 
In the general case, we will refer to them as \emph{generalized minors},
or simply as minors if there will be no danger of confusion. 

\subsection{The twist maps}
\label{sec:twist}

To define the twist maps, we will need the involutive automorphism
$x \mapsto x^\theta$ of the group~$G$ which is uniquely determined by 
\begin{equation}
\label{eq:theta}
a^\theta = a^{-1} \quad (a \in H) \ , \quad 
x_i (t)^\theta = x_{\overline i} (t) \ , \quad 
x_{\overline i} (t)^\theta = x_i (t) \ .
\end{equation} 
Notice that the involution~$\theta$ preserves total
nonnegativity. 
For the type $A_r\,$, if $x$ is a matrix with determinant~$1$, then 
the matrix $x^\theta$ is formed by signless
cofactors of~$x$; in other words, the $(i,j)$-entry of $x^\theta$
is simply the minor of~$x$ obtained by deleting the $i$th row and
the $j$th column. 

\begin{definition}
\label{def:twist}{\rm
For any $u,v\in W$, the twist map 
$\zeta^{u,v}:x\mapsto x'$ is defined by}
\beal
\label{eq:zeta-u,v-x}
x' = \left(
[\overline{\overline{u^{-1}}}x]_-^{-1}\,
\overline{\overline{u^{-1}}} \ x \ \overline{v^{-1}} \,
[x\overline{v^{-1}}]_+^{-1} 
\right)^\theta \ .  
\eea
\end{definition}

\begin{theorem}
\label{th:zeta-regularity}
The right-hand side of~{\rm (\ref{eq:zeta-u,v-x})} is well defined for 
any~$x\in G^{u,v}$, and the twist map $\zeta^{u,v}$
establishes a biregular isomorphism 
between $G^{u,v}$ and  $G^{u^{-1}, v^{-1}}$.
The inverse isomorphism is $\zeta^{u^{-1},v^{-1}}$.    
\end{theorem} 

The specific choice of representatives for $u$ and $v$ in 
(\ref{eq:zeta-u,v-x}) is essential for the following
important property.

\begin{theorem}
\label{th:zeta-positivity}
The twist map $\zeta^{u,v}$ restricts to a bijection 
$G^{u,v}_{>0} \to G^{u^{-1}, v^{-1}}_{>0}$.
\end{theorem}

\begin{example}
\label{example:SL2-cells}{\rm
Let $G=SL_2(\CC)$. 
Then $W=S_2$ consists of two permutations: $e$ and~$\wnot\,$;
thus $G$ is partitioned into four double Bruhat cells.
Table~\ref{table:SL2} shows the conditions under which a matrix
$
x\!=\! \mat{x_{1 1}}{x_{1 2}}{x_{2 1}}{x_{2 2}}
$
with determinant~1 belongs to each  of these cells,
or to the corresponding 
totally positive variety~$G^{u,v}_{>0}\,$. 
The table also shows the formulas defining each twist map~$\zeta^{u,v}$.  

\begin{table}[ht]
\[
\begin{array}{c|cccc}
& u=e      & u=e          & u=\wnot        & u=\wnot \\ 
& v=e      & v=\wnot        & v=e          & v=\wnot \\[.05in]
\hline
\\[-.15in]
G^{u,v} 
& \begin{array}{c}x_{12}=0\\     x_{21}=0\end{array} 
& \begin{array}{c}x_{12}\neq 0\\ x_{21}=0\end{array} 
& \begin{array}{c}x_{12}=0   \\   x_{21}\neq 0\end{array} 
& \begin{array}{c}x_{12}\neq 0 \\ x_{21}\neq 0\end{array} 
\\[.15in]
\hline
\\[-.15in]
G^{u,v}_{>0} 
& \begin{array}{c}x_{11}>0\\x_{12}=0\\ x_{21}=0\end{array} 
& \begin{array}{c}x_{11}>0\\x_{12}>0\\ x_{21}=0\end{array} 
& \begin{array}{c}x_{11}>0\\x_{12}=0\\x_{21}>0\end{array} 
& \begin{array}{c}x_{11}>0\\x_{12}>0\\x_{21}>0\end{array} 
\\[.22in]
\hline
\\[-.1in]
\zeta^{u,v}(x)
& \mat{x_{11}^{-1}}{0}{0}{x_{11}}
& \mat{x_{12}^{-1}}{x_{11}^{-1}}{0}{x_{12}}
& \mat{x_{21}^{-1}}{0}{x_{11}^{-1}}{x_{21}}
& \mat{x_{11}x_{12}^{-1}x_{21}^{-1}}{x_{21}^{-1}}{x_{12}^{-1}}{x_{22}}
\!\!\! \\
\end{array}
\]
\vspace{.1in}
\caption{Double Bruhat decomposition and the twist maps for $SL_2$}
\label{table:SL2}
\end{table}
}
\end{example}

\vspace{-.3in}

\subsection{Formulas for factorization parameters}
\label{sec:a,t_k}

To give explicit formulas for factorization parameters,
we will need some more notation.
First, we will write
\begin{equation}
\label{abs-i}
|i| = |\overline i| = i 
\end{equation}
for any $i \in [1,r]$. 
Let us fix a pair $(u,v) \in W \times W$ and a double reduced word 
$\ii = (i_1, \ldots, i_m)\in R(u,v)$. 
Recall that $\ii$ is a shuffle of a reduced
word for $u$ written in the alphabet 
$[\overline 1, \overline r]$ and a reduced word for $v$ written in the 
alphabet $[1,r]$.
In particular, the length 
$m$ of $\ii$ is equal to $\l (u) + \l (v)$.
We will add $r$~additional entries 
$i_{m+1}, \ldots, i_{m+r}$ at the end of $\ii$ by setting
\begin{equation}
\label{eq:extra-i}
i_{m+j} = \overline j \quad (j \in [1,r]) \ .
\end{equation}
For $k\in [1,m+r]$, we denote
\[
k^- = \max\{l:l<k, |i_l| = |i_k|\} 
\]
unless $|i_l| \neq |i_k|$ for all $l<k$, in which case we set 
$k^-=0$.
For $k\in [0,m+r]$, let 
\begin{equation}
\label{eq:epsilon}
\varepsilon_k =
\left\{
\begin{array}{ll}
1 &\textrm{if } i_k \in [1,r] \textrm{ or } k=0\,;\\[.1in]
0 &\textrm{if } i_k \in [\overline 1,\overline r] \,.
\end{array}
\right.
\end{equation}
For $k,l\in [0,m+r]$, let
\begin{equation}
\label{eq:chi}
\chi(k,l) =
\left\{
\begin{array}{ll}
1 &\textrm{if } k < l \textrm{ and } \varepsilon_k = \varepsilon_l 
\,;\\[.1in]
1/2 &\textrm{if } k = l \, ; \\[.1in]
0 &\textrm{if $k>l$ or $\varepsilon_k \neq \varepsilon_l$}\, . 
\end{array}
\right.
\end{equation}
For $k\in [1,m+r]$, we denote 
\beal
\label{eq:v_leq k}
u_{\geq k} 
= s_{|i_m|}^{1 - \varepsilon_m} s_{|i_{m-1}|}^{1 - \varepsilon_{m-1}} 
\cdots s_{|i_{k}|}^{1 - \varepsilon_{k}} \ , \,\,
v_{<k} = s_{|i_1|}^{\varepsilon_1} s_{|i_2|}^{\varepsilon_2} \cdots 
s_{|i_{k-1}|}^{\varepsilon_{k-1}} \ ;
\eea
where, by convention, $u_{\geq k} = e$ and $v_{<k}=v$ for $k>m$.
(For example, if 
$\ii=
\overline 2~~1~~\overline 3~~3~~2~~\overline 1~~\overline 2~~1~\overline 1$,
then, say, $u_{\geq 7} = s_1 s_2 \,$ and $v_{<7} = s_1 s_3 s_2\,$.)
Finally, for each $k\in [1,m+r]$, 
let us define a regular function 
$\Delta_k = \Delta_{k,\ii}$ on $G$ by 
\begin{equation}
\label{eq:Q-factors}
\Delta_k = \Delta_{k,\ii} = 
\Delta_{u_{\geq k} \omega_{|i_k|}, v_{<k} \omega_{|i_k|}} \ .
\end{equation}

With all this notation in mind, we now formulate our first main
result: a solution to the factorization problem of
Section~\ref{sec:factorization-problem}. 

\begin{theorem}
\label{th:t-through-x}
Let $\ii = (i_1, \ldots, i_m)$ be 
a double reduced word for $(u,v)$, and suppose an element 
$x \in G^{u,v}$ can be factored as
$x=a\, x_{i_1} (t_1) \cdots x_{i_m} (t_m)$ 
with $a \in H$ and all $t_k$ nonzero complex numbers.
Then the factorization parameters 
are determined by the following formulas,
where we denote $x' = \zeta^{u,v} (x)$:
\begin{equation}
\label{eq:t-through-x}
t_k =  \prod_{l= 1}^{m+r} 
\Delta_{l,\ii} (x')^{(\chi(k,l^-) - \chi(k,l)) a_{|i_l|,|i_k|}} \ ; 
\end{equation}
\beal
\label{eq:a-via-Q} 
a^{\omega_i} = 
\displaystyle\prod_{1 \leq k \leq m+r, \ |i_k| = i} 
\Delta_{k,\ii} (x')^{\varepsilon_k - \varepsilon_{k^-}} \ .  
\eea 
\end{theorem}

Formulas~(\ref{eq:a-via-Q}) can be restated as
the following closed expression for the element~$a$: 
\begin{equation}
\label{eq:a-through-x}                                                       
a = \prod_{k=1}^{m+1} 
\left[\overline {\overline {u_{\geq k}^{-1}}} x' \overline {v_{<k}}\right]_0^
{\varepsilon_{k}-\varepsilon_{k-1}} \ . 
\end{equation}

Since (\ref{eq:t-through-x}) and (\ref{eq:a-via-Q}) express
the $m+r$ independent parameters $t_1,\dots,t_m$ and 
$ a^{\omega_1},\dots,a^{\omega_r}$ as Laurent monomials in the $m+r$ 
minors $\Delta_{1,\ii} (x'),  \dots, \Delta_{m+r,\ii} (x')$,
we obtain the following important corollary.

\begin{theorem}
\label{th:inverse-monomial}
Under the assumptions of Theorem~\ref{th:t-through-x},
the parameters $t_1,\dots,t_m$ and $ a^{\omega_1},\dots,a^{\omega_r}$
are related to the minors 
$\Delta_{1,\ii} (x'),  \dots, \Delta_{m+r,\ii} (x')$
by an invertible monomial transformation.
\end{theorem}

The inverse of this monomial transformation 
can be computed explicitly: one can show that it is given by 
\begin{equation}
\label{eq:inverse-monomial}
\Delta_{k,\ii} (x') =
a^{- u u_{\geq k} \omega_{|i_k|}} 
\doublesubscript{\prod}{l < k}{\varepsilon_l = 0}
t_l^{u_{\geq l}^{-1} u_{\geq k} \omega_{|i_k|} (h_{|i_l|})}
\doublesubscript{\prod}{l \geq k}{\varepsilon_l = 1}
t_l^{v_{< l+1}^{-1} v_{<k} \omega_{|i_k|} (h_{|i_l|})} \ .
\end{equation}
Since we will not use this formula, we will not prove it in this paper.

\subsection{Total positivity criteria}
\label{sec:TP-criteria}

Theorem~\ref{th:t-through-x}
implies a family of criteria for total positivity
that generalize the ones in~\cite{BFZ,BZ}.
Each of these criteria asserts that 
a point $x \in G^{u,v}$ belongs to the totally positive 
variety $G^{u,v}_{>0}$ if and only if a particular collection 
of $\l (u) + \l (v) + r$ minors evaluated at $x$ are all positive. 

For every double reduced word $\,\ii = (i_1, \ldots, i_m) \in R(u,v)$,
let $F(\ii)$ denote the following collection of $m+r$ minors:
\begin{equation}
\label{eq:F(i)} 
F(\ii)
=\{\Delta_{k, \ii^*} \,:\, k\in[1,m+r]\}
\end{equation}
(cf.\ (\ref{eq:Q-factors})),
where $\ii^* = (i_m, \dots, i_1) \in R(u^{-1}, v^{-1})$ 
is $\ii$ written backwards. 

\begin{theorem}
\label{th:TP-criteria}
Let $\ii\in R(u,v)$. An element $x \in G^{u,v}$ is totally nonnegative 
if and only if $\Delta(x)>0$ for any minor $\Delta\in F(\ii)$. 
\end{theorem}

As a consequence, for any two double reduced words $\ii,\ii'\in R(u,v)$,
the positivity of all minors from $F(\ii)$ at a given $x\in G^{u,v}$
is equivalent to the positivity of all minors from $F(\ii')$. 
This phenomenon has the following algebraic explanation.
Let 
\beal
\label{eq:F(u,v)}
F(u,v)=\displaystyle\bigcup_{\ii\in R(u,v)} F(\ii)\ .
\eea
This can be restated as 
\beal
\label{eq:F(u,v)-explicit}
F(u,v) = \{
\Delta_{u' \omega_i, v' \omega_i} \,:\, i\in [1,r] \,,\,
u' \preceq u \,,\, v' \preceq v^{-1} \}\ ,
\eea
where $u'\preceq u$ stands for $\l(u)=\l(u')+\l({u'}^{-1}u)$
(the weak order on~$W$).

\begin{theorem}
\label{th:TP-base}
For any $\ii\in R(u,v)$, the collection $F(\ii)$ is a 
transcendence basis for the field of rational
functions~$\CC(G^{u,v})$.  
Furthermore, every minor in $F(u,v)$ can be expressed as a ratio of two
polynomials in the variables $\Delta\in F(\ii)$
with nonnegative integer coefficients. 
\end{theorem}

Combining 
Theorems~\ref{th:double}, \ref{th:zeta-regularity},
and~\ref{th:inverse-monomial} yields the following result. 

\begin{theorem}
\label{th:Laurent}
For any $\ii\in R(u,v)$, 
every minor in $F(u,v)$ is a Laurent polynomial with integer
coefficients in the variables $\Delta\in F(\ii)$. 
\end{theorem}

We suggest the following common refinement of 
Theorems~\ref{th:TP-base} and~\ref{th:Laurent}.

\begin{conjecture}
\label{conj:Laurent} 
For any $\ii\in R(u,v)$, 
every minor in $F(u,v)$ is a Laurent 
polynomial in the variables $\Delta\in F(\ii)$
with nonnegative integer coefficients. 
\end{conjecture}

Note that Theorems~\ref{th:TP-base} and~\ref{th:Laurent} 
do not automatically imply Conjecture~\ref{conj:Laurent},
since there do exist subtraction-free rational 
expressions that are Laurent polynomials although \emph{not} 
with nonnegative coefficients 
(for example, think of $(p^3+q^3)/(p+q)=p^2-pq+q^2$).

It is not hard to derive the following special case of
Conjecture~\ref{conj:Laurent} from~\cite[Theorem~3.7.4]{BFZ}.

\begin{theorem}
Conjecture~\ref{conj:Laurent} holds for $G=SL_{r+1}\,$,
when either $u$ or $v$ is the identity element~$e$. 
\end{theorem}

\subsection{Fundamental determinantal identities}

The subtraction-free rational expressions in
Theorem~\ref{th:TP-base} can be computed by an explicitly
described algorithm. 
This algorithm is based on repeated application of the following
generalized determinantal identities.

The first group of identities follow from~\cite[Corollary~6.6]{BZ}.
They correspond to pairs of simple roots that generate
a root subsystem of type $A_2$ or~$B_2$.
There are similar identities for subsystems of type $G_2$
(see~\cite[(4.8)--(4.11)]{BZ}) although we will not reproduce them here.

\begin{theorem}
\label{th:minors-Plucker}
Let $u, v \in W$ and $i, j \in [1,r]$. 
\vspace{.1in}

\noindent \emph{1.} 
If $a_{ij}\!=\!a_{ji}\!=\!-1$ 
and $\l (vs_i s_j s_i) = \l (v) + 3$, then
\[
\Delta_{u \omega_i, vs_i \omega_i} \Delta_{u \omega_j, vs_j
  \omega_j} = 
\Delta_{u \omega_i, v\omega_i} \Delta_{u \omega_j, vs_i s_j \omega_j} + 
\Delta_{u \omega_i, vs_j s_i \omega_i} \Delta_{u \omega_j, v \omega_j} \ . 
\]
\noindent \emph{2.} 
If $a_{ij}\! =\! -2$, $a_{ji}\! =\! -1$, and 
$\l (vs_i s_j s_i s_j) = \l (v) + 4$, then 
\bealn
\Delta_{u \omega_j, vs_i s_j\omega_j} 
\Delta^2_{u \omega_i, vs_j s_i \omega_i}
\Delta_{u \omega_j, vs_j \omega_j} 
\\[.1in]
=
\Delta_{u \omega_j, vs _j s_i s_j \omega_j}
\Delta^2_{u \omega_i, v s_j s_i \omega_i} 
\Delta_{u \omega_j, v \omega_j} 
\\[.1in]
\qquad\qquad\qquad\qquad
+(\Delta_{u \omega_i, v \omega_i} 
\Delta_{u \omega_j, vs _j s_i s_j \omega_j}\! +\!
\Delta_{u \omega_i, v s_i s_j s_i \omega_i} 
\Delta_{u \omega_j, vs_j \omega_j})^2
\eean
and
\bealn
\Delta_{u \omega_i, vs_i \omega_i} \Delta_{u \omega_i, vs_j s_i \omega_i}
\Delta_{u \omega_j, vs_j \omega_j} \\[.1in]
=\Delta^2_{u \omega_i, v s_j s_i \omega_i} 
 \Delta_{u \omega_j, v \omega_j} 
\\[.1in]
\qquad\qquad
+\Delta_{u \omega_i, v \omega_i} 
(\Delta_{u \omega_i, v \omega_i} 
 \Delta_{u \omega_j, vs _js_i s_j \omega_j} 
+\Delta_{u \omega_i, v s_i s_j s_i \omega_i} 
 \Delta_{u \omega_j, vs_j \omega_j})\ .
\eean
\noindent \emph{3.} 
Each of the above identities has a companion identity,
obtained by ``transposing'' all participating minors, i.e., 
by replacing every $\Delta_{\gamma,\delta}$ by $\Delta_{\delta,\gamma}\,$.
\end{theorem}

We also make use of the following new identity. 

\begin{theorem}
\label{th:minors-Dodgson}
Suppose $u,v \in W$ and $i \in [1,r]$
are such that $\l (us_i) = \l (u) + 1$ and $\l (vs_i) = \l (v) + 1$.
Then
\beal
\label{eq:minors-Dodgson}
\Delta_{u \omega_i, v \omega_i} \Delta_{us_i \omega_i, v s_i \omega_i}
= \Delta_{us_i \omega_i, v \omega_i} \Delta_{u \omega_i, v s_i \omega_i}
+ \prod_{j \neq i} \Delta_{u \omega_j, v \omega_j}^{- a_{ji}} \ .
\eea
\end{theorem}

The proof of Theorem~\ref{th:minors-Dodgson} is given in
Section~\ref{sec:det-identities}. 

For the type~$A_r\,$, the identities of 
Theorems~\ref{th:minors-Plucker} and~\ref{th:minors-Dodgson} 
become certain 3-term determinantal identities
known since early 19th century. We discuss their attribution
in Section~\ref{sec:applications-to-TP}.   


\section{Preliminaries}
\label{sec:prerequisites}

In what follows, we retain the notation and terminology introduced in 
Section~\ref{sec:main}.

\subsection{Involutions}
\label{sec:involutions}

Following \cite{BZ}, we define involutive anti-automorphisms
$x \mapsto x^T$ (the ``transpose'')
and $x \mapsto x^\iota$ of the group $G$ by setting
\begin{equation}
\label{eq:T}
a^T = a \quad (a \in H) \ , \quad x_i (t)^T = x_{\overline i} (t) \ , 
\quad x_{\overline i} (t)^T = x_i (t) 
\end{equation} 
and
\begin{equation}
\label{eq:iota}
a^\iota = a^{-1} \quad (a \in H) \ , \quad x_i (t)^\iota = x_i (t) \ , 
\quad x_{\overline i} (t)^\iota = x_{\overline i} (t) \ .
\end{equation}
These two involutive anti-automorphisms commute with each other 
and with the involutive anti-automorphism $x \mapsto x^{-1}$ of~$G$.
Hence these three maps generate the group isomorphic to $(\ZZ/2\ZZ)^3$;
in particular, any composition of them is again an involution. 
Notice that the involutions $x \mapsto x^T$ and $x \mapsto x^\iota$
preserve total nonnegativity, while $x \mapsto x^{-1}$ does not.
Informally, $x^\iota$ is a ``totally nonnegative version''
of~$x^{-1}$.

In the notation just introduced, the involution $x\mapsto
x^\theta$ that was defined by (\ref{eq:theta}) is given by
$x^\theta = (x^\iota)^T = (x^T)^\iota$.

The involutions $x \mapsto (x^{-1})^\iota$ and $x \mapsto x^T$
obviously preserve $G_0 = N_- H N$, and we have 
\beal
\label{eq:T-Gauss}
[(x^{-1})^\iota]_0  = [x^T]_0 = [x]_0 \ .
\eea

All three involutions $x \mapsto x^{-1}$, $x \mapsto x^T$ and 
$x \mapsto x^\iota$ act on $W$ by $w \mapsto w^{-1}$.
The relations between these involutions and the special
representatives introduced in Section~\ref{sec:minors}
are summarized in the following proposition. 

\begin{proposition}
\label{pr:bar-involutions}
We have
\beal
\label{eq:barT}
\overline w^{\ -1} = \overline w^{\,T}
                   =\overline{\overline{w^{-1}}}
                   = \overline{\overline{w}}^{\,\iota} 
, \quad 
\overline {w^{-1}} = \overline w^{\,\iota}
                   = \overline{\overline{w}}^{\,-1} 
                   = \overline{\overline{w}}^{\,T} \ .
\eea
\end{proposition}

\proof
Since all three involutions are 
antiautomorphisms, it is enough to check (\ref{eq:barT}) for $w = s_i$,
in which case it follows by a calculation in~$SL_2\,$. \endproof

\subsection{Commutation relations}
\label{sec:commutation}

For convenience of exposition, we collect here some known
commutation relations in $G$ that will be used in our proofs. 
Recall that $x_i (t)$ and $x_{\overline i} (t)$ are defined 
by~(\ref{eq:x,y}), and $\alpha_1, \ldots, \alpha_r$ are the simple roots
of~$\gg$.

First of all, for every $a \in H$, we have 
\beal
\label{eq:H-conjugation}
a x_i (t) = x_i (a^{\alpha_i} t) a \ , \quad
a x_{\overline i} (t) = x_{\overline i} (a^{- \alpha_i} t) a \ .
\eea

The following relations between the elements $x_i (t)$ 
can be found, e.g., in~\cite[Section~3]{BZ} 
(some of them appeared earlier in~\cite{lusztig-quantum}).
If $a_{ij} = a_{ji} = 0$, then
\beal\label{eq:2-commutation}
x_i (t_1) x_j (t_2) = x_j (t_2) x_i (t_1) 
\eea
for any $t_1$ and $t_2$. 
If $a_{ij} = a_{ji} = -1$, then
\begin{equation}
\label{eq:3-commutation}
x_i (t_1) x_j (t_2) x_i (t_3)
= x_j \left(\displaystyle\frac{t_2 t_3}{t_1 + t_3}\right) 
x_i (t_1 + t_3) x_j\left(\displaystyle\frac{t_1 t_2}{t_1 + t_3}\right) 
\end{equation}
whenever $t_1 + t_3 \neq 0$. 
If $a_{ij} = -2$ and $a_{ji} = -1$, then
\begin{equation}
\label{eq:4-commutation}
x_i (t_1) x_j (t_2) x_i (t_3) x_j (t_4) = 
x_j (t_2 t_3^2 t_4 q^{-1}) 
x_i (q p^{-1}) 
x_j (p^2 q^{-1})
x_i (t_1 t_2 t_3 p^{-1}) \ , 
\end{equation}
where
$$ p = t_1 t_2 + (t_1+t_3) t_4, \,\, 
q = t_1^2 t_2 + (t_1+t_3)^2 t_4 \ ;$$
this relation holds whenever $p \neq 0$ and $q \neq 0$. 
In the case when $a_{ij} = -3, \,  a_{ji} = -1$
(i.e., when $\alpha_i$ and $\alpha_j$ generate a root subsystem of 
type $G_2$), there is also a relation similar to (\ref{eq:3-commutation})
and (\ref{eq:4-commutation}).
This relation is given in~\cite[(3.6) -- (3.10)]{BZ};
we will not reproduce it here. 
Each of the relations (\ref{eq:2-commutation})--(\ref{eq:4-commutation})  
has a counterpart for the elements
$x_{\overline i} (t)$; it can be obtained by applying 
the antiautomorphism $x \mapsto x^T$ (cf.~(\ref{eq:T})).

In conclusion, let us describe the commutation relations 
between the elements $x_i (t)$ and $x_{\overline j} (t')$.
If $i \neq j$, then $[e_i, f_j] = 0$ in $\gg$, hence 
\beal
\label{eq:2-commutation-mixed}
x_i (t) x_{\overline j} (t') = x_{\overline j} (t') x_i (t) 
\eea
for any $t$ and $t'$. 
To handle the case $i=j$, we will need the following notation.
For a nonzero $t\in\CC$ and $i\in [1,r]$, we denote 
\beal
\label{eq:t^h}
t^{h_i} = \varphi_i \mat{t}{0}{0}{t^{-1}} \ ,
\eea
where $\varphi_i: SL_2 \to G$ is defined by (\ref{eq:phi_i});
alternatively, $t^{h_i}$ is an element of $H$ uniquely determined by 
the condition that 
$(t^{h_i})^\gamma = t^{\gamma (h_i)}$ for any weight~$\gamma \in P$. 
Then we have
\beal
\label{eq:5-commutation-generic}
x_i (t) x_{\overline i} (t') = 
x_{\overline i} 
\left(\displaystyle\frac{t'}{1+tt'}\right) (1+tt')^{h_i} 
x_i \left(\displaystyle\frac{t}{1+tt'}\right) 
\eea
whenever $1 + tt' \neq 0$. 
This relation can be first checked for $SL_2$
by a simple matrix calculation, and then extended to $G$
by applying the homomorphism $\varphi_i\,$.
By the same method, we verify the relations 
\beal
\label{eq:5-commutation-special}
x_i (t) x_{\overline i} (-t^{-1}) = 
x_{\overline i} (t^{-1}) t^{h_i} \overline {\overline {s_i}}
\eea
and
\beal
\label{eq:commutation-s_i}
x_i (t) \overline {s_i} = 
x_{\overline i} (t^{-1}) t^{h_i} x_i (- t^{-1})\ , \\[.1in]
\overline {\overline {s_i}}\,  x_{\overline i} (t) = 
x_{\overline i} (-t^{-1}) t^{h_i} x_i (t^{-1}) \ .
\eea

\subsection{Generalized determinantal identities}
\label{sec:det-identities}

We start with some identities for the
``principal minors" $\Delta^{\omega_i}$.
The definition (\ref{eq:Delta-principal}) implies that,
for any $x \in G\,$, $x^- \in N_-\,$, $x^+ \in N$, and $a \in H$, we have
\beal
\label{eq:principal-minor-invariance}
\Delta^{\omega_i} (x^- x) = \Delta^{\omega_i} (x x^+) = 
\Delta^{\omega_i} (x) \ , \\[.1in]
\Delta^{\omega_i} (a x) = \Delta^{\omega_i} (x a) = 
a^{\omega_i} \Delta^{\omega_i} (x) \ .
\eea
In view of (\ref{eq:T-Gauss}), we also have
\beal
\label{eq:principal-minor-T}
\Delta^{\omega_i}( (x^{-1})^\iota)  = \Delta^{\omega_i} (x^T) = 
\Delta^{\omega_i} (x) \ .
\eea
The following property is less obvious.

\begin{proposition}
\label{pr:principal-minors-x_j}
For any $x \in G, \, j \neq i$, and $t \in \CC$, we have
\beal
\label{eq:principal-minors-x_j}
\Delta^{\omega_i} (x x_{\overline j}(t)) = 
\Delta^{\omega_i} (x_j (t) x) =  \Delta^{\omega_i} (x) \ .
\eea
\end{proposition}

\proof
It is possible to deduce the proposition from the commutation 
relations given in Section~\ref{sec:commutation}
but we prefer another proof based on representation theory.
The group $G$ acts by right translations in the space $\CC [G]$ 
of regular functions on~$G$.
It is well known that every $f \in \CC [G]$ generates a finite-dimensional
subrepresentation of $\CC [G]$. 
In view of (\ref{eq:principal-minor-invariance}), the function
$\Delta^{\omega_i}$ is a highest weight vector of weight $\omega_i$ in 
$\CC [G]$.
Since $\omega_i (h_j) = 0$ for $j \neq i$, it follows that
$\Delta^{\omega_i}$ has weight $0$ 
with respect to the subgroup $\varphi_j (SL_2)$ of $G$ 
(cf.~(\ref{eq:phi_i})).
Therefore, $\Delta^{\omega_i}$ generates a trivial representation
of $\varphi_j (SL_2)$.
In particular, 
$\Delta^{\omega_i} (x x_{\overline j}(t)) = \Delta^{\omega_i} (x)$,
as desired. 
The equality $\Delta^{\omega_i} (x_j (t) x) =  \Delta^{\omega_i} (x)$
now follows from (\ref{eq:principal-minor-T}).
\endproof

Our next proposition justifies the validity of 
Definition~\ref{def:general minors}. 

\begin{proposition}
\label{prop:principal-minors-s_j}
For any $x \in G$ and any $j \neq i$, we have
\beal
\label{eq:principal-minors-s_j}
\Delta^{\omega_i} (x \overline {s_j}) = 
\Delta^{\omega_i} (\overline {\overline {s_j}} x) =  \Delta^{\omega_i} (x) \ .
\eea
\end{proposition}

\proof
Follows from (\ref{eq:principal-minor-invariance}),
(\ref{eq:principal-minors-x_j}), and~(\ref{eq:s_i-via-x_i}). 
\endproof

The extension of principal minors from the open subset $G_0$ to the whole of
$G$ is given as follows. 
The Bruhat decomposition theorem implies 
that every $x \in G$ can be written as 
\begin{equation}
\label{eq:B_-WB_+}
x = x^- a \overline w x^+
\end{equation} 
for some $x^- \in N_-, \, a \in H, \, w \in W$, and $x^+ \in N$;
moreover, the elements $a \in H$ and $w \in W$
are uniquely determined by~$x$.  

\begin{proposition}
\label{pr:principal-minor-G}
If $x \in G$ is expressed in the form {\rm (\ref{eq:B_-WB_+})}, 
then 
\beal
\label{eq:principal-minor_G}
\Delta^{\omega_i} (x) =               
\left\{
\begin{array}{ll}
a^{\omega_i} &\textrm{if } w \omega_i = \omega_i; \\[.05in]
0 \ , & \textrm{otherwise} \, .
\end{array}
\right.
\eea
\end{proposition} 

\proof
By (\ref{eq:principal-minor-invariance}), we have
$$\Delta^{\omega_i} (x^- a \overline w x^+) = 
a^{\omega_i} \Delta^{\omega_i} (\overline w) \ .$$
Thus, to prove (\ref{eq:principal-minor_G}) we only need to show that
\beal
\label{eq:principal-minor_W}
\Delta^{\omega_i} (\overline w) = 
\left\{
\begin{array}{ll}
1 &\textrm{if } w \omega_i = \omega_i;\\[.05in]
0 \ , &\textrm{otherwise} \, .
\end{array}
\right.
\eea
The formula is obvious for $w = e$, the identity element of $W$.
Hence we can assume that $\l (w) \geq 1$ and write 
$w$ as $u s_j$  for some $u \in W$ and $j \in [1,r]$ with 
$\l (u) = \l (w) - 1$. 
Since $\Delta^{\omega_i}$ is a regular function on $G$, we have
\beal
\label{eq:perturbation}
\Delta^{\omega_i} (\overline w) = 
\Delta^{\omega_i} (\overline u \, \overline {s_j}) = \lim_{t \to 0} 
\Delta^{\omega_i} (\overline u x_j (t) \overline {s_j}) \ .
\eea
Substituting into (\ref{eq:perturbation}) the expression for 
$x_j (t) \overline {s_j}$ given by (\ref{eq:commutation-s_i}) 
and using (\ref{eq:principal-minor-invariance}), we obtain
$$
\Delta^{\omega_i} (\overline w) = \lim_{t \to 0} 
t^{\omega_i (h_j)} \Delta^{\omega_i} (\overline u x_{\overline j} (t^{-1}))
\ . $$
Since $\l (u) = \l (w) - 1$, the root $u (\alpha_j)$ is positive,
implying that 
$\overline u x_{\overline j} (t^{-1}) \overline u^{-1} \in N_-\,$. 
Again using (\ref{eq:principal-minor-invariance}), we obtain
$$
\Delta^{\omega_i} (\overline w) = \lim_{t \to 0} 
t^{\delta_{ij}} \Delta^{\omega_i} (\overline u) \ .
$$
It follows that 
$$
\Delta^{\omega_i} (\overline w) = 
\left\{
\begin{array}{ll}
\Delta^{\omega_i} (\overline u) &\textrm{if } j \neq i;\\[.05in]
0 \ , &\textrm{otherwise} \, .
\end{array}
\right.
$$
This implies (\ref{eq:principal-minor_W}) by induction on $\l (w)$.
\endproof

\smallskip  

As a corollary, we obtain the following useful characterization
of the set~$G_0\,$.

\begin{corollary}
\label{cor:G_0-description}
An element $x \in G$ admits the Gaussian decomposition if and only if 
$\Delta^{\omega_i} (x) \neq 0$ for any $i \in [1,r]$.
\end{corollary}

In subsequent proofs, we will also make use
of the following identities. 

\begin{proposition}
\label{pr:minors-Gauss}
For any $x = [x]_- [x]_0 [x]_+ \in G_0$ and any $w \in W$, we have 
\beal
\label{eq:minors-Gauss1}
\Delta_{\omega_i, w \omega_i} ([x]_+) = 
\displaystyle\frac{\Delta_{\omega_i,w\omega_i}(x)}{\Delta^{\omega_i}(x)}
\ ,
\eea
\beal
\label{eq:minors-Gauss2}
\Delta_{w \omega_i, \omega_i} ([x]_-) = 
\displaystyle\frac{\Delta_{w \omega_i,\omega_i}(x)}{\Delta^{\omega_i} (x)} \ .
\eea
\end{proposition}

\proof
Using (\ref{eq:Delta-general}) and (\ref{eq:principal-minor-invariance}),
we obtain: 
$$
\Delta_{\omega_i, w \omega_i} (x) = 
\Delta^{\omega_i} (x \overline w) = 
[x]_0^{\omega_i} \Delta^{\omega_i} ([x]_+ \overline w) = 
\Delta^{\omega_i} (x) \Delta_{\omega_i, w \omega_i} ([x]_+) \ ,
$$
which proves (\ref{eq:minors-Gauss1}).
The proof of (\ref{eq:minors-Gauss2}) is similar.
\endproof

The transformation $(-\wnot)$ permutes fundamental weights.
We will use the notation $i \mapsto i^*$ for the induced 
permutation of the index set $[1,r]$, so that 
\beal
\label{eq:i^*}
\omega_{i^*} = - \wnot (\omega_i) \,\, (i \in [1,r]) \ . 
\eea

\begin{proposition}
\label{pr:minors-T}
For any $x \in G$ and any $u,v \in W$, we have
\beal
\label{eq:minors-T}
\Delta_{u\omega_i, v \omega_i} (x) = 
\Delta_{v \omega_i, u \omega_i} (x^T) =
\Delta_{v \wnot \omega_{i^*}, u \wnot \omega_{i^*}} (x^\iota)  \ .
\eea
\end{proposition}

\proof
Using (\ref{eq:principal-minor-T}) and (\ref{eq:barT}), we
obtain: 
$$
\Delta_{u \omega_i, v \omega_i} (x) = 
\Delta^{\omega_i} (\overline {\overline {u^{\ -1}}} x \overline v) = 
\Delta^{\omega_i} ((\overline {\overline {u^{\ -1}}} x \overline v)^T) = 
\Delta^{\omega_i} (\overline {\overline {v^{\ -1}}} x^T \overline u) =  
\Delta_{v \omega_i, u \omega_i} (x^T) \ ,
$$
which proves the first equality in (\ref{eq:minors-T}).

To prove the second equality, let us introduce the antiautomorphism
$\eta$ of $G$ by
$ \eta (x) = \overline {\overline {\wnot}} x^\iota \overline {\wnot}\,$.
Since $\eta$ preserves $H$ and interchanges $N_-$ and $N$,
it follows that $\eta$ preserves $G_0\,$, and 
$[\eta (x)]_0 = \eta ([x]_0) = 
\overline {\overline {\wnot}} [x]_0^{-1} \overline {\wnot}$
for any $x \in G_0\,$.
Using (\ref{eq:W-action}) and (\ref{eq:i^*}), we conclude that
$\Delta^{\omega_i} (\eta (x)) =  
\Delta^{\omega_{i^*}} (x)$.
Hence 
\bealn
\Delta_{u \omega_i, v \omega_i} (x) = 
\Delta^{\omega_i} (\overline {\overline {u^{\ -1}}} x \overline v) = 
\Delta^{\omega_{i^*}} (\eta (\overline {\overline {u^{\ -1}}} x \overline v))
=\Delta^{\omega_{i^*}} (\overline {\overline {\wnot}} \,
{\overline {\overline v}}^{-1} \ x^\iota \ 
{\overline {u^{-1}}}^{\ -1} \overline {\wnot}) \\[.1in]
= \Delta^{\omega_{i^*}} (\overline {\overline {\wnot v^{\ -1}}} x^\iota
\overline {u \wnot}) = 
\Delta_{v \wnot \omega_{i^*}, u \wnot \omega_{i^*}} (x^\iota) \ ,
\eean
as claimed.
\endproof


\emph{Proof of Theorem~\ref{th:minors-Dodgson}.}
First of all, since $\overline {\overline {s_i u^{-1}}} = 
\overline {\overline {s_i}} \, \overline {\overline {u^{-1}}}$
and $\overline {v s_i} = \overline v \, \overline {s_i}$, 
the definition (\ref{eq:Delta-general}) implies that it is enough
to prove (\ref{eq:minors-Dodgson}) in the case when $u = v = e$,
the identity element. 
Thus, we only need to show that
\beal
\label{eq:Dodgson-principal}
\Delta_{\omega_i, \omega_i} \Delta_{s_i \omega_i, s_i \omega_i}
- \Delta_{s_i \omega_i, \omega_i} \Delta_{\omega_i, s_i \omega_i} =
\prod_{j \neq i} \Delta_{\omega_j, \omega_j}^{- a_{ji}} \ .
\eea
As in the case of Proposition~\ref{pr:principal-minors-x_j},
our proof of (\ref{eq:Dodgson-principal}) will rely on  
representation theory. 
Consider the representation $\rho$ of the group $G \times G$ 
in $\CC [G]$ given by
$$
\rho (x_1, x_2) f (x) = f(x_1^T x x_2) \ .
$$
Let us denote the left- and right-hand sides of 
(\ref{eq:Dodgson-principal}) by $f_1$ and $f_2$, respectively.
We first verify that the function $f_2 \in \CC [G]$ 
has the following properties: 

\noindent (1) $f_2$ is a highest weight vector in the representation
$\rho$, i.e., it is invariant under the subgroup 
$N \times N \subset G \times G$;

\noindent (2) $f_2$ has weight $(\gamma, \gamma)$,
where $\gamma = 2 \omega_i - \alpha_i$; that is,
$\rho (a_1, a_2) f_2 = (a_1 a_2)^\gamma f_2$ for any $a_1, a_2 \in H$;

\noindent (3) $f_2 (e) = 1$
(here $e$ stands for the identity element of~$G$).

Property (3) is trivial, while (1) follows from 
(\ref{eq:principal-minor-invariance}). 
Also by (\ref{eq:principal-minor-invariance}), $f_2$ has weight 
$(- \sum_{j \neq i} a_{ji} \omega_j, \, 
- \sum_{j \neq i} a_{ji} \omega_j)$.
To prove (2), it is enough to show that 
$- \sum_{j \neq i} a_{ji} \omega_j = \gamma$; 
but this follows from the equality
\beal
\label{eq:weight}
\sum_{j \in [1,r]} a_{ji} \omega_j =  \alpha_i \ ,
\eea
which can be taken as a definition of the Cartan matrix.

Properties (1)--(3) uniquely determine the restriction of $f_2$
to~$G_0\,$. 
Since $G_0$ is dense in $G$, and $f_2$ is regular,
these properties uniquely determine $f_2\,$.
It remains to show that $f_1$ satisfies~(1)--(3).

The normalization condition (3) follows from 
Proposition~\ref{pr:principal-minor-G}; 
indeed, in view of (\ref{eq:principal-minor_G}), 
$$
\Delta_{\omega_i, \omega_i} (e) =  
\Delta_{s_i \omega_i, s_i \omega_i} (e) = 1 \ , \,\,
\Delta_{s_i \omega_i, \omega_i}(e) =  
\Delta_{\omega_i, s_i \omega_i} (e) = 0 \ .
$$

To prove that $f_1$ satisfies (2), notice that for any $u, \, v \in W$, 
the function $\Delta_{u \omega_i, v \omega_i}$
has weight $(u \omega_i, v \omega_i)$ (this follows from 
(\ref{eq:Delta-general}) and (\ref{eq:W-action})). 
Hence both summands in $f_1$ have weight 
$(\omega_i + s_i \omega_i, \omega_i + s_i \omega_i) = (\gamma, \gamma)$.

To prove that $f_1$ satisfies (1), we first notice that, in view of 
(\ref{eq:minors-T}), we have $f_1 (x^T) = f_1 (x)$ for any $x \in G$,
hence $\rho (x_1, x_2) f_1 = \rho (x_2, x_1) f_1$ for any $x_1, x_2 \in G$.
Therefore, it suffices to show that $f_1$ is invariant under 
the action of $N$ by right translations. 
Let $E_1, \ldots, E_r$ be the infinitesimal right translation operators
on $\CC [G]$ defined by
$$
E_j (f) (x) = \frac{d}{dt} f(x x_j (t))|_{t=0} \ ;
$$
each $E_j$ is a derivation of the ring $\CC [G]$.
It is enough to show that $E_j f_1 = 0$ for all~$j$.
If $j \neq i$, then $E_j$ annihilates
all four minors that appear in $f_1$ (this follows from the fact
that $\overline {s_i} x_j (t) \overline {s_i}^{\ -1} \in N$),
hence $E_j f_1 = 0$. 
It remains to prove that $E_i f_1 = 0$.
Clearly, we have
\beal
\label{eq:E_i = 0}
E_i \Delta_{\omega_i, \omega_i} = E_i \Delta_{s_i \omega_i, \omega_i} = 0
\ . \eea
We claim that
\beal
\label{eq:E_i}
E_i \Delta_{\omega_i, s_i \omega_i} = \Delta_{\omega_i, \omega_i}, \,\,
E_i \Delta_{s_i \omega_i, s_i \omega_i} = \Delta_{s_i \omega_i, \omega_i} 
\ . \eea
Combining (\ref{eq:E_i = 0}) and (\ref{eq:E_i}) and using the Leibniz rule,
we obtain
$$
E_i f_1 \!=\! 
E_i (\Delta_{\omega_i, \omega_i} \Delta_{s_i \omega_i, s_i \omega_i}
- \Delta_{s_i \omega_i, \omega_i} \Delta_{\omega_i, s_i \omega_i})\! =\! 
\Delta_{\omega_i, \omega_i} \Delta_{s_i \omega_i, \omega_i} - 
\Delta_{\omega_i, \omega_i} \Delta_{s_i \omega_i, \omega_i} \!=\! 0\,,
$$
as required.
We will deduce (\ref{eq:E_i}) from the following lemma which is
a standard fact in the representation theory of $SL_2$.

\begin{lemma}
\label{lem:SL_2}
Suppose $f \in \CC [G]$ is such that $E_i f = 0$ and 
$f(x t^{h_i}) = t^k f(x)$ for some $k \geq 0$. 
Let $f' \in \CC [G]$ be given by $f' (x) = f(x \overline {s_i})$.
Then $E_i^k (f') = k! f$. 
\end{lemma}

The first equality in (\ref{eq:E_i}) follows by applying this lemma to
$f = \Delta_{\omega_i, \omega_i}$ (in this case, $k = 1$ and 
$f' = \Delta_{\omega_i, s_i \omega_i}$). 
Similarly, the second equality in (\ref{eq:E_i}) follows by applying 
Lemma~\ref{lem:SL_2} to $f = \Delta_{s_i \omega_i, \omega_i}$ 
(in this case, $k = 1$ and $f' = \Delta_{s_i \omega_i, s_i \omega_i}$). 
This completes the proof of Theorem~\ref{th:minors-Dodgson}.
\endproof

\subsection{Affine coordinates in Schubert cells}
\label{sec:p-coordinates}

For every $w \in W$, the corresponding \emph{Schubert cell}
$(BwB)/B \subset G/B$ is the image of the Bruhat cell $B w B$
under the natural projection of $G$ onto the flag variety~$G/B$.

Let the subgroups $N_+(w) \subset N$ and 
$N_-(w) \subset N_-$ be defined by 
\begin{equation}
\label{eq:N+(u)}
N_+(w) = N \cap {\tilde w} N_- {\tilde w}^{-1}\ ,\quad 
N_-(w) = N_- \cap {\tilde w}^{-1} N {\tilde w} \ ,
\end{equation}
where $\tilde w$ is any representative of $w$ in $G$; 
since $H$ normalizes $N$ and $N_-\,$, these subgroups do not depend on 
the choice of~$\tilde w$. 
The following proposition is essentially well known 
(cf.~\cite[Corollary~23.60]{fulton-harris}).

\begin{proposition} 
\label{pr:Bruhat cell}
An element $x \in G$ lies in the Bruhat cell $B w B$ if and only if,
for some (equivalently, any) representative $\tilde w \in G$ of $w$, 
we have ${\tilde w}^{-1} x \in G_0$ and $[{\tilde w}^{-1} x]_- \in N_- (w)$.
Furthermore, the element
\beal
\label{eq:y+}
y_+ = \pi_+ (x) = \tilde w [{\tilde w}^{-1} x]_- {\tilde w}^{-1} \in N_+ (w)
\eea
does not depend on the choice of $\tilde w$, and the correspondence
$\pi_+ : x \mapsto y_+$ induces a biregular isomorphism between
the Schubert cell $(BwB)/B$ and $N_+ (w)$. 
\end{proposition}

Using the transpose map $x \mapsto x^T$, one obtains a counterpart 
of Proposition~\ref{pr:Bruhat cell} for the opposite Bruhat cell $B_- w B_-$.

\begin{proposition} 
\label{pr:opposite Bruhat cell}
An element $x \in G$ lies in $B_- w B_-$ if and only if,
for some (equivalently, any) representative $\tilde w \in G$ of $w$, 
we have $x {\tilde w}^{-1} \in G_0$ and $[x {\tilde w}^{-1}]_+ \in N_+ (w)$.
Furthermore, the element
\beal
\label{eq:y-}
y_- = \pi_- (x) = {\tilde w}^{-1} [x {\tilde w}^{-1}]_+ {\tilde w} 
\in N_- (w)
\eea
does not depend on the choice of $\tilde w$, and the correspondence
$\pi_- : x \mapsto y_-$ induces a biregular isomorphism between
the ``opposite Schubert cell" $B_- \backslash (B_- w B_-)$ and $N_-(w)$.  
\end{proposition}

The group $N_- (w)$ is a unipotent Lie group of dimension $\l = \l (w)$,
hence it is isomorphic to the affine space $\CC^\l$ as an algebraic variety.
We will associate with any $\ii = (i_1, \ldots, i_\l) \in R(w)$
the following system of affine coordinates on $N_- (w)$.
For $(p_1, \ldots, p_\l) \in \CC^\l$, we set
\beal
\label{eq:N(w)-coordinates}
y_\ii (p_1, \ldots, p_\l) = {\overline {\overline w}}^{-1} \cdot 
\overline {\overline {s_{i_1}}} x_{\overline {i_1}} (p_1) \cdots 
\overline {\overline {s_{i_\l}}} x_{\overline {i_\l}} (p_\l) \ .
\eea
Also, let us define
\begin{equation}
\label{eq:w_k}
w_k = w_{k,\ii} = s_{i_\l} s_{i_{\l-1}} \cdots s_{i_k}
\end{equation}
for $k\in [1,\l+1]$,
so that  $w_1 = w^{-1}$ and $w_{\l+1} = e$. 

\begin{proposition}
\label{pr:N(w)-coordinates}
The map $(p_1, \ldots, p_\l) \mapsto y = y_\ii (p_1, \ldots, p_\l)$
is a biregular isomorphism between $\CC^\l$ and $N_- (w)$.
The inverse map is given by 
\begin{equation}
\label{eq:N(w)-coordinates-det}
p_k = \Delta_{w_k \omega_{i_k}, w_{k+1} \omega_{i_k}} (y) \ .
\end{equation}
\end{proposition}

\proof
We can rewrite (\ref{eq:N(w)-coordinates}) as
$$y_\ii (p_1, \ldots, p_\l) = 
\prod_{k=1}^\l \overline {w_{k+1}} x_{\overline {i_k}} (p_k)
\overline {w_{k+1}}^{\ -1} \ . $$
Each factor 
$\overline {w_{k+1}} x_{\overline {i_k}} (p_k) \overline {w_{k+1}}^{\ -1}$
belongs to the root subgroup in $G$ corresponding to the root
$- w_{k+1} (\alpha_{i_k})$, and these are all the root subgroups in
$N_- (w)$ (cf. \cite[VI, 1.6]{bourbaki}). 
This implies the first statement in Proposition~\ref{pr:N(w)-coordinates}.
To prove (\ref{eq:N(w)-coordinates-det}), we set
$\ii' = (i_1, \ldots, i_{k-1})$ and $\ii'' = (i_{k+1}, \ldots, i_\l)$
so that $\ii = (\ii', i_k, \ii'')$.
Let $y' = y_{\ii'} (p_1, \ldots, p_{k-1})$ and
$y'' = y_{\ii''} (p_{k+1}, \ldots, p_\l)$.
In view of (\ref{eq:N(w)-coordinates}), we have
\beal
\label{eq:three-factors}
\overline {\overline {w_{k}^{\ -1}}} y \overline {w_{k+1}} = 
y' \, (\overline {\overline {s_{i_k}}} x_{\overline {i_k}} (p_k)) \,
(\overline {\overline {w_{k+1}^{\ -1}}} y'' \overline {w_{k+1}}) \ .
\eea
In this decomposition, the first factor $y'$ belongs to
$N_- (s_{i_1} \cdots s_{i_{k-1}}) \subset N_-$, while the last factor
$\overline {\overline {w_{k+1}^{\ -1}}} y'' \overline {w_{k+1}}$ belongs to
$w_{k+1}^{-1} \ N_- (w_{k+1}^{-1}) \ w_{k+1} \subset N$.
Using (\ref{eq:principal-minor-invariance}) and 
(\ref{eq:commutation-s_i}), we conclude that
\bealn
\Delta_{w_k \omega_{i_k}, w_{k+1} \omega_{i_k}} (y) = 
\Delta^{\omega_{i_k}} (\overline {\overline {w_k^{\ -1}}}  y 
\overline {w_{k+1}}) =
\Delta^{\omega_{i_k}} (\overline {\overline {s_{i_k}}} x_{\overline {i_k}}
(p_k)) \\[.1in]
=\Delta^{\omega_{i_k}} (x_{\overline {i_k}} (-p_k^{-1}) p_k^{h_{i_k}} 
x_{i_k} (p_k^{-1})) = p_k \ ,
\eean
as claimed. \endproof

Note for future use that a similar argument allows us to prove that, for
any $i \in [1,r]$, $k\in [1,\l+1]$, and $y \in N_- (w)$, we have
\begin{equation}
\label{eq:trivial-minors-N(w)}
\Delta_{w_k \omega_{i}, w_{k} \omega_{i}} (y) = 1 \ .
\end{equation}
This follows from a decomposition similar to (\ref{eq:three-factors}):
$$\overline {w_{k}}^{\ -1} y \overline {w_{k}} = 
y'\, (\overline {w_{k}}^{\ -1} y'' \overline {w_{k}}) \ ,$$
where $y' \in N_- (s_{i_1} \cdots s_{i_{k-1}}) \subset N_-$, and 
$\overline {w_{k}}^{\ -1} y'' \overline {w_{k}} \in 
w_{k}^{-1} \ N_- (w_{k}^{-1}) \ w_{k} \subset N$.

As a corollary of Proposition~\ref{pr:N(w)-coordinates}, we obtain
defining equations for $N_- (w)$ as a subvariety in $N_-\,$.
Notice that $N_- = N_- (\wnot)$.
Hence, for every $\jj = (j_1, \ldots, j_n) \in R(\wnot)$,
any element $y \in N_-$ can be 
uniquely written as $y = y_\jj (p_1, \ldots, p_n)$
for some $(p_1, \ldots, p_n) \in \CC^n$ (here $n = \l (\wnot)$). 
Let us choose $\jj$ so that its first $n-\l$ indices form a reduced word 
$\jj_1$ for $\wnot w^{-1}$, while the last $\l$ indices 
form a reduced word $\jj_2$ for $w$.
Then write $(\wnot)_k = s_{j_n} s_{j_{n-1}} \cdots s_{j_k}$ for 
$k \in [1,n+1]$, in agreement with~(\ref{eq:w_k}).
Finally, let us denote
$$N'_- (w) = N_- \cap {\tilde w}^{-1} N_- {\tilde w} \ ,$$
where $\tilde w$ is any representative of $w$ in $G$ (cf.~(\ref{eq:N+(u)})).
The following proposition is an immediate consequence of 
Proposition~\ref{pr:N(w)-coordinates} and the
definition~(\ref{eq:N(w)-coordinates}).

\begin{proposition}
\label{pr:N(w)-equations}
Every $y \in N_-$ is uniquely written as $y = y_1 y_2$ with 
$y_1 \in N'_- (w)$ and $y_2 \in N_- (w)$.
In the above notation, if $y = y_\jj (p_1, \ldots, p_n)$ then 
$$y_1 = {\overline {\overline w}}^{-1} y_{\jj_1} (p_1, \ldots, p_{n-\l}) 
\overline {\overline w} \ , \,\, 
y_2 = y_{\jj_2} (p_{n-\l+1}, \ldots, p_n) \ . $$
Hence $y$ lies in $N_-(w)$ if and only if
\begin{equation}
\label{eq:N(w)-equations}
\Delta_{(\wnot)_{k} \omega_{j_k}, (\wnot)_{k+1} \omega_{j_k}} (y) = 0
\end{equation}
for $k = 1, \ldots, n-\l$.
\end{proposition}

\subsection{$y$-coordinates in double Bruhat cells}
\label{sec:y-coordinates}

Let us fix a pair $(u,v) \in W \times W$ and consider 
the open subset $G^{u,v}_0 = G^{u,v} \cap G_0$ 
consisting of the elements $x$ in the double Bruhat cell $G^{u,v}$ 
that admit Gaussian decomposition~(\ref{eq:Gaussian-decomp}). 
In view of Propositions~\ref{pr:Bruhat cell} and 
\ref{pr:opposite Bruhat cell}, the restrictions 
$\pi_+ : G^{u,v} \to N_+ (u)$ and $\pi_- : G^{u,v} \to N_- (v)$
are well defined. 
Let us also introduce the map $\pi_0 : G^{u,v}_0 \to H$ by
\begin{equation}
\label{eq:y0}
y_0 = \pi_0 (x) = [x]_0\ ,
\end{equation}
thus obtaining the map 
\[
\pi^{u,v} = (\pi_+\,, \pi_0\,, \pi_-): G^{u,v}_0 \to 
N_+ (u) \times H \times N_- (v) \ .
\]
For $x \in G^{u,v}_0$, we will write $\pi^{u,v}(x) = (y_+\,, y_0\,, y_-)$ 
and call this triple the $y$-\emph{coordinates} of~$x$. 

\begin{example}
\label{example:sl2-y}
{\rm
Let $G=SL_2(\CC)$, and let $u=v=\wnot$
(cf.\ Example~\ref{example:SL2-cells}).
A matrix
$
x\!=\! \left( 
{\begin{array}{cc}
{{x}_{1 1}} & {{x}_{1 2}} \\
{{x}_{2 1}} & {{x}_{2 2}}
\end{array}}
 \right)
$
with determinant~1 belongs to $G^{u,v}_0$ if and only if
$x_{11}\neq 0$, $x_{12}\neq 0$, $x_{21}\neq 0$.
Using formulas (\ref{eq:y+}), (\ref{eq:y-}), and (\ref{eq:y0}), 
we see that the $y$-coordinates of $x$ are given by
\begin{equation}
\label{eq:gl2yyy}
y_+\!=\!  \left( 
{\begin{array}{rc}
1 & x_{11} x_{21}^{-1} \\ [.1in]
0 & 1
\end{array}}
 \right) 
\ ,\ \  
y_0\!=\!  \left( 
{\begin{array}{cc}
x_{11} & 0 \\[.05in]
0 & x_{11}^{-1}
\end{array}}
 \right) 
\ ,\ \ 
y_-\!=\!  \left( 
{\begin{array}{cr}
1 & 0 \\[.05in]
x_{11}x_{12}^{-1} & 1
\end{array}}
 \right)
\ .
\end{equation}
}
\end{example}

Our use of the term ``coordinates'' for the triple 
$(y_+\,, y_0\,, y_-)$ is justified by the following statement.

\begin{proposition}
\label{pr:y-coordinates}
The map $\pi^{u,v}$ is a biregular isomorphism
$$
G^{u,v}_0 \to (N_+ (u) \cap G_0 u^{-1}) \times H 
\times (N_- (v) \cap v^{-1} G_0)\ .
$$
The inverse isomorphism $(y_+, y_0, y_-) \mapsto x$ is given by
\beal
\label{eq:y-to-x}
[x]_- = [y_+ \tilde u]_-, \,\, [x]_0 = y_0, \,\, 
[x]_+ = [\tilde v y_-]_+ \ ,
\eea
where $\tilde u$ and $\tilde v$ are arbitrary representatives 
of $u$ and~$v$.
\end{proposition} 

\proof
By Proposition~\ref{pr:Bruhat cell}, any $x \in BuB$ can be written as
$x = y_+ \tilde u b$, where $y_+ = \pi_+ (x)$, and $b \in B$. 
It follows that $x \in G_0$ if and only if $y_+ \tilde u \in G_0$,
and if this is the case then $[x]_- = [y_+ \tilde u]_-$.
Similarly, an element $x \in B_- v B_-$ lies in $G_0$ if and only if
$\tilde v y_- \in G_0$, and then $[x]_+ = [\tilde v y_-]_+$.
It follows that $\pi^{u,v}$ is an embedding of $G^{u,v}_0$ 
into $(N_+ (u) \cap G_0 u^{-1}) \times H \times (N_- (v) \cap v^{-1} G_0)$,
and that the inverse map is given by (\ref{eq:y-to-x}).
The same argument shows that if the triple 
$(y_+\,, y_0\,, y_-)$ lies in 
$(N_+ (u) \cap G_0 u^{-1}) \times H \times (N_- (v) \cap v^{-1} G_0)$,
then the element $x$ given by (\ref{eq:y-to-x}) lies in $G^{u,v}_0$,
and we are done. \endproof

The following proposition is immediate from definitions.

\begin{proposition}
\label{pr:T-properties}
We have $N_+ (u)^T = N_- (u^{-1})$ and 
$(G^{u,v}_0)^T = G^{v^{-1},u^{-1}}_0$.
If $x \in  G^{u,v}_0$ has $y$-coordinates $(y_+, y_0, y_-)$ 
then $x^T$ has $y$-coordinates $(y_-^T, y_0, y_+^T)$. 
\end{proposition}

This proposition shows that the transpose map ``interchanges"
the coordinates $y_+$ and $y_-$, so that any statement about
$y_-$ has a counterpart for $y_+$. 
For instance, Proposition~\ref{pr:Bruhat cell} is a counterpart of 
Proposition~\ref{pr:opposite Bruhat cell} in this sense.

\begin{proposition}
\label{pr:diagonal via y}
Suppose $x \in  G^{u,v}_0$ has the $y$-coordinates $(y_+, y_0, y_-)$. 
Then 
\beal
\label{eq:diagonal via y}
[\overline {\overline {u^{-1}}} x]_0^{-1} = [y_+ \overline u]_0 y_0^{-1} 
\ , \,\, [x \overline {v^{-1}}]_0^{-1} =  
y_0^{-1} [\overline{\overline{v}} y_-]_0 \ .
\eea
\end{proposition}

\proof
By Proposition~\ref{pr:Bruhat cell}, $x = y_+ \overline u b$ for some 
$b \in B$.
It follows that $y_0 = [y_+ \overline u]_0 [b]_0$.
On the other hand, 
$\overline {\overline {u^{-1}}} x = {\overline u}^{-1} y_+ \overline u b
\in N_- b$, hence 
$$[\overline {\overline {u^{-1}}} x]_0 = [b]_0 = 
y_0 [y_+ \overline u]_0^{-1} \ .$$
This proves the first equality in (\ref{eq:diagonal via y});
the second one follows by Proposition~\ref{pr:T-properties}.
\endproof

It will be of special importance for us to specialize 
Proposition~\ref{pr:y-coordinates} to the case when
$(u,v) = (e, w)$, where $e$ is the identity element of $W$, and 
$w \in W$ is arbitrary.
Then we have $G^{e,w}_0 = G^{e,w} = H N^w$ where
\begin{equation}
\label{eq:N^w}
N^w = N \cap B_- w B_- \ .
\end{equation}
Specializing Proposition~\ref{pr:y-coordinates} to this case, we obtain
the following statement.

\begin{proposition}
\label{pr:y-for-N}
For any $w \in W$, the map $\pi_-: B_- w B_- \to N_- (w)$
restricts to a biregular isomorphism $N^w \to N_- (w) \cap w^{-1} G_0$.
The inverse isomorphism $N_- (w) \cap w^{-1} G_0 \to N^w$ 
is given by $y \mapsto [\tilde w y]_+$,
where $\tilde w$ is an arbitrary representative of $w$.
\end{proposition}

Using Proposition~\ref{pr:T-properties}, we see that 
Proposition~\ref{pr:y-coordinates} is equivalent to its special case
given by Proposition~\ref{pr:y-for-N} combined with the following
decomposition:
\begin{equation}
\label{eq:Guv0}
G^{u,v}_0  =  (N^{u^{-1}})^T H N^v \ .
\end{equation}

\subsection{Factorization problem in Schubert cells}
\label{sec:ChAnsatz-Nw}

In this section we recall some results from \cite{BZ} 
concerning a version of the factorization problem
for the variety $N^w = N \cap B_- w B_-$ (cf. (\ref{eq:N^w})). 
We will need the following analogue of Theorem~\ref{th:double} which is  
essentially due to G.~Lusztig~\cite{lusztig-reductive} 
(cf. also~\cite[Proposition~1.1]{BZ}).

\begin{proposition}
\label{pr:bireg-N^w}
For any $w \in W$ and any reduced word
$\ii = (i_1, \ldots, i_\l) \in R(w)$, the map 
$(t_1, \ldots, t_\l) \mapsto x_{i_1} (t_1) \cdots x_{i_\l} (t_\l)$ 
is a biregular isomorphism between $\CC_{\neq 0}^\l$ and a Zariski open subset
of $N^{w}$. 
\end{proposition}

We will give explicit formulas for the inverse of the product map in 
Proposition~\ref{pr:bireg-N^w}. 

\begin{theorem}
\label{th:Ansatz-N^w}
Let $\ii = (i_1, \ldots, i_\l) \in R(w)$,
and let $x = x_{i_1} (t_1) \cdots x_{i_\l} (t_\l) \in N^w$
with all $t_k$ nonzero complex numbers.
Then the $t_k$ are recovered from $x$ by
\begin{equation}
\label{eq:Ansatz-N^w}
t_k = 
\displaystyle\frac{1}{\Delta_{w_k \omega_{i_k}, \omega_{i_k}}(y) 
\Delta_{w_{k+1} \omega_{i_k}, \omega_{i_k}}(y)}
\displaystyle\prod_{j \ne i_k} 
\Delta_{w_{k+1} \omega_j, \omega_j}(y)^{-a_{j,i_k}} \ , 
\end{equation}
where $w_k$ is given by {\rm (\ref{eq:w_k})}, and
$y = \pi_- (x) \in N_- (w)$ (cf.~{\rm(\ref{eq:y-})}). 
\end{theorem}

This theorem is a reformulation of \cite[Theorems~1.4, 6.2]{BZ}. 
Here we present a new proof which is in some sense more elementary 
than the one in \cite{BZ}, and also provides additional information
that we will need later on. 

\proof
There is nothing to prove if $w = e$, so we will assume that
$\l (w) = \l \geq 1$. 
Let $y = \pi_- (x)$ and $z = \overline {\overline w} y$.
By Proposition~\ref{pr:y-for-N}, $z \in G_0$, and $x = [z]_+$. 
Let us write $i_1 = i$, and denote $w' = s_i w$, 
$\ii' = (i_2, \ldots, i_\l) \in R(w')$, 
$x' = x_{i} (-t_1) x = x_{i_2} (t_2) \cdots x_{i_\l}(t_\l)\in N^{w'}$,
$y' = \pi_- (x')\in N_- (w')$, and $z' = \overline {\overline {w'}} y'$.
Here is the key lemma.

\begin{lemma}
\label{lem:p-through-t}
In the notation just introduced,
let us write $y' = y_{\ii'} (p_2, \ldots, p_\l)$,
in accordance with Proposition~\ref{pr:N(w)-coordinates}.
Then $y = y_{\ii} (p_1, p_2, \ldots, p_\l)$,
where $p_1$ is given by
\beal
\label{eq:p1}
p_1 = \Delta_{s_i \omega_i, \omega_i} 
(x_{\overline i} ([z']_0^{- \alpha_i} t_1^{-1}) [z']_-^{-1}) \ .
\eea
Furthermore, we have
\beal
\label{eq:t1}
t_1 = [z']_0^{\omega_i - \alpha_i} [z]_0^{- \omega_i} \ .
\eea

\end{lemma}

\proof
Let us temporarily denote $\tilde y = y_{\ii} (p_1, p_2, \ldots, p_\l)$,
and $\tilde z = \overline {\overline w} \tilde y$, where $p_1$ is given by
(\ref{eq:p1}); our goal is to show that $\tilde y = y$ and $\tilde z = z$.
By Proposition~\ref{pr:y-for-N}, it suffices to show that $[\tilde z]_+ = x$, 
or equivalently that $\tilde z x^{-1} \in B_-$. 

By Proposition~\ref{pr:N(w)-equations} (applied to $w = s_i$),
formula (\ref{eq:p1}) implies that
\beal
\label{eq:y''}
x_{\overline i} ([z']_0^{- \alpha_i} t_1^{-1}) [z']_-^{-1} = 
\overline {s_i}  y'' {\overline {s_i}}^{\ -1} x_{\overline i} (p_1) \ ,
\eea
where $y'' \in N_-$. 
Using (\ref{eq:H-conjugation}) and (\ref{eq:5-commutation-special}),
we can rewrite the left-hand side of (\ref{eq:y''}) as follows:
$$x_{\overline i} ([z']_0^{- \alpha_i} t_1^{-1}) [z']_-^{-1} = 
[z']_0 x_{\overline i} (t_1^{-1}) [z']_0^{-1} [z']_-^{-1} = 
[z']_0 x_{\overline i} (t_1^{-1}) x' (z')^{-1} = $$
$$[z']_0 x_{\overline i} (t_1^{-1}) x_i (-t_1) x (z')^{-1} =
[z']_0 \overline {s_i} t_1^{-h_i} x_{\overline i} (- t_1^{-1}) x (z')^{-1} 
\ .$$
Substituting this expression into (\ref{eq:y''}) and using the fact
that $\tilde z = {\overline {s_i}}^{-1} x_{\overline i} (p_1) z'$,
we can rewrite (\ref{eq:y''}) as follows:
\beal
\label{eq:z xinv}
\tilde z x^{-1} = 
(y'')^{-1} {\overline {s_i}}^{-1} [z']_0 \overline {s_i} 
t_1^{-h_i} x_{\overline i} (- t_1^{-1}) \ .
\eea
It follows that $\tilde z x^{-1} \in B_-$, hence 
$\tilde y = y$ and $\tilde z = z$.
Since the left-hand side of (\ref{eq:z xinv}) is equal to
$[z]_-[z]_0\,$, it follows that 
\begin{equation}
\label{eq:z-diagonal}
[z]_0 = {\overline {s_i}}^{\ -1} [z']_0 \overline {s_i} \
t_1^{-h_i} \ .
\end{equation}
Finally, (\ref{eq:t1}) follows from (\ref{eq:z-diagonal}) by applying
the character $a \mapsto a^{\omega_i}$ to both sides and using
(\ref{eq:W-action}). 
\endproof

Note that (\ref{eq:p1}) can be simplified as follows:
\begin{equation}
\label{eq:p_1-alt}
p_1 =  \Delta_{\omega_i, {w'}^{-1} \omega_i} (x {y'}^{\ -1})/ 
\Delta_{\omega_i, {w}^{-1} \omega_i} (x) \ ;
\end{equation}
since we will not need this formula, the proof is left to the reader.

Continuing with the proof of Theorem~\ref{th:Ansatz-N^w},
let us define, for $k = 1,\dots, \l$:
$$x^{(k)} = x_{i_k} (t_k) \cdots x_{i_\l} (t_\l) \in N^{w_k^{-1}} \ ,
\,\, y^{(k)} = \pi_- (x^{(k)}) \ , 
\,\, z^{(k)} = \overline {\overline{w_k^{-1}}} y^{(k)} \ .$$
Applying (\ref{eq:t1}) with $x$ replaced by $x^{(k)}$ yields
\beal
\label{eq:t-through-z}
t_k = [z^{(k+1)}]_0^{\omega_{i_k} - \alpha_{i_k}} 
[z^{(k)}]_0^{- \omega_{i_k}} \ .
\eea
On the other hand, combining the definition (\ref{eq:Delta-general})
with (\ref{eq:weight}), we can rewrite (\ref{eq:Ansatz-N^w}) as follows:
\beal
\label{eq:Ansatz-N^w-diagonal}
t_k = 
[\overline {\overline{w_{k+1}^{-1}}} y]_0^{\omega_{i_k} - \alpha_{i_k}} 
[\overline {\overline{w_k^{-1}}} y]_0^{- \omega_{i_k}}  
\eea
Comparing (\ref{eq:Ansatz-N^w-diagonal}) with 
(\ref{eq:t-through-z}), we see that Theorem~\ref{th:Ansatz-N^w}
would follow from the equality 
$[\overline {\overline{w_k^{-1}}} y]_0 = [z^{(k)}]_0$.
The latter is obtained by observing that
$\overline {\overline{w_k^{-1}}} y =
\overline {\overline{w_k^{-1}}} {\overline {\overline w}}^{-1} z
= \tilde y z^{(k)}$, where
$\tilde y=y_{(i_1,\dots,i_{k-1})}(p_1,\dots,p_{k-1})\in N_-$
(this $\tilde y$ was denoted by $y'$ in (\ref{eq:three-factors})).
\endproof

\subsection{Totally positive bases for $N_-(w)$}

Although most of the results in this section were obtained in
\cite{BZ}, we prefer to give independent proofs here; 
in some cases, this will allow us to refine the statements
in~\cite{BZ}. 

We start with the following general definition.

\begin{definition}
\label{def:TP-base} {\rm 
Let $F$ be a finite collection of functions on a set $X$. 
A subset $\B \subset F$ is called a 
\emph{totally positive base} for $F$ if $\B$ is a minimal 
(with respect to inclusion) subset of $F$ with the property
that every $f \in F$ is a \emph{subtraction-free expression} 
(i.e., a ratio of two polynomials with nonnegative integer
coefficients) in the elements of~$\B$. 
}\end{definition} 


For every $w \in W$, let us denote 
\beal
\label{eq:F(w)}
F(w) = \{ \Delta_{w'' \omega_i, w' \omega_i} (y)
\,:\,
i \in [1,r]\,,\ w' \preceq w'' \preceq w^{-1}\}\ .
\eea
(As earlier in (\ref{eq:F(u,v)}), 
$w' \preceq w''$ stands for $\l(w'')=\l(w')+\l({w'}^{-1} w'')$.)
To every reduced word $\,\ii = (i_1, \ldots, i_m) \in R(w)$
we associate three collections of regular functions on the group $N_- (w)$. 
\begin{eqnarray}
\label{eq:FFF}
\begin{array}{rcl}
F_1 (\ii) &=& \{\Delta_{w_k \omega_{i_k}, \omega_{i_k}} \,:\, 
1 \leq k \leq m \} \ ,\\[.1in]
F_2 (\ii) &=& \{\Delta_{w^{-1} \omega_{i_k}, w_k \omega_{i_k}} \,:\,
2 \leq k \leq m+1 \} \ ,\\[.1in]
F (\ii) &=& 
\{\Delta_{w_k \omega_{i}, w_l \omega_{i}}  \,:\, 
i \in [1,r], \, 1 \leq k \leq l \leq m+1\} \ ,
\end{array}
\end{eqnarray}
where $w_k=s_{i_m} \cdots s_{i_k}$ (cf.\ (\ref{eq:w_k})).

\begin{theorem}
\label{th:TP-base-y}
For any reduced word $\,\ii = (i_1, \ldots, i_m) \in R(w)$,
each of the collections $F_1 (\ii)$ and  $F_2 (\ii)$ of regular functions on
$N_- (w)$ is a transcendence basis for $\CC(N_- (w))$ and a 
totally positive base for~$F(w)$.
\end{theorem}

\proof
Let us first deal with $F_1 (\ii)$. 
The most important part of the proof is to show that every minor in $F(w)$
is a subtraction-free expression in the minors from $F_1(\ii)$. 
Since we obviously have 
\[
F(w) = \bigcup_{\ii \in R(w)} F(\ii)\ ,
\] 
this statement will directly follow from
Lemmas~\ref{lem:TP-base-BZ} and \ref{lem:TP-base-ii} below. 

\begin{lemma}
\label{lem:TP-base-BZ}
For any two reduced words $\ii, \ii' \in R(w)$, every minor in $F_1(\ii')$
is a subtraction-free expression in the minors from $F_1(\ii)$. 
\end{lemma}

\proof
This is an immediate corollary of \cite[Corollary~6.7]{BZ}.
The proof in \cite{BZ} is based on repeated
applications of determinantal identities of
Theorem~\ref{th:minors-Plucker}. 
\endproof

\begin{lemma}
\label{lem:TP-base-ii}
Every minor in $F(\ii)$
is a subtraction-free expression in the minors from $F_1(\ii)$. 
\end{lemma}

\proof
We need to show that every minor 
$\Delta_{w_k \omega_{i}, w_l \omega_{i}}$ for $i \in [1,r]$ and 
$1 \leq k \leq l \leq m+1$ is a subtraction-free expression
in the minors $\Delta_{w_k \omega_{i}, \omega_{i}}$. 
Recall that, by convention, $w_{m+1} = e$, so the  
statement is trivial for $l = m+1$.
By (\ref{eq:trivial-minors-N(w)}), it also holds 
for $k = l$, since the corresponding minor equals~$1$. 
Thus we may assume that $1 \leq k < l \leq m$;
increasing $k$ and $l$ if necessary, we can also assume without loss
of generality that $i_k = i_l = i$.

Let us arrange all the pairs $(k,l)$ with $1 \leq k \leq l \leq m+1$ 
in the following order: $(k',l')<(k,l)$ if 
either $l' > l$, or $l' = l, \, k' > k$. 
Using induction with respect to this linear order, it is enough to show 
that for every $(k,l)$ such that $1 \leq k < l \leq m$ and $i_k = i_l = i$,
the minor $\Delta_{w_k \omega_{i}, w_l \omega_{i}}$     is a 
subtraction-free expression in the minors 
$\Delta_{w_{k'} \omega_{j}, w_{l'} \omega_{j}}$ with $j \in [1,r]$
and $(k',l')<(k,l)$. 
The latter follows from the identity (\ref{eq:minors-Dodgson})
applied to $u = w_{k+1}$ and $v = w_{l+1}\,$.
Indeed, this identity
can be rewritten as
$$\Delta_{w_k \omega_{i}, w_l \omega_{i}} = \frac
{\Delta_{w_k \omega_{i}, w_{l+1} \omega_{i}} 
\Delta_{w_{k+1} \omega_{i}, w_l \omega_{i}} + \prod_{j \neq i}
\Delta_{w_{k+1} \omega_{j}, w_{l+1} \omega_{j}}^{- a_{ji}}}
{\Delta_{w_{k+1} \omega_{i}, w_{l+1} \omega_{i}}} \ ,$$
providing a desired subtraction-free expression. 
\endproof

Lemma~\ref{lem:TP-base-ii} implies in particular that each minor
$\Delta_{w_k \omega_{i_k}, w_{k+1} \omega_{i_k}} (y)$ is a rational 
function of the minors from $F_1 (\ii)$. 
By Proposition~\ref{pr:N(w)-coordinates}, it follows that 
$F_1 (\ii)$ is a transcendence basis for $\CC(N_- (w))$,
hence it is a totally positive base for $F(w)$. 

To prove that $F_2(\ii)$ has the same properties, we will apply 
the anti-automorphism $\tau_w$ of $G$ given by 
\begin{equation}
\label{eq:tau-w}
\tau_w (y) = \overline w (y^{-1})^{\theta} {\overline w}^{\ -1} \ ,
\end{equation} 
where $\theta$ was defined in~(\ref{eq:theta}). 
In view of (\ref{eq:barT}), if $\tau_w (y) = y'$, then
\begin{equation}
\label{eq:tau-T}
y = \tau_{w^{-1}} (y') \ , \ \ 
 \tau_w (y^T) = {y'}^T \ .
\end{equation}
A straightforward check shows that 
\beal
\label{eq:tau-property1}
\tau_{w^{-1}} (N_+ (w)) = N_+ (w^{-1})\ ,\ \ 
\tau_w (N_- (w)) = N_- (w^{-1})  \ .
\eea

\begin{lemma}
\label{lem:minors-tau}
Let $y \in G$, and let $y' =  \tau_{w} (y)$ for some $w \in W$.
For any $w', w'' \in W$ such that $w' \preceq w^{-1}$ and 
$w'' \preceq w^{-1}$, we have
\beal
\label{eq:minors-tau}
\Delta_{w'' \omega_i, w' \omega_i} (y) = 
\Delta_{w w' \omega_i, w w'' \omega_i} (y') \ .
\eea
\end{lemma}

\proof
In view of (\ref{eq:braid}) and (\ref{eq:barT}), the conditions 
$w' \preceq w^{-1}$ and $w'' \preceq w^{-1}$ imply that 
$$
\overline {\overline {(w w')^{\ -1}}} = 
{\overline {\overline {w'}}}^{\ -1} \ \overline {w}^{\ -1} \ ,
\,\, \overline {w w''} = \overline {w} \ \overline {\overline {w''}} \ .
$$
Combining this with (\ref{eq:principal-minor-T}), we obtain: 
\beal
\Delta_{w w' \omega_i, w w'' \omega_i} (y') =
\Delta^{\omega_i} (\overline {\overline {(w w')^{\ -1}}} y' 
\overline {w w''}) = 
\Delta^{\omega_i} ({\overline {\overline {w'}}}^{\ -1} \ 
{\overline w}^{\ -1} y' \ \overline {w} \ \overline {\overline {w''}}) 
\\[.1in]
= \Delta^{\omega_i} ({\overline {\overline {w'}}}^{\ -1} (y^{-1})^\theta 
\ \overline {\overline {w''}}) = 
\Delta^{\omega_i} (\overline {\overline {{w''}^{\ -1}}} y \overline {w'}) =
\Delta_{w'' \omega_i, w' \omega_i} (y) \ ,
\eea
as claimed. 
\endproof 

By Lemma~\ref{lem:minors-tau}, the antiautomorphism $\tau_w$
transforms $F(w)$ into $F(w^{-1})$, and $F_2 (\ii)$ into $F_1 (\ii^*)$,
where $\ii^*=(i_m, \ldots, i_1) \in R(w^{-1})$ is $\ii$ written backwards. 
Thus the fact that $F_2 (\ii)$ is a transcendence basis for $\CC(N_- (w))$
and a totally positive base for $F(w)$, follows from the same properties
for $F_1 (\ii)$ that we already proved. 
This completes the proof of Theorem~\ref{th:TP-base-y}. \endproof


\subsection{Total positivity in $y$-coordinates}

Let $N_{\geq 0} \subset N$ denote the multiplicative semigroup
generated by the elements $x_i (t)$ for $i\in [1,r]$ and $t > 0$.
For every $w \in W$, let us denote (cf.~\cite{BZ})
\beal
\label{eq:N^w>0}
N^w_{> 0} = N_{\geq 0} \cap N^w =  N_{\geq 0} \cap B_-w B_- \ .
\eea
The following analogue of Theorem~\ref{th:LusztigTP} is due to 
G.~Lusztig~\cite{lusztig-reductive} (cf. Proposition~\ref{pr:bireg-N^w}).

\begin{proposition}
\label{pr:TP-N^w}
For any $w \in W$ and any reduced word
$\ii = (i_1, \ldots, i_\l) \in R(w)$, the map 
$(t_1, \ldots, t_\l) \mapsto x_{i_1} (t_1) \cdots x_{i_\l} (t_\l)$ 
restricts to a bijection $\RR_{> 0}^\l \to N^{w}_{> 0}$.
\end{proposition}

We will use Theorem~\ref{th:TP-base-y} to obtain 
the following criteria for total positivity.

\begin{theorem}
\label{th:TP-criteria-y}
Let $x \in N^w$, let $y = \pi_- (x) \in N_- (w)$, and let 
$\ii = (i_1, \ldots, i_\l) \in R(w)$.
Then the following conditions are equivalent:
\smallskip

\noindent {\rm(1)} $x \in N^w_{>0}\,$;
\smallskip

\noindent {\rm(2)} $\Delta (y) > 0$ for any $\Delta \in F(w)$;
\smallskip

\noindent {\rm(3)} $\Delta (y) > 0$ for any $\Delta \in F_1 (\ii)$;
\smallskip

\noindent {\rm(4)} $\Delta (y) > 0$ for any $\Delta \in F_2 (\ii)$.
\end{theorem}

\proof
The equivalence of (2), (3) and (4) is immediate from 
Theorem~\ref{th:TP-base-y}. 
Let us show the equivalence of (1) and (3).
By Proposition~\ref{pr:TP-N^w}, every 
$x \in N^w_{>0}$ is of the form 
$x = x_{i_1} (t_1) \cdots x_{i_\l} (t_\l))$ for some 
$t_1, \ldots, t_\l > 0$. 
By Theorem~\ref{th:Ansatz-N^w}, each $t_k$ is a monomial
in $\l$ variables $\{\Delta(y): \Delta \in F_1 (\ii)\}$.
It follows that the monomial transformation from 
$\{\Delta(y): \Delta \in F_1 (\ii)\}$ to $\{t_1, \ldots, t_\l\}$ 
is invertible (an explicit expression for the inverse 
transformation was given in \cite[Theorem~4.3]{BZ} but we will not need
it here). 
Thus every $\Delta (y)$ with $\Delta \in F_1 (\ii)$ is a 
Laurent monomial in $t_1, \ldots, t_\l$. 
Hence $\Delta (y) > 0$, and $(1) \Rightarrow (3)$ is proved. 

To prove $(3) \Rightarrow (1)$, suppose that
$\Delta (y) > 0$ for $\Delta \in F_1 (\ii)$. 
Let us define $t_1, \ldots, t_\l$ via (\ref{eq:Ansatz-N^w}),
and let $\tilde x = x_{i_1} (t_1) \cdots x_{i_\l} (t_\l) \in N^w_{>0}$.
Setting ${\tilde y} = \pi_- (\tilde x)$,
we see that $\Delta (y) = \Delta ({\tilde y})$ 
for any $\Delta \in F_1 (\ii)$. 
By Lemma~\ref{lem:TP-base-ii}, we have $\Delta (y) = \Delta ({\tilde y})$ 
for any $\Delta \in F (\ii)$. 
In particular, 
$\Delta_{w_k \omega_{i_k}, w_{k+1} \omega_{i_k}} (y)= 
\Delta_{w_k \omega_{i_k}, w_{k+1} \omega_{i_k}} ({\tilde y})$ for 
$k = 1, \ldots, \l$.
Using (\ref{eq:N(w)-coordinates-det}), we conclude that
$y = {\tilde y}$ and so $x = \tilde x \in N^w_{>0}\,$,
proving $(3) \Rightarrow (1)$.
\endproof

\begin{corollary}
\label{cor:tau-positivity}
The map $\tau_w:  N_- (w) \to N_- (w^{-1})$ 
(cf.\ {\rm(\ref{eq:tau-property1})}) restricts to a bijection 
$\pi_- (N^w_{> 0}) \to \pi_- (N^{w^{-1}}_{> 0})$.
\end{corollary} 

\proof
We have already observed that by Lemma~\ref{lem:minors-tau},
$\tau_w$ transforms $F(w)$ into $F(w^{-1})$.
The corollary then follows from Theorem~\ref{th:TP-criteria-y}.
\endproof


We will now show that using $y$-coordinates 
(i.e., passing from a double Bruhat cell $G^{u,v}$ to 
the open subset $G^{u,v}_0$)
will not create problems in the study of totally positive varieties.

\begin{proposition}
\label{pr:TP-Gauss}
We have 
$N_{\geq 0} = N \cap G_{\geq 0}$, $N^T_{\geq 0} = N_- \cap G_{\geq 0}$ and 
$G_{\geq 0} = N_{\geq 0}^T H_{> 0} N_{\geq 0}$.
In particular, $G_{\geq 0} \subset G_0\,$, i.e.,  
any totally nonnegative element in $G$ admits the Gaussian decomposition.
Furthermore, for any $u, v \in W$, the totally positive variety
$G^{u,v}_{>0}$ decomposes as
\beal
\label{eq:Guv>0}
G^{u,v}_{>0}  =  (N^{u^{-1}}_{>0})^T H_{>0} N^v_{>0} \ ,
\eea
in the notation of~{\rm (\ref{eq:N^w>0})}. 
\end{proposition}

\proof
By the definition of $G_{\geq 0}$, every totally nonnegative element 
$x \in G$ has the form (cf. (\ref{eq:productmap}))
$x = x_\ii (a; t_1, \ldots, t_m)$, where $\ii = (i_1, \ldots, i_m)$
is some word in the alphabet $[1,r] \cup [\overline 1, \overline r]$,
the $t_k$ are positive real numbers, and $a \in H_{>0}$.
We say that $\ii$ is \emph{unmixed} if all the indices from 
$[\overline 1, \overline r]$ precede those from $[1, r]$. 
By repeated application of
the commutation relations (\ref{eq:H-conjugation}), 
(\ref{eq:2-commutation-mixed}) and (\ref{eq:5-commutation-generic}),
we can transform $x$ to the form $x = x_{\ii'} (a'; t'_1, \ldots, t'_m)$
for an unmixed word $\ii'$, $a' \in H_{>0}$, and all $t'_k > 0$. 
This proves the decomposition $G_{\geq 0} = N_{\geq 0}^T H_{> 0} N_{\geq 0}$.
The equalities $N_{\geq 0} = N \cap G_{\geq 0}$ and  
$N_{\geq 0}^T = N_- \cap G_{\geq 0}$ follow from this decomposition
of $G_{\geq 0}$ and the uniqueness of the Gaussian decomposition.
Finally, (\ref{eq:Guv>0}) is proved by the same argument combined
with (\ref{eq:Guv0}). \endproof

Combining Propositions~\ref{pr:TP-Gauss} and \ref{pr:T-properties}
with Theorem~\ref{th:TP-criteria-y}, we obtain the following description
of the totally positive variety $G^{u,v}_{>0}$ in terms of $y$-coordinates.

\begin{theorem}
\label{th:TP-through-y+-0}
An element $x \in G^{u,v}_0$ lies in $G^{u,v}_{>0}$ if and only if 
its $y$-coordinates $(y_+, y_0, y_-)$ satisfy the following properties:
\begin{itemize}
\item $\Delta (y_-) > 0$ for any $\Delta \in F(v)$;
\item $\Delta (y_+^T) > 0$ for any $\Delta \in F(u^{-1})$;
\item $y_0 \in H_{>0}$.
\end{itemize}
\end{theorem}


\section{Proofs of main results}
\label{sec:proofs-general}

This section contains proofs of the main results in Section~\ref{sec:main}.

\subsection{Proofs of Theorems~\ref{th:Guv-affine}, 
\ref{th:double}, and \ref{th:LusztigTP}}

\emph{Proof of Theorem~\ref{th:Guv-affine}.} 
We will explicitly construct a desired biregular isomorphism of $G^{u,v}$
with a Zariski open subset of $\CC^{r+\l(u)+\l(v)}$
with the help of a ``twisted" version of $y$-coordinates (cf. 
Section~\ref{sec:y-coordinates}).
We fix a representative $\tilde u$ of $u$, and associate to any 
$x \in G^{u,v}$ a triple $(y_{(+)}, y_{(0)}, y_{(-)})$ given by
\beal
\label{eq:y-for-Guv}
y_{(+)} = \pi_+ (x), \,\,  y_{(-)} = \pi_- (x^{-1}), \,\,
y_{(0)} = [\tilde u^{-1} x]_0 \ ;
\eea
in view of Propositions~\ref{pr:Bruhat cell} and 
\ref{pr:opposite Bruhat cell}, this triple is well defined and belongs to 
$N_+ (u) \times H \times N_- (v^{-1})$. 
Our statement is a consequence of the following.

\begin{proposition}
\label{pr:Guv-embedding}
The correspondence $x \mapsto (y_{(+)}, y_{(0)}, y_{(-)})$ given by 
{\rm(\ref{eq:y-for-Guv})} is a biregular isomorphism of $G^{u,v}$
with the Zariski open subset of $N_+ (u) \times H \times N_- (v^{-1})$
consisting of triples such that $y_{(-)} y_{(+)} \in v G_0 u^{-1}$. 
\end{proposition}

\proof
The proof is essentially the same as that of 
Proposition~\ref{pr:y-coordinates};
the inverse of the correspondence (\ref{eq:y-for-Guv})
is given by
$$x = y_{(+)} \tilde u [\tilde v^{-1} y_{(-)} y_{(+)} \tilde u]_+^{-1}  
y_{(0)} \ ,$$
where $\tilde v$ is any representative of $v$. \endproof


\emph{Proof of Theorem~\ref{th:double}.}
Actually, we will only prove at this point the following weak
version of this theorem. 

\begin{proposition}
\label{product-embedding}
For every $u,v\in W$ and $\ii= (i_1, \ldots, i_m)\in R(u,v)$, 
the map $x_\ii$ restricts to an injective regular map
$H \times\CC_{\neq 0}^m \to G^{u,v}$.
\end{proposition}

Before proving this proposition, let us make some comments.
In view of Theorem~\ref{th:Guv-affine}, we can think of $G^{u,v}$
as a Zariski open subset of $\CC^{m+r}$. 
On the other hand, $H \times\CC_{\neq 0}^m$ can be identified 
with the complex torus $\CC_{\neq 0}^{m+r}$
(for instance, by using the coordinates $t_1, \ldots, t_m$ and 
$a^{\omega_1}, \ldots, a^{\omega_r}$). 
Thus, Proposition~\ref{product-embedding} allows us to
think of the map $x_\ii: H \times\CC_{\neq 0}^m \to G^{u,v}$
as of a regular embedding $\varphi: \CC_{\neq 0}^{m+r} \to \CC^{m+r}$.
It is easy to see that such an embedding $\varphi$ is always a birational
isomorphism.
Theorem~\ref{th:double} makes a stronger claim that in our situation
the image of $\varphi$ is Zariski open, and $\varphi$ is a biregular
isomorphism of $\CC_{\neq 0}^{m+r}$ onto its image.
This will follow from the explicit form of the inverse map $x_\ii^{-1}$
given by Theorem~\ref{th:t-through-x} 
(to be proven in Section~\ref{sec:t-through-x}).

\proof
First let us show that $x_\ii (H \times\CC_{\neq 0}^m) \subset G^{u,v}$.
We will show that $x_\ii (H \times\CC_{\neq 0}^m) \subset B_- v B_-$;
the inclusion $x_\ii (H \times\CC_{\neq 0}^m) \subset B u B$
is proved similarly (or deduced from the previous one with the help
of the transpose map).
Let 
$$\{k_1 < \cdots < k_\l\} = 
\{k \in [1,m]: i_k \in [1, r]\} \ ,$$ 
so that $(i_{k_1}, \ldots, i_{k_\l}) \in R(v)$.  
Let us use the fact that, for every $i \in [1,r]$ and nonzero $t \in \CC$,
we have $x_i (t) \in B_- s_i B_-$ and $x_{\overline i} (t) \in B_-$
(cf.~(\ref{eq:5-commutation-special})).
It follows that if $(a; t_1, \ldots, t_m) \in H \times\CC_{\neq 0}^m$ then
$$x_\ii (a; t_1, \ldots, t_m)  \in 
B_- \cdot B_- s_{i_{k_1}} B_-  \cdot B_- 
\cdots B_- \cdot B_- s_{i_{k_\l}} B_- \cdot B_- 
= B_- v B_- \ ,$$
as desired; the last equality follows from the well known fact that 
$$B_- w' B_- \cdot B_- w'' B_- = B_- w'w'' B_-$$ 
whenever $\l (w' w'') = \l (w') + \l (w'')$ (cf.~\cite[IV.2.4]{bourbaki}).

It remains to show that the map $x_\ii: H \times\CC_{\neq 0}^m \to G^{u,v}$
is injective. 
There is nothing to prove if $u = v = e$, so we can assume that
$m = \l(u) + \l(v) \geq 1$. 
Suppose that $i_m = i \in [1,r]$ (the case when 
$i_m \in [\overline 1, \overline r]$ is treated in the same way). 
Denote $v' = v s_i$ so that $\ii' = (i_1, \ldots, i_{m-1}) \in R(u,v')$.
Now suppose
$$x_\ii (a; t_1, \ldots, t_m) = x_\ii (a'; t'_1, \ldots, t'_m) \ ,$$
where $(a; t_1, \ldots, t_m)$ and $(a'; t'_1, \ldots, t'_m)$ belong to
$H \times\CC_{\neq 0}^m$. 
Multiplying both sides of the last equality on the right by 
$x_i (- t'_m)$, we obtain: 
\beal
\label{eq:i=i'}
x_\ii (a; t_1, \ldots, t_{m-1}, t_m - t'_m) = 
x_{\ii'} (a'; t'_1, \ldots, t'_{m-1}) \ .
\eea
If $t_m \neq t'_m\,$, then the left-hand side of (\ref{eq:i=i'}) lies in
$G^{u,v}$, while the right-hand side lies in $G^{u,v'}$.
But this is impossible because the double Bruhat cells are disjoint.
Thus, $t_m = t'_m$, and the desired injectivity follows by induction on $m$.
\endproof

\emph{Proof of Theorem~\ref{th:LusztigTP}.}
If a double reduced word $\ii \in R(u,v)$ is unmixed,
i.e., all the indices from $[\overline 1, \overline r]$ precede those 
from $[1, r]$, then our statement follows by combining 
(\ref{eq:Guv>0}) with Proposition~\ref{pr:TP-N^w}.
The statement for an arbitrary $\ii \in R(u,v)$ can be reduced to the case
of an unmixed $\ii$ by the argument used in the proof of
Proposition~\ref{pr:TP-Gauss}, i.e, by repeated application of
the commutation relations (\ref{eq:H-conjugation}), 
(\ref{eq:2-commutation-mixed}) and (\ref{eq:5-commutation-generic}).
\endproof

\subsection{Proofs of Theorems~\ref{th:zeta-regularity} and
\ref{th:zeta-positivity}} 

\emph{Proof of Theorem~\ref{th:zeta-regularity}.}
The fact that the right-hand side of~(\ref{eq:zeta-u,v-x}) 
is well defined for any $x\in G^{u,v}$ follows from 
Propositions~\ref{pr:Bruhat cell} and \ref{pr:opposite Bruhat cell}.
Let us show that $x' = \zeta^{u,v} (x)\in G^{u^{-1}, v^{-1}}$.
Using (\ref{eq:y+}), we can rewrite $x'$ as
$$x' = {\overline{\overline u}}^{-1} \, 
\left(
y_+^{-1}\, [x\overline{v^{-1}}]_- [x\overline{v^{-1}}]_0
\right)^\theta \ ,$$
where $y_+ = \pi_+ (x)$. 
It follows that ${\overline{\overline u}} x' \in G_0$, and
\beal 
\label{eq:y+ for x'}
[\overline{\overline u} x']_- = (y_+^{-1})^\theta \ .
\eea
Hence $[\overline{\overline u} x']_- \in N_+ (u)^\theta = N_- (u^{-1})$,
and we conclude from Proposition~ \ref{pr:Bruhat cell} that 
$x' \in B u^{-1} B$. 
The inclusion $x' \in B_- v^{-1} B_-$ is proved in a similar way
(or by using the transpose map); the counterpart of (\ref{eq:y+ for x'})
is given by
\beal 
\label{eq:y- for x'}
[x' \overline v]_+ = (y_-^{-1})^\theta \ ,
\eea
where $y_- = \pi_- (x)$ (cf. (\ref{eq:y-})).

We have proved that $x' \in G^{u^{-1}, v^{-1}}$.
To complete the proof of Theorem~\ref{th:zeta-regularity}, 
it suffices to show that 
$\zeta^{u^{-1}, v^{-1}} ( \zeta^{u,v} (x)) = x$ for any $x \in G^{u,v}$.
Notice that (\ref{eq:y+ for x'}) and (\ref{eq:y- for x'})
can be rewritten as
\beal
\label{eq:x'+- symmetric}
[\overline{\overline u} x']_- = \overline{\overline u} \, 
\left([\overline{\overline {u^{-1}}} x]_-^{-1}\right)^\theta \,
{\overline{\overline u}}^{-1} \ , \,\,
[x'\overline v]_+ = {\overline v}^{-1} \, 
\left([x \overline {v^{-1}}]_+^{-1}\right)^\theta \,
\overline v \ .
\eea
The desired equality $\zeta^{u^{-1}, v^{-1}}(x') = x$
follows by subsituting these expressions into the expression
for $\zeta^{u^{-1}, v^{-1}}(x')$ obtained from~(\ref{eq:zeta-u,v-x}). 
\endproof

The following proposition shows that the twist map 
respects the Gaussian decomposition.

\begin{proposition}
\label{pr:twist-Gauss}
The twist map $\zeta^{u,v}: x \mapsto x'$ sends the open subset
$G^{u,v}_0$ to $G^{u^{-1}, v^{-1}}_0$, and we have
\beal
\label{eq:x'-diagonal}
[x']_0 
= [\overline {\overline {u^{-1}}} x]_0^{-1} 
[x]_0 [x \overline {v^{-1}}]_0^{-1} \ .
\eea
\end{proposition}

\proof
To show that $x$ and $x'$ belong to $G_0$ simultaneously, 
it suffices to rewrite (\ref{eq:zeta-u,v-x}) as 
$$x' = \left(
[\overline{\overline{u^{-1}}}x]_0  [\overline{\overline{u^{-1}}}x]_+ \,
x^{-1} \, [x \overline{v^{-1}}]_- [x \overline{v^{-1}}]_0
\right)^\theta \ .$$
In view of (\ref{eq:T-Gauss}), this also implies 
(\ref{eq:x'-diagonal}). \endproof

Let us now describe the twist map in terms of $y$-coordinates.
Recall the definition (\ref{eq:tau-w}) of
the anti-automorphism $\tau_w$ of the group~$G$.

\begin{proposition}
\label{pr:twist-y}
Suppose $x \in G^{u,v}_0$ has $y$-coordinates $(y_+, y_0, y_-)$.
Then the $y$-coordinates $(y'_+, y'_0, y'_-)$ of 
$x'=\zeta^{u,v} (x)$ are given by 
\beal
\label{eq:zeta-u,v}
y'_+ = \tau_{u^{-1}} (y_+)\ ,\\[.1in]
y'_0 = [y_+ \overline u]_0 \, y_0^{-1} \, 
[\overline{\overline{v}} y_-]_0\ ,\\[.1in]
y'_- = \tau_{v} (y_-) \ . 
\eea
\end{proposition}

\proof
The desired expressions for $y'_+$ and $y'_-$ follow from
(\ref{eq:x'+- symmetric}); the expression for $y'_0$ follows from
(\ref{eq:x'-diagonal}) combined with (\ref{eq:diagonal via y}).
\endproof

\emph{Proof of Theorem~\ref{th:zeta-positivity}.}
Let $x \in G^{u,v}_{> 0}$, and let $x'=\zeta^{u,v} (x)$. 
By Proposition~\ref{pr:TP-Gauss}, $x \in G^{u,v}_0$ so $x$ 
has well-defined $y$-coordinates $(y_+, y_0, y_-)$.
By Proposition~\ref{pr:twist-y}, the $y$-coordinates 
$(y'_+, y'_0, y'_-)$ of $x'$ are given by (\ref{eq:zeta-u,v}). 
By Theorem~\ref{th:TP-through-y+-0}, the triple $(y_+, y_0, y_-)$
satisfies the properties given there, and 
it suffices to check that $(y'_+, y'_0, y'_-)$ satisfies
the same properties with $(u,v)$ replaced by $(u^{-1}, v^{-1})$.
In view of (\ref{eq:zeta-u,v}) and (\ref{eq:minors-tau}), 
if $\Delta (y_-) > 0$ for any $\Delta \in F(v)$ then 
$\Delta (y'_-) > 0$ for any $\Delta \in F(v^{-1})$.
Similarly, using (\ref{eq:tau-T}) we obtain that 
if $\Delta (y_+^T) > 0$ for any $\Delta \in F(u^{-1})$ then 
$\Delta ({y'_+}^T) > 0$ for any $\Delta \in F(u)$.

It remains to show that $y'_0 \in H_{> 0}$.
Applying the character $a \to a^{\omega_i}$ to both sides 
of the equality $y'_0 = [y_+ \overline u]_0 \ y_0^{-1} \ 
[\overline{\overline{v}} y_-]_0$ in (\ref{eq:zeta-u,v}),
we obtain
$$(y'_0)^{\omega_i} = \Delta_{u \omega_i, \omega_i} (y_+^T) \, 
y_0^{-\omega_i} \, \Delta_{v^{-1} \omega_i, \omega_i} (y_-) \ .$$
Since $\Delta_{v^{-1} \omega_i, \omega_i} \in F(v)$ and  
$\Delta_{u \omega_i, \omega_i} \in F(u^{-1})$ it follows that
$(y'_0)^{\omega_i} > 0$ for any $i \in [1,r]$.
Therefore, $y'_0 \in H_{> 0}$, as desired.
\endproof

\subsection{Proof of Theorem~\ref{th:t-through-x}}
\label{sec:t-through-x}

First notice that the equivalence of (\ref{eq:a-via-Q})
and (\ref{eq:a-through-x}) follows by applying the character
$a \mapsto a^{\omega_i}$ to both sides of (\ref{eq:a-through-x})
and simplifying the result.
In proving (\ref{eq:t-through-x}) and (\ref{eq:a-through-x}),
we will follow the same strategy that was used in the proof of
Theorem~\ref{th:LusztigTP}: first treat the case when
$\ii$ is unmixed, and then extend the result to the general case with 
the help of commutation relations (\ref{eq:H-conjugation}),
(\ref{eq:2-commutation-mixed}) and (\ref{eq:5-commutation-generic}).

Let us first assume that $\ii\in R(u,v)$ is unmixed, i.e, all 
the indices from $[\overline 1, \overline r]$ precede those from $[1,r]$. 
Repeatedly using (\ref{eq:H-conjugation}), we conclude that 
in this case $x= x_\ii (a; t_1, \dots, t_m) \in G^{u,v}_0$, and 
the components in the Gaussian decomposition of $x$ are given by
\beal
\label{eq:x-unmixed-Gauss}
[x]_- = x_{i_1} (a^{-\alpha_{|i_1|}} t_1) \cdots 
x_{i_{\l(u)}} (a^{-\alpha_{|i_{\l (u)}|}} t_{\l (u)}) \ , \\[.1in]
[x]_0 = a \ , \,\, [x]_+ = x_{i_{\l(u)+1}} (t_{\l(u)+1}) \cdots 
x_{i_m} (t_m) \ .
\eea
Since by Theorem~\ref{th:zeta-regularity}, 
$x = \zeta^{u^{-1},v^{-1}} (x')$, formula (\ref{eq:x'-diagonal})
implies that
\beal
\label{eq:a-unmixed}
a = [x]_0 = [\overline {\overline{u}} x']_0^{-1} [x']_0 
[x' \overline {v}]_0^{-1} \ .
\eea
This proves (\ref{eq:a-through-x}) since
a simple inspection shows that
the right-hand side of (\ref{eq:a-unmixed}) is equal to 
that of (\ref{eq:a-through-x}) when $\ii$ is unmixed.

Turning to the proof of (\ref{eq:t-through-x}),
let us first consider the case $\l (u) < k \leq m$.
Let $(y_+, y_0, y_-)$ be the $y$-coordinates of $x$,
and $(y'_+, y'_0, y'_-)$ be the $y$-coordinates of $x'$.
By (\ref{eq:Ansatz-N^w}), we have
\begin{equation}
\label{eq:t-through-y-}
t_k = 
\displaystyle\frac{1}{\Delta_{v^{-1} v_{<k} \omega_{i_k}, \omega_{i_k}}(y_-) 
\Delta_{v^{-1} v_{<k+1} \omega_{i_k}, \omega_{i_k}}(y_-)}
\displaystyle\prod_{j \ne i_k} 
\Delta_{v^{-1} v_{<k+1} \omega_j, \omega_j}(y_-)^{-a_{j,i_k}} \ .
\end{equation}
Using (\ref{eq:zeta-u,v}), (\ref{eq:minors-tau}), (\ref{eq:minors-Gauss1}),
(\ref{eq:y-to-x}), and (\ref{eq:diagonal via y}), we can rewrite
$\Delta_{v^{-1} v_{<k} \omega_j, \omega_j}(y_-)$ as follows:
\beal
\label{eq:Delta(y-)}
\Delta_{v^{-1} v_{<k} \omega_j, \omega_j}(y_-) = 
\Delta_{v \omega_j, v_{<k} \omega_j}(y'_-) = 
\Delta_{\omega_j, v_{<k} \omega_j}(\overline {\overline {v^{-1}}} y'_-) = 
\\ [.1in]
\Delta_{\omega_j, v_{<k} \omega_j}([\overline {\overline {v^{-1}}} y'_-]_+)
[\overline {\overline {v^{-1}}} y'_-]_0^{\omega_j} = 
\Delta_{\omega_j, v_{<k} \omega_j} (x') [y']_0^{- \omega_j} 
[\overline {\overline {v^{-1}}} y'_-]_0^{\omega_j} = \\ [.1in]
\Delta_{\omega_j, v_{<k} \omega_j} (x')/\Delta_{\omega_j, v \omega_j} (x') \ .
\eea
Substituting the expressions given by (\ref{eq:Delta(y-)}) into
(\ref{eq:t-through-y-}), we express $t_k$ as a Laurent monomial
in the minors $\Delta_l (x')$ given by (\ref{eq:Q-factors}).
Using the notation from Section~\ref{sec:a,t_k}, 
this monomial can be written as follows:
\beal
\label{eq:t_k-monomial}
\quad t_k = \displaystyle\frac{1}{\Delta_k (x') \Delta_{k^+} (x')} 
\displaystyle\prod_{l: l^- < k < l} \Delta_l (x')^{-a_{|i_l|,i_k}} 
\displaystyle\prod_{j \in [1,r]} \Delta_{m+j} (x')^{a_{j,i_k}} \ ,
\eea
where the index $k^+$ is defined by $(k^+)^- = k$. 
Formula (\ref{eq:t-through-x}) now follows by simple inspection 
which shows that, for $\ii$ unmixed, the right-hand side of
(\ref{eq:t-through-x}) is equal to the one of (\ref{eq:t_k-monomial}).

The proof of (\ref{eq:t-through-x}) for $1 \leq k \leq \l(u)$ 
is practically the same as above. 
In this case, the counterpart of (\ref{eq:t_k-monomial}) is given by
\begin{equation}
\label{eq:at-monomial}
a^{- \alpha_{|i_k|}} t_k = 
\displaystyle\frac{1}{\Delta_k (x') \Delta_{k^+} (x')} 
\displaystyle\prod_{l: l^- < k < l} \Delta_l (x')^{-a_{|i_l|,|i_k|}} 
\displaystyle\prod_{l: l^- = 0} \Delta_l (x')^{a_{|i_l|,|i_k|}} \ .
\end{equation}
To deduce (\ref{eq:t-through-x}) from (\ref{eq:at-monomial}), 
first notice that in view of 
(\ref{eq:a-via-Q}) and (\ref{eq:weight}), we have
$$ a^{\alpha_{|i_k|}} = 
\displaystyle\prod_{1 \leq l \leq m+r} 
\Delta_l (x')^{(\varepsilon_l - \varepsilon_{l^-})a_{|i_l|,|i_k|}} \ .$$
Thus in order to check (\ref{eq:t-through-x}), it suffices to show that,
for $\ii$ unmixed and $1 \leq k \leq \l (u)$, 
the right-hand side of (\ref{eq:at-monomial}) is equal to 
$$\prod_{l= 1}^{m+r} 
\Delta_l (x')^{(\varepsilon_{l^-} + \chi(k,l^-) - 
\varepsilon_{l} - \chi(k,l)) \ a_{|i_l|,|i_k|}} \ ; $$
this is again checked by direct inspection. 
                                                                            
Now let us prove (\ref{eq:t-through-x}) and (\ref{eq:a-through-x}) 
for an arbitrary double reduced word for $u$ and $v$. 
Every such word can be obtained from an unmixed one
by a sequence of \emph{mixed moves} of the form
$$ \cdots \overline j \, i \cdots \leadsto 
\cdots i \, \overline j \cdots \ .$$
It therefore suffices to prove the following statement.

\begin{lemma}
\label{lem:2-move-transition}
Suppose a reduced word $\ii' \in R(u,v)$ is obtained from $\ii \in R(u,v)$ 
by a mixed move. 
If {\rm (\ref{eq:t-through-x})} and {\rm (\ref{eq:a-through-x})} hold 
for $\ii$, then they also hold for~$\ii'$.
\end{lemma}

\proof
Suppose $\ii'$ is obtained from $\ii$ by interchanging
$i_k = \overline j$ and $i_{k+1} = i$. 
By (\ref{eq:H-conjugation}), (\ref{eq:2-commutation-mixed}) and 
(\ref{eq:5-commutation-generic}), the factorization parameters
that appear in two factorizations
$$x = x_\ii (a; t_1, \ldots, t_m)  = x_{\ii'} (a'; t'_1, \ldots, t'_m)$$
of an element $x \in G^{u,v}$, are related as follows.
If $i \neq j$, then 
\beal
\label{eq:trivial-2-move-transition}
a' = a  \ ,  \,\, t'_l = t_l \,\, (l \notin\{ k, k+1\}) \ , 
t'_k = t_{k+1}, \,\, t'_{k+1} = t_{k} \ .
\eea
For $i = j$, a straightforward calculation using (\ref{eq:H-conjugation}),
(\ref{eq:5-commutation-generic}) and (\ref{eq:epsilon}) shows that 
\beal
\label{eq:mixed-2-move-transition1}
(a')^{\omega_p} = a^{\omega_p} (1 + t_k t_{k+1})^{- \delta_{i,p}} \ , 
\eea
\beal
\label{eq:mixed-2-move-transition2}
t'_l = t_l (1 + t_k t_{k+1})^{(2 \varepsilon_l - 1) a_{i,|i_l|}}
\,\, (l < k) \ ,
\eea
\beal
\label{eq:mixed-2-move-transition3}
t'_k = t_{k+1} (1 + t_k t_{k+1}), \,\, 
t'_{k+1} = t_{k} (1 + t_k t_{k+1})^{-1} \ , 
\eea
\beal
\label{eq:mixed-2-move-transition4}
t'_l = t_l \,\, (l > k+1) \ .
\eea
We need to show the following: if we substitute the parameters 
$t_l$ and $a^{\omega_p}$ given by 
(\ref{eq:t-through-x}) and (\ref{eq:a-via-Q}) into
(\ref{eq:trivial-2-move-transition})--(\ref{eq:mixed-2-move-transition4}),
then the resulting $t'_l$ and $(a')^{\omega_p}$ 
satisfy the same formulas (\ref{eq:t-through-x}) 
and (\ref{eq:a-via-Q}) with $\ii$ replaced by $\ii'$.
This is immediate from the definitions when $i \neq j$, 
so let us assume $i = j$. 
By the definition (\ref{eq:Q-factors}), we have 
$\Delta_{l,\ii} = \Delta_{l,\ii'}$ for $l \neq k+1$, 
so we will denote this minor simply by~$\Delta_l\,$.
The key calculation is now as follows.

\begin{lemma}
\label{lem:1+tt'}
In the above notation, if $t_k$ and $t_{k+1}$ satisfy 
{\rm (\ref{eq:t-through-x})}, then
\beal
\label{eq:1+tt'-encoded}
1 + t_k t_{k+1} = \displaystyle\frac{\Delta_{k+1,\ii} (x')
\Delta_{k+1,\ii} (x')} {\Delta_{k} (x') \Delta_{(k+1)^+} (x')} \ .
\eea
\end{lemma}

\proof
Let us denote $u' = u_{\geq k+2}$ and $v' = v_{<k}$
(this is unambiguous since these expressions are the same for $\ii$
and $\ii'$, cf.~(\ref{eq:v_leq k})). 
In view of (\ref{eq:Q-factors}), we have
\beal 
\label{eq:exchanged minors}
\Delta_{k+1,\ii} = \Delta_{u' \omega_{i}, v' \omega_{i}} \ , \,\,
\Delta_{k+1,\ii'} = \Delta_{u's_i \omega_{i}, v's_i \omega_{i}} \ , \\ [.1in]
\Delta_{k} = \Delta_{u's_i \omega_{i}, v' \omega_{i}} \ , \,\,
\Delta_{(k+1)^+} = \Delta_{u'\omega_{i}, v's_i \omega_{i}} \ ,
\eea
so (\ref{eq:1+tt'-encoded}) takes the form
\beal
\label{eq:1+tt'}
1 + t_k t_{k+1} = 
\displaystyle\frac{\Delta_{u' \omega_i, v' \omega_i} (x')
\Delta_{u' s_i \omega_i, v' s_i \omega_i} (x')}
{\Delta_{u' s_i \omega_i, v' \omega_i} (x')
\Delta_{u' \omega_i, v' s_i \omega_i} (x')} \ .
\eea
On the other hand, if $t_k$ and $t_{k+1}$ are given by
(\ref{eq:t-through-x}), then 
$$t_k t_{k+1} = \prod_{l= 1}^{m+r} 
\Delta_l (x')^{(\chi(k,l^-) + \chi(k+1,l^-) - 
\chi(k,l) - \chi(k+1,l)) \ a_{|i_l|, i}}  \ .$$
By (\ref{eq:chi}), we have
\begin{equation}
\chi(k,l) + \chi(k+1,l) =
\left\{
\begin{array}{ll}
1 &\textrm{if } l > k+ 1 \,;\\[.1in]
1/2 &\textrm{if } l \in \{k,k+1\} \, ; \\[.1in]
0 &\textrm{if } l < k \,.
\end{array}
\right.
\end{equation}
It follows that
\begin{eqnarray*}
\begin{array}{rcl}
t_k t_{k+1} &=& \Delta_{k} (x')^{-1} \Delta_{(k+1)^+} (x')^{-1}
\prod_{l: l^- < k < l} \Delta_l (x')^{- a_{|i_l|, i}} \\[.1in]
&=&(\Delta_{u' s_i \omega_i, v' \omega_i} (x')
\Delta_{u' \omega_i, v' s_i \omega_i} (x'))^{-1} 
\prod_{p \neq i} \Delta_{u'\omega_p, v' \omega_p} (x')^{- a_{pi}}
\ .
\end{array}
\end{eqnarray*}
Therefore, (\ref{eq:1+tt'}) becomes a consequence of 
(\ref{eq:minors-Dodgson}), and we are done. \endproof

If we substitute the expressions given by 
(\ref{eq:t-through-x}), (\ref{eq:a-via-Q}) and (\ref{eq:1+tt'-encoded})
into the formulas 
(\ref{eq:trivial-2-move-transition})--(\ref{eq:mixed-2-move-transition4}),
then an easy simplification shows that they will be given by 
(\ref{eq:t-through-x}) and (\ref{eq:a-via-Q}) for $\ii'$. 
This completes the proofs of Lemma~\ref{lem:2-move-transition} and
Theorem~\ref{th:t-through-x}.
\endproof
 
\subsection{Proofs of Theorems~\ref{th:TP-criteria}
and~\ref{th:TP-base}}
\label{sec:ProofsTP}

We start by recalling a well-known property of reduced words in 
Coxeter groups (cf.~\cite{bourbaki,humphreys}). 
To state it, we will need the following notion.

We call a \emph{$d$-move} the transformation of a reduced word
that replaces $d$ consecutive entries
$i,j,i,j, \ldots$ by $j,i,j,i,\ldots$, 
for some $i$ and $j$ such that $d$ is the order of~$s_i s_j\,$. 
Note that, for given $i$ and $j$, 
the value of $d$ can be determined from the Cartan matrix as follows: 
if $a_{ij}a_{ji} = 0$ (resp.\ $1,2,3$), then $d=2$ (resp.\ $3,4,6$).  

\begin{proposition}
\label{pr:d-moves}
Every two reduced words for the same element of a Coxeter group 
can be obtained from each other by a sequence of $d$-moves. 
\end{proposition}

Applying this proposition to the group $W \times W$,
we conclude that every two double reduced words $\ii, \ii' \in
R(u,v)$ can be obtained from each other by a sequence of the
following operations: $d$-moves for each of the alphabets $[1,r]$
and $[\overline 1, \overline r]$, and also mixed moves (cf.\
Section~\ref{sec:t-through-x}) and their inverses.  
\medskip

\emph{Proof of Theorem~\ref{th:TP-base}.}
Let us first prove that $F(\ii)$ is a 
transcendence basis for the field~$\CC(G^{u,v})$.  
By Theorem~\ref{th:Guv-affine}, $F(\ii)$ is of cardinality
$\dim\,G^{u,v}$. 
It is therefore enough to show that $F(\ii)$
generates~$\CC(G^{u,v})$. 
In view of Theorem~\ref{th:zeta-regularity}, 
it suffices to show that the collection of ``twisted" minors 
$\Delta_{k,\ii} (x')$ (cf.~(\ref{eq:Q-factors})) generates~$\CC(G^{u,v})$. 
By Theorem~\ref{th:double}, the field~$\CC(G^{u,v})$ is generated
by the factorization parameters $t_k$ and $a^{\omega_i}$, 
and the claim follows by Theorem~\ref{th:t-through-x}.

The second statement of the theorem is a consequence of the
following lemma. 

\begin{lemma}
\label{lem:set-diff}
Suppose a double reduced word $\ii'$ is obtained from $\ii$
by a $d$-move in one of the alphabets $[1,r]$ and 
$[\overline 1, \overline r]$, or by a mixed move, or by 
the inverse of a mixed move. 
Then each element of the set difference $F(\ii') \setminus
F(\ii)$ is a subtraction-free expression in the elements of
$F(\ii)$.  
\end{lemma}

\proof 
For $d$-moves in $[1,r]$ or $[\overline 1, \overline r]$, 
the desired subtraction-free expressions
can be obtained from generalized Pl\"ucker relations in 
Theorem~\ref{th:minors-Plucker} (including the omitted relations
of type~$G_2$); 
this part of the argument is essentially borrowed
from~\cite[Proposition~6.10]{BZ}. 
For mixed moves and their inverses, the statement
follows in the same way from Theorem~\ref{th:minors-Dodgson}
(cf.\ the proof of Lemma~\ref{lem:2-move-transition} above). 
Lemma~\ref{lem:set-diff} and Theorem~\ref{th:TP-base} are proved.
\endproof

\emph{Proof of Theorem~\ref{th:TP-criteria}.} 
It will suffice to show that the following are equivalent:

\noindent (1) $x \in G^{u,v}_{>0}$;

\noindent (2) $\Delta (x) > 0$ for any $\Delta \in F(u,v)$;

\noindent (3) $\Delta (x) > 0$ for any $\Delta \in F(\ii)$.

The equivalence of (2) and (3) follows from Theorem~\ref{th:TP-base}. 
Let us show that $(1) \Rightarrow (3)$.
By Theorem~~\ref{th:zeta-positivity}, if $x \in G^{u,v}_{>0}$ then
$x' \in G^{u^{-1},v^{-1}}_{>0}$. 
The condition (3) now follows by applying
Theorems~\ref{th:LusztigTP} and \ref{th:inverse-monomial}
to $x'$ and the reduced word $\ii^* \in R(u^{-1}, v^{-1})$ opposite
to $\ii$.

It remains to show that $(2) \Rightarrow (1)$.
First of all, by Corollary~\ref{cor:G_0-description}, 
(2) implies that $x \in G^{u,v}_{0}$; 
moreover, $[x]_0 \in H_{>0}\,$.
By Proposition~\ref{pr:T-properties} and 
Theorems~\ref{th:TP-through-y+-0} and \ref{th:TP-criteria-y}, 
is suffices to show that $y_- = \pi_- (x)$ satisfies
$\Delta_{v^{-1} \omega_i, v' \omega_i} (y_-) > 0$ for
$i \in [1,r]$ and any $v' \preceq v^{-1}$. 
Using (\ref{eq:y-to-x}) and (\ref{eq:minors-Gauss1}), 
we obtain (cf.~(\ref{eq:Delta(y-)})): 
\beal
\label{eq:minors-y_-}
\Delta_{v^{-1} \omega_i, v' \omega_i} (y_-)  = 
\displaystyle\frac{\Delta_{\omega_i, v' \omega_i} (x)}
{\Delta_{\omega_i, v^{-1} \omega_i} (x)} \ .
\eea
By (2), the right-hand side of (\ref{eq:minors-y_-})
is positive, and the proof is complete.
\endproof


\section{$GL_n$ theory}
\label{sec:gln}

Throughout this section, $G=GL_n (\CC)$ is the group of invertible
$n\times n$ matrices with complex entries.
In this case, the problems under consideration  become quite
natural questions in ``classical" linear algebra, so we will 
formulate them here---and state our main results---in an elementary
and self-contained way.
We will not give any proofs though, 
since these results can be easily derived from the type~$A$
specializations of the corresponding statements in
Section~\ref{sec:main}; pointers to these statements are provided,
wherever appropriate.

\subsection{Bruhat cells and double Bruhat cells for $GL_n\,$}

Our first object of interest are the double Bruhat cells.
Let us introduce them for the group $G = GL_n (\CC)$.
We will need some notation.
Let $B$ (resp. $B_-$) be the subgroup of 
upper-triangular (resp. lower-triangular) matrices in $G$.
Let $W = S_n$ be the symmetric group acting on the set 
$[1,n] = \{1, \dots, n\}$; we will think of $W$ as a subgroup
of $G$ by identifying a permutation $w$ with the matrix
$w = (\delta_{i, w(j)})$.
The double cosets $BwB$ and $B_- w B_-$ are called Bruhat cells
(with respect to $B$ and $B_-$, respectively).
The group $G$ has two Bruhat decompositions into a disjoint union of 
Bruhat cells (see, e.g., \cite[Section~2.4]{alperin-bell}):
\[
G = \bigcup_{u \in W} B u B = \bigcup_{v \in W} B_- v B_-  \ . 
\]
The \emph{double Bruhat cells}~$G^{u,v}$ are defined by 
\[
G^{u,v} = B u B  \cap B_- v B_- \ ;
\]
thus $G$ is the disjoint union of all $G^{u,v}$ for $(u,v)\in
W\times W$. 

As an algebraic variety, a double Bruhat cell $G^{u,v}$  
is biregularly isomorphic to a Zariski
open subset of an affine space of dimension $n+\l(u)+\l(v)$,
where $\l (u)$ is the number of inversions of a permutation $u$;
cf.\ Theorem~\ref{th:Guv-affine}.
(In other words, $G^{u,v}$ is isomorphic, as an algebraic variety, 
to a subset of $\CC^{n+\l(u)+\l(v)}$ 
obtained by excluding common zeroes of a finite set of polynomials.) 

Each Bruhat cell (hence each double Bruhat cell) 
can be described explicitly by a set of conditions specifying
vanishing and non-vanishing of certain minors. 
Let us denote by $\Delta_{I,J}$ the minor with the row set $I$ and the 
column set $J$; here $I$ and $J$ are two subsets of the same size
in $[1,n]$, and the minor is viewed as a (regular) function on~$G$.
(This notation corresponds to that of Definition~\ref{def:general minors}, 
as follows: the function $\Delta_{u \omega_i, v \omega_i}$ in 
(\ref{eq:Delta-general}) becomes the minor $\Delta_{u([1,i]), v([1,i])}$.)
The following description of Bruhat cells is probably the most ``economical."

\begin{proposition}
\label{pr:Bruhat-GL}
A matrix $x \in G$ belongs to the Bruhat cell $BwB$ if and only if
it satisfies the following conditions:
\begin{itemize}
\item $\Delta_{w([1,i]), [1,i]} \neq 0$ for $i = 1, \dots, n-1$;
\smallskip

\item $\Delta_{w([1,i-1] \cup \{j\}), [1,i]} = 0$ for all $(i,j)$ such
that $1 \leq i < j \leq n$ and $w(i) < w(j)$.
\end{itemize}
\end{proposition}
 
This proposition can be proved by specializing 
Propositions~\ref{pr:Bruhat cell} and 
\ref{pr:minors-Gauss} and Corollary~\ref{cor:G_0-description}. 
Notice that in our present situation the subgroup 
$N_- (w)$ (cf.\ (\ref{eq:N+(u)})) consists of all unipotent 
lower-triangular matrices $y$ such that 
$y_{ij} = 0$ whenever $w(i) > w(j)$.

The transpose map $x \mapsto x^T$ transforms a minor $\Delta_{I,J}$
into $\Delta_{J,I}$ and sends a Bruhat cell $BwB$ to $B_- w^{-1} B_-\,$.
Thus Proposition~\ref{pr:Bruhat-GL} implies a similar description 
of the opposite Bruhat cells $B_- w B_-\,$.
Combining the two sets of conditions yields an explicit description 
of the double Bruhat cells.

\subsection{Factorization problem for $GL_n\,$}

In the situation under consideration, 
the maximal torus $H = B \cap B_-$ in $G$ 
is the subgroup of invertible diagonal matrices.
Thus $H$ is naturally identified with $\CC_{\neq 0}^n$ 
by taking the diagonal entries as coordinates.
This allows us to state the factorization problem 
of Section~\ref{sec:factorization-problem} in a more symmetric form,
as follows.

Let $E_{i,j}$ denote the $n \times n$ matrix whose 
$(i,j)$-entry is equal to~1 while all other entries are~0;
let $I \in G$ denote the identity matrix.
For $i=1,\dots,n-1$, let
\beal
\label{eq:jacobi}
x_i(t)=I+tE_{i,i+1}=
\left(\begin{array}{cccccc}
1      & \cdots & 0      & 0      & \cdots & 0      \\
\cdots & \cdots & \cdots & \cdots & \cdots & \cdots \\
0      & \cdots & 1      & t      & \cdots & 0      \\
0      & \cdots & 0      & 1      & \cdots & 0      \\
\cdots & \cdots & \cdots & \cdots & \cdots & \cdots \\
0      & \cdots & 0      & 0      & \cdots & 1 \\
\end{array}\right)
\eea
and
\beal
x_{\bar i}(t)=I+tE_{i+1,i}=
\left(\begin{array}{cccccc}
1      & \cdots & 0      & 0      & \cdots & 0      \\
\cdots & \cdots & \cdots & \cdots & \cdots & \cdots \\
0      & \cdots & 1      & 0      & \cdots & 0      \\
0      & \cdots & t      & 1      & \cdots & 0      \\
\cdots & \cdots & \cdots & \cdots & \cdots & \cdots \\
0      & \cdots & 0      & 0      & \cdots & 1 \\
\end{array}\right)\ .
\eea
Also, for $i=1,\dots, n$ and $t \neq 0$, let
\beal
\label{eq:diag-jacobi}
x_\circi(t)=I+(t-1)E_{i,i}=
\left(\begin{array}{cccccc}
1      & \cdots & 0      & 0      & \cdots & 0      \\
\cdots & \cdots & \cdots & \cdots & \cdots & \cdots \\
0      & \cdots & t      & 0      & \cdots & 0      \\
0      & \cdots & 0      & 1      & \cdots & 0      \\
\cdots & \cdots & \cdots & \cdots & \cdots & \cdots \\
0      & \cdots & 0      & 0      & \cdots & 1 \\
\end{array}\right)\ .
\eea
The matrices defined in (\ref{eq:jacobi})--(\ref{eq:diag-jacobi}) 
are called \emph{elementary Jacobi matrices}.
It is easy to see that these matrices generate $G$ as a group.

Consider the alphabet of $3n-2$ symbols 
\beal
\label{eq:alphabet}
\A=\{1, \dots, n-1, 
\circled{\,1}\,, \dots, \circled{\,n}\,, 
\overline{1}, \dots, \overline{n-1}\}\ .
\eea
The formulas (\ref{eq:jacobi})--(\ref{eq:diag-jacobi})
associate a matrix $x_i(t) \in G$
to any symbol $i\in\A$ and any $t\in\CC_{\neq 0}\,$.
An analogue of the product map (\ref{eq:productmap})
is now defined as follows: to any sequence
$\ii = (i_1, \dots, i_l)$ of symbols in~$\A$,
we associate the map $x_\ii : \CC_{\neq 0}^l \to G$ defined by
\begin{equation}
\label{eq:t(x)}
x_{\ii}(t_1,\dots,t_l)=x_{i_1}(t_1)\cdots x_{i_l}(t_l) \ .
\end{equation}
(Thus the difference with (\ref{eq:productmap}) is that now 
the factor $a \in H$ is split into elementary factors,
which are allowed to be spread along the factorization.)

For instance, the sequence
$\,\ii=\circled{\,1}\, \bar 1 \circled{\,2}\, 1$
gives rise to the map
\begin{equation}
\label{eq:GL2-xxxx}
(t_1,t_2,t_3,t_4)\,\mapsto\, 
\mat{t_1}{0}{0}{1}
\mat{1}{0}{t_2}{1}
\mat{1}{0}{0}{t_3}
\mat{1}{t_4}{0}{1} =  \mat{t_1}{t_1t_4}{t_2}{t_2t_4\!+\!t_3} .
\end{equation}

The matrix $x = x_\ii (t_1, \dots, t_l)$ has a simple combinatorial
description in terms of \emph{planar networks}.
This description (cf.\ \cite{brenti} and references therein)
generalizes the one in~\cite[Section~2.4]{BFZ}, and provides
combinatorial formulas for the minors of~$x$ as polynomials with
nonnegative coefficients in the variables $t_1,\dots,t_l\,$.

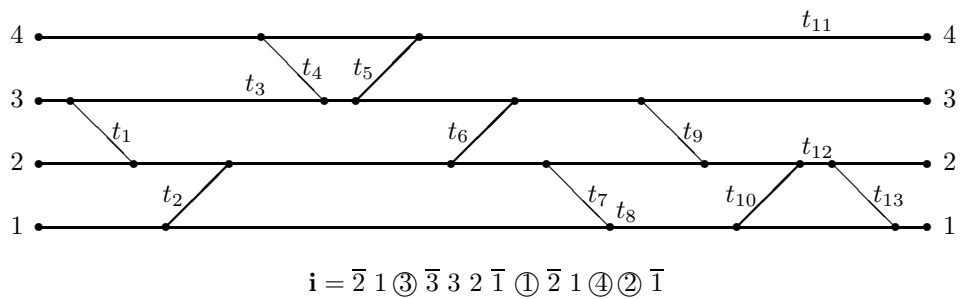
\begin{figure}[hb]
\setlength{\unitlength}{1.2pt} 
\begin{center}
\begin{picture}(280,75)(0,-15)
\thicklines

\put(0,0){\line(1,0){280}}
\put(0,20){\line(1,0){280}}
\put(0,40){\line(1,0){280}}
\put(0,60){\line(1,0){280}}

\put(40,0){\line(1,1){20}}
\put(100,40){\line(1,1){20}}
\put(130,20){\line(1,1){20}}
\put(220,0){\line(1,1){20}}

\thinlines
\put(10,40){\line(1,-1){20}}
\put(70,60){\line(1,-1){20}}
\put(160,20){\line(1,-1){20}}
\put(190,40){\line(1,-1){20}}
\put(250,20){\line(1,-1){20}}

  \put(285,-2){${1}$}
  \put(285,18){${2}$}
  \put(285,38){${3}$}
  \put(285,58){${4}$}
  \put(-9,-2){${1}$}
  \put(-9,18){${2}$}
  \put(-9,38){${3}$}
  \put(-9,58){${4}$}

  \put(0,0){\circle*{2.5}}
  \put(0,20){\circle*{2.5}}
  \put(0,40){\circle*{2.5}}
  \put(0,60){\circle*{2.5}}
  \put(280,0){\circle*{2.5}}
  \put(280,20){\circle*{2.5}}
  \put(280,40){\circle*{2.5}}
  \put(280,60){\circle*{2.5}}

  \put(40,0){\circle*{2.5}}
  \put(180,0){\circle*{2.5}}
  \put(220,0){\circle*{2.5}}
  \put(270,0){\circle*{2.5}}

  \put(30,20){\circle*{2.5}}
  \put(60,20){\circle*{2.5}}
  \put(130,20){\circle*{2.5}}
  \put(160,20){\circle*{2.5}}
  \put(210,20){\circle*{2.5}}
  \put(240,20){\circle*{2.5}}
  \put(250,20){\circle*{2.5}}

  \put(10,40){\circle*{2.5}}
  \put(90,40){\circle*{2.5}}
  \put(100,40){\circle*{2.5}}
  \put(150,40){\circle*{2.5}}
  \put(190,40){\circle*{2.5}}

  \put(70,60){\circle*{2.5}}
  \put(120,60){\circle*{2.5}}

  \put(39,8){$t_2$}
  \put(99,48){$t_5$}
  \put(129,28){$t_6$}
  \put(216.5,8){$t_{10}$}


  \put(23,28){$t_1$}
  \put(83,48){$t_4$}
  \put(173,8){$t_7$}
  \put(203,28){$t_9$}
  \put(263,8){$t_{13}$}


  \put(65,43){$t_3$}
  \put(182,3){$t_8$}
  \put(240.5,63){$t_{11}$}
  \put(240.5,23){$t_{12}$}
\put(85,-20)
{$\ii=\overline 2~~ 1\circled{\,3}\,\, \overline 3~~3~~2~~\overline 1\,
\circled{\,1}\,\,\overline 2~~1\circled{\,4}\circled{\,2}~\overline 1$}

\end{picture}
\end{center}
\caption{Planar network}
\label{fig:planar-network}
\end{figure}

The planar network~$\Gamma(\ii)$ associated to a sequence~$\ii$ of symbols 
from~$\A$ (see Figure~\ref{fig:planar-network})
is constructed as a concatenation of
``elementary'' networks that correspond to the parameters 
$t_1,\dots,t_l$ (in this order).
Each unbarred, barred, or circled entry $i_k$ of $\ii$ corresponds to a
fragment of one of the following three kinds, respectively:

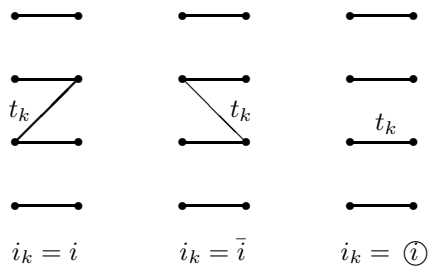
\begin{figure}[ht]
\setlength{\unitlength}{1.2pt} 
\begin{center}
\begin{picture}(50,70)(0,-10)
\thicklines

\put(5,0){\line(1,0){20}}
\put(5,20){\line(1,0){20}}
\put(5,40){\line(1,0){20}}
\put(5,60){\line(1,0){20}}

\put(5,20){\line(1,1){20}}

  \put(5,0){\circle*{2.5}}
  \put(5,20){\circle*{2.5}}
  \put(5,40){\circle*{2.5}}
  \put(5,60){\circle*{2.5}}

  \put(25,0){\circle*{2.5}}
  \put(25,20){\circle*{2.5}}
  \put(25,40){\circle*{2.5}}
  \put(25,60){\circle*{2.5}}

  \put(3,28){$t_k$}

  \put(4,-17){$i_k=i$}

\end{picture}
\begin{picture}(50,65)(0,-10)
\thicklines

\put(5,0){\line(1,0){20}}
\put(5,20){\line(1,0){20}}
\put(5,40){\line(1,0){20}}
\put(5,60){\line(1,0){20}}

  \put(5,0){\circle*{2.5}}
  \put(5,20){\circle*{2.5}}
  \put(5,40){\circle*{2.5}}
  \put(5,60){\circle*{2.5}}

  \put(25,0){\circle*{2.5}}
  \put(25,20){\circle*{2.5}}
  \put(25,40){\circle*{2.5}}
  \put(25,60){\circle*{2.5}}
\thinlines

\put(5,40){\line(1,-1){20}}
  \put(20,28){$t_k$}

  \put(4,-17){$i_k=\bar i$}

\end{picture}
\begin{picture}(50,65)(0,-10)
\thicklines

\put(5,0){\line(1,0){20}}
\put(5,20){\line(1,0){20}}
\put(5,40){\line(1,0){20}}
\put(5,60){\line(1,0){20}}

  \put(5,0){\circle*{2.5}}
  \put(5,20){\circle*{2.5}}
  \put(5,40){\circle*{2.5}}
  \put(5,60){\circle*{2.5}}

  \put(25,0){\circle*{2.5}}
  \put(25,20){\circle*{2.5}}
  \put(25,40){\circle*{2.5}}
  \put(25,60){\circle*{2.5}}

  \put(13,24){$t_k$}

  \put(2,-17){$i_k=\circled{\,i}$}

\end{picture}
\end{center}

\end{figure}

\noindent
(a diagonal edge connects horizontal levels $i$ and $i+1$;
in the examples above, $i=2$).
These fragments are the combinatorial equivalents of the elementary
matrices (\ref{eq:jacobi})--(\ref{eq:diag-jacobi}).
Each fragment has a distinguished edge whose weight is~$t_k\,$;
all other edges have weight~1. 
All edges  are presumed oriented left-to-right.

We number the sources and sinks of the network~$\Gamma(\ii)$ 
bottom-to-top, and 
define the weight of a path in~$\Gamma(\ii)$  
to be the product of the weights of all edges in the path.
One easily checks that the sum of these weights, over all paths
that connect a given source $i$ to a given sink $j$,
is nothing but the matrix element $x_{ij}$ of 
$x=  x_{\ii}(t_1,\dots,t_l)$. 

This observation can be generalized. 
Let us define the weight of a family of paths in~$\Gamma(\ii)$ 
to be the product of the weights of all paths in the family. 
Then the minors of $x$ are computed as follows.

\begin{proposition}
\label{pr:planar-minors}
A minor $\Delta_{I,J}(x)$ equals the sum of
weights of all families of vertex-disjoint paths in $\Gamma(\ii)$
connecting the sources labeled by~$I$ with the sinks labeled by~$J$.  
\end{proposition}

For example, in Figure~\ref{fig:planar-network} we have
$x_{21}=t_7t_8 + t_{12}t_{13}+t_6 t_9 t_{12}t_{13}$ and 
$\Delta_{12,12}(x) =t_8 t_{12} (1+t_6t_9)$.

We will be especially interested in a particular class of sequences~$\ii$ 
which we call factorization schemes (they are analogues of double
reduced words of Section~\ref{sec:factorization-problem}).

\begin{definition} {\rm 
Let $u$ and $v$ be two permutations in $W=S_n\,$.
A \emph{factorization scheme} of type $(u,v)$
is a word~$\ii$ of length $n + \l(u)+\l(v)$
in the alphabet~$\A$ which is an arbitrary shuffle of 
three words of the following kind: 
\begin{itemize}
\item a reduced word for $v$;
\item a reduced word for $u$, with all entries barred;
\item a permutation of the symbols $\circled{1}\,, \dots, \circled{n}\,$.
\end{itemize}
These three words will be called, respectively, the $E$-\emph{part}, 
the $F$-\emph{part}, and the $H$-\emph{part} of a factorization scheme~$\ii$. 
}\end{definition}


For example, let
\beal
\label{eq:4312-4213}
u=4312=s_2 s_3 s_1 s_2 s_1\in S_4 \ ,\\[.1in]
v= 4213=s_1 s_3 s_2 s_1\in S_4     \ ,
\eea
Then
\begin{equation}
\label{eq:example-fact-scheme}
\begin{array}{lcl}
\ii & = & 
\overline 2~~ 1\circled{\,3}\,\, \overline 3~~3~~2~~\overline 1\,
\circled{\,1}\,\,\overline 2~~1\circled{\,4}\circled{\,2}~\overline 1
\end{array}
\end{equation}
is a factorization scheme of type~$(u,v)$.

The following result is an analogue of Theorem~\ref{th:double}.

\begin{theorem}
\label{th:birational-2}
Let $u,v\in W=S_n$, and let $l= n + \l(u)+\l(v)$.
For any factorization scheme $\ii$ of type~$(u,v)$, 
the product map $x_\ii$ given by {\rm (\ref{eq:t(x)})} is a biregular 
isomorphism between $\CC_{\neq 0}^l$ and a Zariski open subset of 
the double Bruhat cell $G^{u,v}$.
\end{theorem}

The \emph{factorization problem} for $GL_n$ can be now formulated as follows:
for a given factorization scheme $\ii$, find explicit formulas for 
the components~$t_k$ in terms of the matrix $x = x_\ii (t_1, \dots, t_l)$.
By Theorem~\ref{th:birational-2}, each $t_k$ is a rational function
in the matrix entries of $x$.
For example, if $\,\ii=\circled{\,1}\, \bar 1 \circled{\,2}\, 1$, 
so that the map $x_\ii$ is given by (\ref{eq:GL2-xxxx}), then 
the solution to the factorization problem is given by 
\beal
\label{eq:GL2-inverse}
t_1=x_{11}\,,\ 
t_2=x_{21}\,,\ 
t_3=\displaystyle\frac{\det(x)}{x_{11}}\,,\ 
t_4=\displaystyle\frac{x_{12}}{x_{11}}\ .
\eea

\subsection{The twist maps for $GL_n\,$}
\label{sec:twist-GL}

As in the general case, our solution to the factorization problem 
for $G = GL_n$ will utilize the ``twist maps'' $\zeta^{u,v}:x\mapsto
x'$, which are defined for any two permutations $u$ and $v$. 
The definition (\ref{eq:zeta-u,v-x})
can be rewritten as 
\begin{equation}
\label{eq:gln-zeta}
x' = d_0 \, [x^T \overline{u}]_+\, \overline{u}^T \, (x^T)^{-1}\,
{\overline{v^{-1}}} \, [{\overline{v^{-1}}}^T x^T]_- \, d_0^{-1} \ ,  
\end{equation}
where the following notation is used.
The matrix $d_0$ is the diagonal $n\times n$
matrix with diagonal entries $1,-1,1,-1,\dots$.
For a matrix $z\in G$, $z^T$ stands for the transpose of~$z$, 
and $z=[z]_- [z]_0 [z]_+$ denotes the Gaussian decomposition of $z$
(also known 
as the $LDU$ decomposition). 
Finally, the matrix $\overline{w}$ is obtained from a 
permutation matrix for $w$ by the following modification:
an entry is changed from $1$ to $-1$ whenever it has an odd number 
of nonzero entries lying below and to the left of it.

By Theorem~\ref{th:zeta-regularity}, the right-hand side 
of~(\ref{eq:gln-zeta}) is well defined for 
any~$x\in G^{u,v}$, and the twist map $\zeta^{u,v}$
establishes a biregular isomorphism 
between $G^{u,v}$ and  $G^{u^{-1}, v^{-1}}$; 
the inverse isomorphism is $\zeta^{u^{-1},v^{-1}}$. 

We give below a few examples of explicitly computed twist maps.

\begin{example}{\rm
Let $G=GL_2(\CC)$ and $u=v=\wnot\,$. Then 
(cf.~Example~\ref{example:SL2-cells}) 
\[
\overline{u}=\overline{v^{-1}}=\mat{0}{-1}{1}{0}\ ,\ \ 
x'=  
\mat{x_{11}x_{12}^{-1}x_{21}^{-1}}{x_{21}^{-1}}{x_{12}^{-1}}
{x_{22}\,\det(x)^{-1}}\ .
\]
}\end{example}

\begin{example}{\rm
Let $G=GL_3(\CC)$ and $u=v=\wnot\,$. Then 
\begin{eqnarray*}
x'=  
\left( 
{\begin{array}{ccc}
{\displaystyle \frac {{{x}_{1  1}}}{{{x}_{3  1}}\,{{x}_{1  3}}}}
 &  {\displaystyle \frac {\Delta_{12,13}}{{{x}_{3  1}}\,{\Delta_{12,23}}}} & 
{\displaystyle \frac {1}{{{x}_{3  1}}}} \\ [2ex]
{\displaystyle \frac {\Delta_{13,12}}{{\,{{x}_{1  3}}\,\Delta_{23,12}}}} &  
{\displaystyle \frac {
x_{33}\Delta_{12,12}-\det(x)
}{
{\Delta_{23,12}}\,{\Delta_{12,23}}}} & {\displaystyle \frac {{{x}_{3  2}}}{
{\Delta_{23,12}}}} \\ [2ex]
{\displaystyle \frac {1}{{{x}_{1  3}}}} & {\displaystyle \frac {
{{x}_{2  3}}}{{\Delta_{12,23}}}} &   {\displaystyle \frac 
{\Delta_{23,23}}{\det(x)}}
\end{array}}
 \right)  
\end{eqnarray*}
}\end{example}

\begin{example}{\rm
Let $G=GL_4(\CC)$ and $u=v=\wnot\,$. Then $x'$ is equal to
\[
\!\left(\!\!\!
\begin{array}{cccc}
\displaystyle \frac{x_{11}}{x_{14}x_{41}} \!&\!
\displaystyle \frac{\Delta_{12,14}}{x_{41}\Delta_{12,34}} \!&\!
\displaystyle \frac{\Delta_{123,134}}{x_{41}\Delta_{123,234}} \!&\!
\displaystyle \frac{1}{x_{41}} 
\\[.3in]
\displaystyle \frac{\Delta_{14,12}}{x_{14}\Delta_{34,12}} \!&\!
\displaystyle \frac{x_{44}\Delta_{12,12}-\Delta_{124,124}}
{\Delta_{34,12}\Delta_{12,34}} \!&\!
\displaystyle \frac{x_{42}\Delta_{123,134}-x_{41}\Delta_{123,234}}
{\Delta_{34,12}\Delta_{123,234}} \!&\!
\displaystyle \frac{x_{42}}{\Delta_{34,12}}
\\[.3in]
\displaystyle \frac{\Delta_{134,123}}{x_{14}\Delta_{234,123}} \!&\!
\displaystyle \frac{x_{24}\Delta_{134,123}-x_{14}\Delta_{234,123}}
{\Delta_{12,34}\Delta_{234,123}} \!&\!
\displaystyle \frac{\Delta_{123,123}\Delta_{34,34}-x_{33}\det(x)}
{\Delta_{123,234}\Delta_{234,123}} \!&\!
\displaystyle \frac{\Delta_{34,23}}{\Delta_{234,123}} 
\\[.3in]
\displaystyle \frac{1}{x_{14}} \!&\!
\displaystyle \frac{x_{24}}{\Delta_{12,34}} \!&\!
\displaystyle \frac{\Delta_{23,34}}{\Delta_{123,234}} \!&\!
\displaystyle \frac{\Delta_{234,234}}{\det(x)} 
\end{array}
\!\!\!\!\right)\!.
\]
}
\end{example}

\begin{example}
\label{example:running-twist}
{\rm 
Let $n=4$, $u=4312$, and $v= 4213$ (cf.~(\ref{eq:4312-4213}).
Then 
\[
\overline u\!=\! 
\left( 
{\begin{array}{rrrr}
0 & 0 & 1 & 0 \\
0 & 0 & 0 & 1 \\
0 & -1 & 0 & 0 \\
1 & 0 & 0 & 0
\end{array}}
 \right) 
\ ,\ \ 
\overline {v^{-1}}\!=\! 
 \left( 
{\begin{array}{rrrr}
0 & 0 & 0 & -1 \\
0 & -1 & 0 & 0 \\
1 & 0 & 0 & 0 \\
0 & 0 & 1 & 0
\end{array}}
 \right) 
\]
and
\begin{eqnarray*}
x'= \left( 
\begin{array}{cccc}
{\displaystyle \frac {{{x}_{1 1}}}{{{x}_{4 1}
}\,{{x}_{1 3}}}} & 
{\displaystyle \frac {\Delta_{12,13}}{{{x}_{4 1}}\,{\Delta_{12,23}}}}
& 
{\displaystyle \frac {{\Delta_{123,123}}}{{{x}_{4 1}}\,{\Delta_{123,234}}
}}
& {\displaystyle \frac {1}{{{x}_{4 1}}}}  
\\[.3in]
{\displaystyle \frac {\Delta_{14,12}}{{\Delta_{34,12}}\,{{x}_{1 3}}}}
& 
{\displaystyle \frac {x_{43}\Delta_{12,12}-\Delta_{124,123}}
{{\Delta_{34,12}}\,{\Delta_{12,23}}}}
&
 {\displaystyle \frac {{\Delta_{123,123}}\,{{x}_{
4 2}}}{{\Delta_{34,12}}\,{\Delta_{123,234}}}} 
&
 {\displaystyle \frac {{{x}
_{4 2}}}{{\Delta_{34,12}}}} 
\\[.3in]
{\displaystyle \frac {1}{{{x}_{1 3}}}}
&
{\displaystyle \frac {{{x}_{2 3}}}{{\Delta_{12,23}}}}
& 
{\displaystyle \frac {\Delta_{23,23}}{{\Delta_{123,234}}}}
&
{\displaystyle \frac {\Delta_{23,23}}{{\Delta_{123,123}}}}
\\[.3in]
0 & {\displaystyle \frac {{{x}_{1 3}}}{{\Delta_{12,23}}}} 
& {\displaystyle \frac {\Delta_{13,23}}{{\Delta_{123,234}}}} &  
\displaystyle \frac{\Delta_{134,234}}{\det(x)}
\end{array}
\right)\ .
\end{eqnarray*}
Note that, in the course of computing 
the matrix elements of~$x'$ above, one has to take into account
the relations
\beac
\label{eq:gl4rel}
x_{14}=0\ ,\ \ 
x_{24}=0\ ,\ \ 
\Delta_{234,123}=0\  
\eea
satisfied by the matrix elements of $x\in G^{u,v}$.
In particular, our computation of~$x'_{44}$ used (\ref{eq:gl4rel})
in conjunction with Gr\"obner bases techniques
(see, e.g., \cite{cox-little-oshea}). 
}
\end{example}

\subsection{Double pseudoline arrangements}

As an essential new ingredient in our solution to the factorization problem
for $GL_n$, we will represent a factorization scheme $\ii$ 
geometrically by the corresponding \emph{double pseudoline arrangement} 
(or double wiring diagram). 
This arrangement is obtained by superimposing two arrangements
naturally associated to the $E$- and $F$-part of $\ii$
(cf.~\cite{BFZ}). 

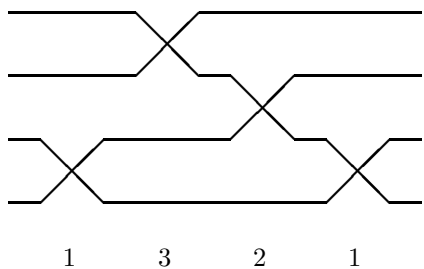
\begin{figure}[ht]
\setlength{\unitlength}{1.2pt} 
\begin{center}
\begin{picture}(130,85)(0,-20)
\thicklines
  \put(0,0){\line(1,0){10}}
  \put(30,0){\line(1,0){70}}
  \put(120,0){\line(1,0){10}}
  \put(0,20){\line(1,0){10}}
  \put(30,20){\line(1,0){40}}
  \put(90,20){\line(1,0){10}}
  \put(120,20){\line(1,0){10}}
  \put(0,40){\line(1,0){40}}
  \put(60,40){\line(1,0){10}}
  \put(90,40){\line(1,0){40}}
  \put(0,60){\line(1,0){40}}
  \put(60,60){\line(1,0){70}}

  \put(10,0){\line(1,1){20}}
  \put(40,40){\line(1,1){20}}
  \put(70,20){\line(1,1){20}}
  \put(100,0){\line(1,1){20}}

  \put(10,20){\line(1,-1){20}}
  \put(40,60){\line(1,-1){20}}
  \put(70,40){\line(1,-1){20}}
  \put(100,20){\line(1,-1){20}}

  \put(17,-20){1}
  \put(47,-20){3}
  \put(77,-20){2}
  \put(107,-20){1}

\end{picture}
\end{center}
\caption{Pseudoline arrangement for the reduced word 1321}
\label{fig:pseudo}
\end{figure}

To be self-contained, let us recall the definition
of a pseudoline arrangement associated to a reduced word.
This is best done by an example.
Consider $v=4213\in S_4$, together with the reduced decomposition
$v=s_1 s_3 s_2 s_1$ (cf.~(\ref{eq:4312-4213})).
The corresponding pseudoline arrangement is given 
in Figure~\ref{fig:pseudo}; to each entry $i$ of $\ii$, 
we associate a crossing at the $i$th level, 
counting from the bottom.

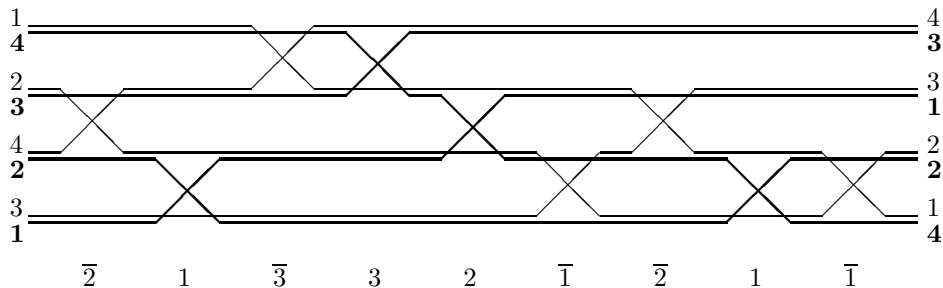
\begin{figure}[ht]
\setlength{\unitlength}{1.2pt} 
\begin{center}
\begin{picture}(280,90)(0,-20)
\thicklines

  \put(0,0){\line(1,0){40}}
  \put(60,0){\line(1,0){160}}
  \put(240,0){\line(1,0){40}}
  \put(0,20){\line(1,0){40}}
  \put(60,20){\line(1,0){70}}
  \put(150,20){\line(1,0){70}}
  \put(240,20){\line(1,0){40}}
  \put(0,40){\line(1,0){100}}
  \put(120,40){\line(1,0){10}}
  \put(150,40){\line(1,0){130}}
  \put(0,60){\line(1,0){100}}
  \put(120,60){\line(1,0){160}}

  \put(40,0){\line(1,1){20}}
  \put(100,40){\line(1,1){20}}
  \put(130,20){\line(1,1){20}}
  \put(220,0){\line(1,1){20}}

  \put(40,20){\line(1,-1){20}}
  \put(100,60){\line(1,-1){20}}
  \put(130,40){\line(1,-1){20}}
  \put(220,20){\line(1,-1){20}}

  \put(47,-20){1}
  \put(107,-20){3}
  \put(137,-20){2}
  \put(227,-20){1}

  \put(283,-6){$\mathbf{4}$}
  \put(283,14){$\mathbf{2}$}
  \put(283,34){$\mathbf{1}$}
  \put(283,54){$\mathbf{3}$}
  \put(-6,-6){$\mathbf{1}$}
  \put(-6,14){$\mathbf{2}$}
  \put(-6,34){$\mathbf{3}$}
  \put(-6,54){$\mathbf{4}$}

\light{
\thinlines
  \put(0,2){\line(1,0){160}}
  \put(180,2){\line(1,0){70}}
  \put(270,2){\line(1,0){10}}
  \put(0,22){\line(1,0){10}}
  \put(30,22){\line(1,0){130}}
  \put(180,22){\line(1,0){10}}
  \put(210,22){\line(1,0){40}}
  \put(270,22){\line(1,0){10}}
  \put(0,42){\line(1,0){10}}
  \put(30,42){\line(1,0){40}}
  \put(90,42){\line(1,0){100}}
  \put(210,42){\line(1,0){70}}
  \put(0,62){\line(1,0){70}}
  \put(90,62){\line(1,0){190}}

  \put(10,22){\line(1,1){20}}
  \put(70,42){\line(1,1){20}}
  \put(160,2){\line(1,1){20}}
  \put(190,22){\line(1,1){20}}
  \put(250,2){\line(1,1){20}}

  \put(10,42){\line(1,-1){20}}
  \put(70,62){\line(1,-1){20}}
  \put(160,22){\line(1,-1){20}}
  \put(190,42){\line(1,-1){20}}
  \put(250,22){\line(1,-1){20}}

  \put(17,-20){$\overline 2$}
  \put(77,-20){$\overline 3$}
  \put(167,-20){$\overline 1$}
  \put(197,-20){$\overline 2$}
  \put(257,-20){$\overline 1$}

  \put(283,2){${1}$}
  \put(283,22){${2}$}
  \put(283,42){${3}$}
  \put(283,62){${4}$}
  \put(-6,2){${3}$}
  \put(-6,22){${4}$}
  \put(-6,42){${2}$}
  \put(-6,62){${1}$}
}

\end{picture}
\end{center}
\caption{Double pseudoline arrangement}
\label{fig:double-pseudo}
\end{figure}

Let us now consider the factorization scheme~$\ii$ defined 
by~(\ref{eq:example-fact-scheme}).
The $E$-part of $\ii$ is 1321, and we already drew the corresponding
arrangement.
The $F$-part of $\ii$ is 23121. 
To construct the double pseudoline arrangement for~$\ii$,
we superimpose the arrangements for 1321 and 23121,
aligning them closely in the vertical direction,
and placing the intersections so that tracing them 
left-to-right would produce the same shuffle of the two reduced words
that appears in~$\ii$.
This results in the double pseudoline  arrangement 
in Figure~\ref{fig:double-pseudo}.

The double pseudoline arrangement that corresponds to 
a factorization scheme~$\ii$ is denoted by $\Arr(\ii)$.
The two subarrangements of~$\Arr(\ii)$ corresponding 
to the reduced words for $v$ and $u$ are called 
the $E$- and $F$-part of $\Arr(\ii)$, 
and their crossing points are referred to as 
$E$- and $F$-crossings, respectively.
These crossings are in an obvious bijection with
the non-circled entries of~$\ii$.

We next label the pseudolines of $\Arr(\ii)$
using the following important convention. 
The pseudolines of the $F$-part of $\Arr(\ii)$
are labelled 1~through~$n$ bottom-up at the right end
of the arrangement (just as in~\cite{BFZ}).
At the same time, the pseudolines of the $E$-part
are labelled bottom-up at the left end. 
See Figure~\ref{fig:double-pseudo}.

Another numbering that we are going to use 
is the bottom-to-top numbering of the $n-1$ horizontal strips
containing the crossings of the arrangement.
We say that the strip between the $j$th and $(j+1)$st horizontal lines,
counting from the bottom, has \emph{level}~$j$, 
and all the $E$- and $F$-crossings contained in this strip 
are of level~$j$.

\begin{figure}[ht]
\setlength{\unitlength}{1.28pt} 
\begin{center}
\begin{picture}(280,90)(0,-20)
\thicklines

  \put(0,0){\line(1,0){40}}
  \put(60,0){\line(1,0){160}}
  \put(240,0){\line(1,0){40}}
  \put(0,20){\line(1,0){40}}
  \put(60,20){\line(1,0){70}}
  \put(150,20){\line(1,0){70}}
  \put(240,20){\line(1,0){40}}
  \put(0,40){\line(1,0){100}}
  \put(120,40){\line(1,0){10}}
  \put(150,40){\line(1,0){130}}
  \put(0,60){\line(1,0){100}}
  \put(120,60){\line(1,0){160}}

  \put(40,0){\line(1,1){20}}
  \put(100,40){\line(1,1){20}}
  \put(130,20){\line(1,1){20}}
  \put(220,0){\line(1,1){20}}

  \put(40,20){\line(1,-1){20}}
  \put(100,60){\line(1,-1){20}}
  \put(130,40){\line(1,-1){20}}
  \put(220,20){\line(1,-1){20}}

  \put(47,-20){1}
  \put(107,-20){3}
  \put(137,-20){2}
  \put(227,-20){1}

\light{
\thinlines
  \put(0,2){\line(1,0){160}}
  \put(180,2){\line(1,0){70}}
  \put(270,2){\line(1,0){10}}
  \put(0,22){\line(1,0){10}}
  \put(30,22){\line(1,0){130}}
  \put(180,22){\line(1,0){10}}
  \put(210,22){\line(1,0){40}}
  \put(270,22){\line(1,0){10}}
  \put(0,42){\line(1,0){10}}
  \put(30,42){\line(1,0){40}}
  \put(90,42){\line(1,0){100}}
  \put(210,42){\line(1,0){70}}
  \put(0,62){\line(1,0){70}}
  \put(90,62){\line(1,0){190}}

  \put(10,22){\line(1,1){20}}
  \put(70,42){\line(1,1){20}}
  \put(160,2){\line(1,1){20}}
  \put(190,22){\line(1,1){20}}
  \put(250,2){\line(1,1){20}}

  \put(10,42){\line(1,-1){20}}
  \put(70,62){\line(1,-1){20}}
  \put(160,22){\line(1,-1){20}}
  \put(190,42){\line(1,-1){20}}
  \put(250,22){\line(1,-1){20}}

  \put(17,-20){$\overline 2$}
  \put(77,-20){$\overline 3$}
  \put(167,-20){$\overline 1$}
  \put(197,-20){$\overline 2$}
  \put(257,-20){$\overline 1$}

}
\bulletcolor{
   \put(62,-20){$\!\!\circled{\,3}$}
  \put(182,-20){$\!\!\circled{\,1}$}
  \put(237,-20){$\!\!\circled{\,4}$}
  \put(246,-20){$\!\!\circled{\,2}$}

  \put(65,41){\circle*{4}}
  \put(185,1){\circle*{4}}
  \put(242,61){\circle*{4}}
  \put(248,21){\circle*{4}}
}

\end{picture}
\end{center}
\caption{Arrangement $\Arr_\bullet(\ii)$}
\label{fig:rigged}
\end{figure}
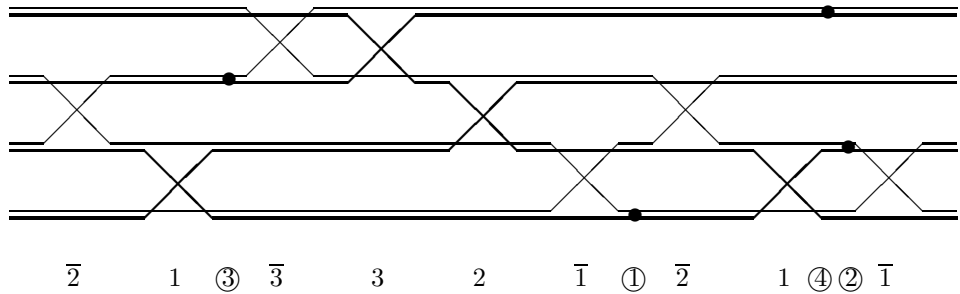

Note that arrangement $\Arr(\ii)$ does not depend on the $H$-part of
the factorization scheme~$\ii$.
In order to include the $H$-part into the picture, we 
associate with each entry~$\circled{\,j}$ a bullet~$\,\bullet\,$ 
placed on the $j$th horizontal line.
The position of a bullet corresponds to the position
of~$\circled{\,j}$ in~$\ii$,
so that when the arrangement is traced left-to-right,
the crossings and bullets appear in the same order 
as the entries of~$\ii$ that they represent.
The resulting ``rigged'' arrangement is denoted by~$\Arr_\bullet(\ii)$.
Figure~\ref{fig:rigged}
shows $\Arr_\bullet(\ii)$ for the factorization
scheme~(\ref{eq:example-fact-scheme}).

To make our terminology uniform, we will refer to the bullets
in~$\Arr_\bullet(\ii)$ as $H$-\emph{crossings} 
(despite the fact that they are not crossings geometrically). 
Thus the total number of crossings in $\Arr_\bullet(\ii)$ 
is~$l = n + l(u)+l(v)$, and they are associated with 
the variables~$t_k$ in the factorization~(\ref{eq:t(x)}).
We will occasionally refer to the crossing in $\Arr_\bullet(\ii)$ 
associated with a factorization parameter~$t_k\,$ 
by simply saying ``crossing~$t_k\,$.'' 
The $H$-crossing lying on the $i$th horizontal line 
will also be denoted by~$d_i\,$. 

\subsection{Solution to the factorization problem}

Let us fix permutations $u,v\in S_n$
and a factorization scheme $\ii$ of type~$(u,v)$;
in this section,  we present our solution to the corresponding
factorization problem. 
As in \cite{BFZ}, the combinatorics needed to formulate the answer
involves not only the crossings of the arrangement $\Arr(\ii)$ 
but also its \emph{chambers}, which can be defined as horizontal
segments between consecutive crossings of the same level. 
More precisely, each horizontal strip with, say,
$k$ crossings breaks down into $k+1$ chambers
(including the ones at the ends of the strip).
Two more chambers are located at the bottom and the top
of the arrangement.
To illustrate, the arrangement in Figure~\ref{fig:double-pseudo}
has 14 chambers; in general, there are $l+1$ of them. 

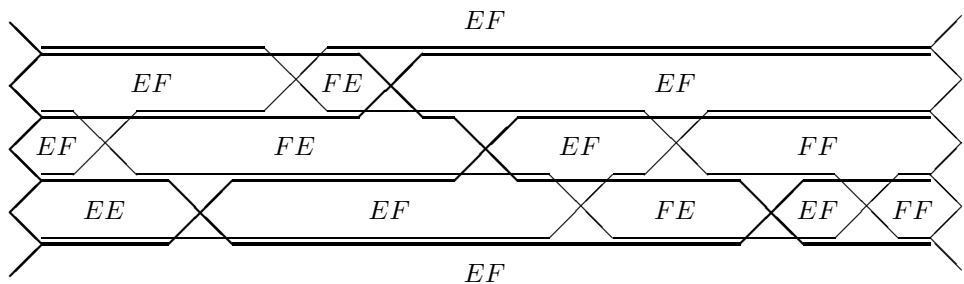
\begin{figure}[t]
\setlength{\unitlength}{1.2pt} 
\begin{center}
\begin{picture}(280,90)(0,-10)
\thicklines

  \put(0,0){\line(1,0){40}}
  \put(60,0){\line(1,0){160}}
  \put(240,0){\line(1,0){40}}
  \put(0,20){\line(1,0){40}}
  \put(60,20){\line(1,0){70}}
  \put(150,20){\line(1,0){70}}
  \put(240,20){\line(1,0){40}}
  \put(0,40){\line(1,0){100}}
  \put(120,40){\line(1,0){10}}
  \put(150,40){\line(1,0){130}}
  \put(0,60){\line(1,0){100}}
  \put(120,60){\line(1,0){160}}

  \put(40,0){\line(1,1){20}}
  \put(100,40){\line(1,1){20}}
  \put(130,20){\line(1,1){20}}
  \put(220,0){\line(1,1){20}}

  \put(40,20){\line(1,-1){20}}
  \put(100,60){\line(1,-1){20}}
  \put(130,40){\line(1,-1){20}}
  \put(220,20){\line(1,-1){20}}

  \put(-10,10){\line(1,1){10}}
  \put(-10,10){\line(1,-1){10}}
  \put(-10,-10){\line(1,1){10}}
  \put(-10,30){\line(1,1){10}}
  \put(-10,30){\line(1,-1){10}}
  \put(-10,50){\line(1,1){10}}
  \put(-10,50){\line(1,-1){10}}
  \put(-10,70){\line(1,-1){10}}

\light{
\thinlines
  \put(0,2){\line(1,0){160}}
  \put(180,2){\line(1,0){70}}
  \put(270,2){\line(1,0){10}}
  \put(0,22){\line(1,0){10}}
  \put(30,22){\line(1,0){130}}
  \put(180,22){\line(1,0){10}}
  \put(210,22){\line(1,0){40}}
  \put(270,22){\line(1,0){10}}
  \put(0,42){\line(1,0){10}}
  \put(30,42){\line(1,0){40}}
  \put(90,42){\line(1,0){100}}
  \put(210,42){\line(1,0){70}}
  \put(0,62){\line(1,0){70}}
  \put(90,62){\line(1,0){190}}

  \put(10,22){\line(1,1){20}}
  \put(70,42){\line(1,1){20}}
  \put(160,2){\line(1,1){20}}
  \put(190,22){\line(1,1){20}}
  \put(250,2){\line(1,1){20}}

  \put(10,42){\line(1,-1){20}}
  \put(70,62){\line(1,-1){20}}
  \put(160,22){\line(1,-1){20}}
  \put(190,42){\line(1,-1){20}}
  \put(250,22){\line(1,-1){20}}

  \put(280,22){\line(1,1){10}}
  \put(280,22){\line(1,-1){10}}
  \put(280,2){\line(1,1){10}}
  \put(280,2){\line(1,-1){10}}
  \put(280,42){\line(1,1){10}}
  \put(280,42){\line(1,-1){10}}
  \put(280,62){\line(1,1){10}}
  \put(280,62){\line(1,-1){10}}
}

  \put(133,-12){$EF$}

  \put(13,8){$EE$}
  \put(103,8){$EF$}
  \put(193,8){$FE$}
  \put(238,8){$EF$}
  \put(268,8){$FF$}

  \put(-2,28){$EF$}
  \put(73,28){$FE$}
  \put(163,28){$EF$}
  \put(238,28){$FF$}

  \put(28,48){$EF$}
  \put(88,48){$FE$}
  \put(193,48){$EF$}

  \put(133,68){$EF$}

\end{picture}
\end{center}
\caption{Types of chambers}
\label{fig:EF-chambers}
\end{figure}

We say that a chamber~$C$ is of \emph{type}~$EF$ if the left endpoint
of~$C$ is an $E$-crossing, 
while its right endpoint is an $F$-crossing. 
Chambers of types $EE$, $FE$ and $FF$ are defined in a similar way.
Figure~\ref{fig:EF-chambers}
shows the types of all 14 chambers of 
the arrangement in Figure~\ref{fig:double-pseudo}. 
Here and in the sequel, we use the following important convention:
on each level, there is a fictitious $E$-crossing 
at the left border of the arrangement, 
and a fictitious $F$-crossing at the right border.
These fictitious crossings determine
the types of the chambers adjacent to the boundary of~$\Arr(\ii)$.

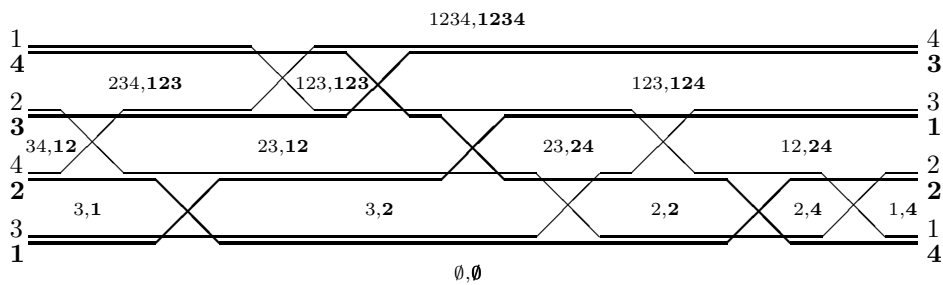
\begin{figure}[t]
\setlength{\unitlength}{1.2pt} 
\begin{center}
\begin{picture}(280,90)(0,-10)
\thicklines

  \put(0,0){\line(1,0){40}}
  \put(60,0){\line(1,0){160}}
  \put(240,0){\line(1,0){40}}
  \put(0,20){\line(1,0){40}}
  \put(60,20){\line(1,0){70}}
  \put(150,20){\line(1,0){70}}
  \put(240,20){\line(1,0){40}}
  \put(0,40){\line(1,0){100}}
  \put(120,40){\line(1,0){10}}
  \put(150,40){\line(1,0){130}}
  \put(0,60){\line(1,0){100}}
  \put(120,60){\line(1,0){160}}

  \put(40,0){\line(1,1){20}}
  \put(100,40){\line(1,1){20}}
  \put(130,20){\line(1,1){20}}
  \put(220,0){\line(1,1){20}}

  \put(40,20){\line(1,-1){20}}
  \put(100,60){\line(1,-1){20}}
  \put(130,40){\line(1,-1){20}}
  \put(220,20){\line(1,-1){20}}

  \put(283,-6){$\mathbf{4}$}
  \put(283,14){$\mathbf{2}$}
  \put(283,34){$\mathbf{1}$}
  \put(283,54){$\mathbf{3}$}
  \put(-6,-6){$\mathbf{1}$}
  \put(-6,14){$\mathbf{2}$}
  \put(-6,34){$\mathbf{3}$}
  \put(-6,54){$\mathbf{4}$}

\light{
\thinlines

  \put(0,2){\line(1,0){160}}
  \put(180,2){\line(1,0){70}}
  \put(270,2){\line(1,0){10}}
  \put(0,22){\line(1,0){10}}
  \put(30,22){\line(1,0){130}}
  \put(180,22){\line(1,0){10}}
  \put(210,22){\line(1,0){40}}
  \put(270,22){\line(1,0){10}}
  \put(0,42){\line(1,0){10}}
  \put(30,42){\line(1,0){40}}
  \put(90,42){\line(1,0){100}}
  \put(210,42){\line(1,0){70}}
  \put(0,62){\line(1,0){70}}
  \put(90,62){\line(1,0){190}}

  \put(10,22){\line(1,1){20}}
  \put(70,42){\line(1,1){20}}
  \put(160,2){\line(1,1){20}}
  \put(190,22){\line(1,1){20}}
  \put(250,2){\line(1,1){20}}

  \put(10,42){\line(1,-1){20}}
  \put(70,62){\line(1,-1){20}}
  \put(160,22){\line(1,-1){20}}
  \put(190,42){\line(1,-1){20}}
  \put(250,22){\line(1,-1){20}}

  \put(283,2){${1}$}
  \put(283,22){${2}$}
  \put(283,42){${3}$}
  \put(283,62){${4}$}
  \put(-6,2){${3}$}
  \put(-6,22){${4}$}
  \put(-6,42){${2}$}
  \put(-6,62){${1}$}
 
}

  \put(134,-10){$_{\light{\emptyset},\emptyset\hspace{-.057in}
\emptyset\hspace{-.057in}\emptyset}$}

  \put(14,10){$_{\light{3},\mathbf{1}}$}
  \put(106,10){$_{\light{3},\mathbf{2}}$}
  \put(196,10){$_{\light{2},\mathbf{2}}$}
  \put(241,10){$_{\light{2},\mathbf{4}}$}
  \put(271,10){$_{\light{1},\mathbf{4}}$}

  \put(-1,30){$_{\light{34},\mathbf{12}}$}
  \put(72,30){$_{\light{23},\mathbf{12}}$}
  \put(162,30){$_{\light{23},\mathbf{24}}$}
  \put(237,30){$_{\light{12},\mathbf{24}}$}

  \put(25,50){$_{\light{234},\mathbf{123}}$}
  \put(84,50){$_{\light{123},\mathbf{123}}$}
  \put(190,50){$_{\light{123},\mathbf{124}}$}

  \put(126,70){$_{\light{1234},\mathbf{1234}}$}

\end{picture}
\end{center}
\caption{Chamber sets}
\label{fig:chamber-sets}
\end{figure}

For every chamber $C$ in $\Arr(\ii)$, let $I(C)$ denote the set 
of labels of the lines of the $F$-part of the arrangement
that pass below~$C$. 
Analogously, $J(C)$ will consist of the labels of lines of the $F$-part 
of $\Arr(\ii)$ that pass below~$C$. 
The sets $I(C)$ and $J(C)$  are called \emph{chamber sets}
for the factorization scheme~$\ii$. 
Figure~\ref{fig:chamber-sets} shows the chamber sets $I(C)$ and $J(C)$
for each chamber of the given arrangement.
Note that if $C$ is a chamber of level~$i$, 
then both $I(C)$ ans $J(C)$ have $i$ elements. 

Our constructions will also involve the ``big'' chambers
formed by the $E$-part and the $F$-part of 
a double pseudoline arrangement,
taken separately. We will refer to these ``big'' chambers
as $E$-chambers and $F$-chambers, respectively.
For example, the arrangement in Figure~\ref{fig:double-pseudo}
has 9 $E$-chambers, which are in obvious bijection
with the 9~chambers in Figure~\ref{fig:pseudo}.

For every chamber~$C$ of the arrangement~$\Arr(\ii)$, we denote
$$
M_C = \Delta_{I(C),J(C)} \ ;
$$
this minor is considered as a regular function on $G$ (with the convention that 
$\Delta_{\emptyset, \emptyset} = 1$). 
For example, if $C$ is the rightmost chamber of level~2
in Figure~\ref{fig:chamber-sets}, 
then $M_C=\Delta_{12,24}\,$.

To each $i = 1, \ldots, n$ we associate a rational function
on $G^{u,v}$ given by
\begin{equation}
\label{eq:Pi}
\Pi_i = \frac{\prod_C M_C}{\prod_{C'} M_{C'}} \ , 
\end{equation}
where $C$ runs over all chambers of level~$i$ and type~$FE$,
while $C'$ runs over all chambers of level~$i$ and type~$EF$.
For example, in Figure~\ref{fig:chamber-sets}
we have
$$
\Pi_2=\frac{\Delta_{23,12}}{\Delta_{34,12}\,\Delta_{23,24}}\ .
$$
Also, by convention, $\Pi_0 = 1$. 

Let $\C$ be a ``big'' $K$-chamber of level~$i$, where $K$
is one of the symbols $E$ and~$F$. 
Let $L$ be the other of these symbols
(i.e., $L=F$ if $K=E$, and $L=E$ if $K=F$). 
We define
\begin{equation}
\label{eq:M-right}
\M^{\rm right}_{\C} =  
\displaystyle \frac{\tilde M \prod_{C'} M_{C'}}{\prod_{C''} M_{C''}} \ ,
\end{equation}
where 
\begin{itemize}
\item 
$C'$ runs over all chambers of level~$i$ and type~$LK$
to the right of~$\C$;
\item
$C''$ runs over all chambers of level~$i$ and type~$KL$
to the right of~$\C$;
\item 
$\tilde M=M_{C}$, where $C$ is the (``small'') chamber at the
right end of~$\C$ (inside~$\C$),
unless $K=E$ and $C$ is stuck to the right border, in which case
$\tilde M=1$. 
\end{itemize}
Analogously,
\begin{equation}
\label{eq:M-left}
\M^{\rm left}_{\C} =  
\displaystyle \frac{\tilde M \prod_{C'} M_{C'}}{\prod_{C''} M_{C''}} \ ,
\end{equation}
where  
\begin{itemize}
\item 
$C'$ runs over all chambers of level~$i$ and type~$KL$
to the left of~$\C$;
\item
$C''$ runs over all chambers of level~$i$ and type~$LK$
to the left of~$\C$;
\item
$\tilde M=M_C$, where $C$ is the (``small'') chamber at the 
left end of~$\C$ (inside~$\C$),
unless $K=F$ and $C$ is stuck to the left border, in which case
$\tilde M=1$ .
\end{itemize}

We are finally prepared to state our solution to the factorization
problem.

\begin{theorem}
\label{th:main-gln} 
Let~$\ii$ be a factorization scheme of type $(u,v)$, and suppose
a matrix $x \in G^{u,v}$ admits the factorization 
$x = x_\ii (t_1, \dots, t_l)$ with all $t_k$ nonzero complex numbers.  
Let $x' = \zeta^{u, v}(x) \in G^{u^{-1}, v^{-1}}$ denote the matrix 
obtained from $x$ by the ``twist'' {\rm (\ref{eq:gln-zeta})}.
Then the factorization parameters $t_k$ are determined as follows.
\begin{itemize} 
\item If $t_k$ corresponds to the $H$-crossing $d_i\,$, then 
\begin{equation}
\label{eq:ansatz-h}
t_k = \displaystyle\frac{\Pi_i (x')}{\Pi_{i-1} (x')}\ ,
\end{equation} 
where $\Pi_i$ and $\Pi_{i-1}$ are given by {\rm (\ref{eq:Pi})}.
\item Let $t_k$ correspond to an $E$- or $F$-crossing of level $i$, and 
let $\A, \B, \C, \D$ be the four ``big'' chambers surrounding this 
crossing, as shown:
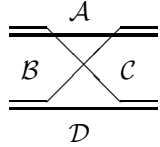
\setlength{\unitlength}{1.4pt} 
\begin{center}
\begin{picture}(40,35)(0,-5)
\thicklines

  \put(0,0){\line(1,0){40}}
  \put(0,20){\line(1,0){40}}

\thinlines
  \put(0,2){\line(1,0){10}}
  \put(30,2){\line(1,0){10}}
  \put(0,22){\line(1,0){10}}
  \put(30,22){\line(1,0){10}}

  \put(10,2){\line(1,1){20}}
  \put(10,22){\line(1,-1){20}}


  \put(16,24){$\A$}
  \put(3,8){$\B$}
  \put(30,8){$\C$}
  \put(16,-9){$\D$}

\end{picture}
\end{center}
Then
\beal
\label{eq:ansatz-ef}
t_k = \displaystyle\frac{\M^{{\rm opp}(d_{i+1})}_{\A}(x')\,
\M^{{\rm opp}(d_{i})}_{\D} (x')}
{\M^{{\rm opp}(d_{i+1})}_{\B}(x') \,\M^{{\rm opp}(d_{i})}_{\C}(x')} \ ,
\eea
where we refer to the notation of {\rm
  (\ref{eq:M-right})--(\ref{eq:M-left})} as follows: 
the superscript ${\rm opp}(d_{i})$ stands for ``left" if 
the $H$-crossing $d_i$ is to the right of~$t_k\,$, and for ``right" if 
$d_i$ is to the left of~$t_k\,$.
\end{itemize}
\end{theorem}

Theorem~\ref{th:main-gln} is obtained as a specialization of 
Theorem~\ref{th:t-through-x}, with the help of the following 
additional commutation relations:
\begin{eqnarray}
\begin{array}{rcl}
\label{eq:commut-rel-2}
x_i(a) x_{\circj}(b) &=& x_{\circj}(b) x_i(a) \ ,\ j\notin \{i,i+1\};
\\[.1in]
x_i(a) x_\circi(b) &=& x_\circi(b) x_i(a/b) \ ;\\[.1in]
x_i(a) x_\circj(b) &=& x_\circj(b) x_i(ab) \ ,\ j=i+1\ ;\\[.1in]
x_{\bar i}(a) x_{\circj}(b) &=& x_{\circj}(b) x_{\bar i} (a) 
\ ,\ j\notin \{i,i+1\}; \\[.1in]
x_{\bar i}(a) x_\circi(b) &=& x_\circi(b) x_{\bar i}(ab) \ ;\\[.1in]
x_{\bar i}(a) x_\circj(b) &=& x_\circj(b) x_{\bar i}(a/b)\ ,\ j=i+1   \ .
\end{array}
\end{eqnarray}

%

\begin{example}\label{example:t9}
{\rm
To illustrate Theorem~\ref{th:main-gln}, let us compute 
the factorization parameter~$t_9$ corresponding 
to the rightmost $F$-crossing of level~2 
in Figure~\ref{fig:rigged}. 
It is given by
\[
t_9 = 
\displaystyle\frac{\M^{\rm right}_{\A}(x')\,
\M^{\rm left}_{\D}(x')}
{\M^{\rm right}_{\B}(x') \,\M^{\rm left}_{\C}(x')} \ ,
\]
where
\begin{eqnarray*}
\begin{array}{rcl}
&\M^{\rm right}(\A)=\Delta_{123,124},\\
\M^{\rm right}(\B)=\Delta_{23,24}\,,
& &
\M^{\rm left}(\C) 
   =\displaystyle\frac{\Delta_{12,24}\,\Delta_{23,12}}{\Delta_{23,24}\,
   \Delta_{34,12}}
     \,,\\
&\ \ \M^{\rm left}(\D) 
   =\displaystyle\frac{\Delta_{2,2}}{\Delta_{3,2}}\,.
\end{array}
\end{eqnarray*}
Hence
\beac
\label{eq:t9}
 t_9=\displaystyle\frac{\Delta_{2,2}\,\Delta_{34,12}\,\Delta_{123,124}}
                       {\Delta_{3,2}\,\Delta_{12,24}\,\Delta_{23,12}} 
                       (x')\ .
\eea
Substituting the twisted matrix~$x'$ from
Example~\ref{example:running-twist} into (\ref{eq:t9})
and simplifying, we finally obtain
\begin{equation}
\label{eq:t9final}
t_9=
\displaystyle
\frac{\Delta_{23,12}(x)\,(x_{43}\Delta_{12,12}(x)-\Delta_{124,123}(x))}
     {x_{23}\,\Delta_{24,12}(x)\,\Delta_{123,123}(x)}\ .
\end{equation}
This answer can be verified directly using 
the combinatorial interpretation of minors~$\Delta_{I,J}$ 
in terms of planar networks (see Proposition~\ref{pr:planar-minors}). From 
Figure~\ref{fig:planar-network} one
obtains: $\Delta_{23,12}=t_3t_7t_8t_9t_{12}\,$,
$x_{43}=t_4\,$,
$\Delta_{12,12}=t_8t_{12}(1+t_6t_9)$,
$\Delta_{124,123}=t_4t_8t_{12}\,$,
$x_{23}=t_6\,$,
$\Delta_{24,12}=t_4t_7t_8t_9t_{12}\,$,
$\Delta_{123,123}=t_3t_8t_{12}\,$, 
implying~(\ref{eq:t9final}).
}
\end{example}

As in Theorem~\ref{th:inverse-monomial}, formulas (\ref{eq:ansatz-h})
and (\ref{eq:ansatz-ef}) imply that the factorization 
parameters $t_1,\dots,t_l$ are related by 
an invertible monomial transformation to the $l$ minors 
$\Delta_{I(C), J(C)} (x')$ of the twisted matrix $x'$
that correspond to the chambers of 
the arrangement~$\Arr(\ii)$, with the bottom chamber excluded.
The inverse transformation has the following description
which can be deduced from (\ref{eq:inverse-monomial})
(since we left the latter formula without proof,
the same is true for our next theorem, although it is not hard to give
it a direct proof).

\begin{theorem}
\label{th:inverse-ansatz}
Formulas~{\rm (\ref{eq:ansatz-h})--(\ref{eq:ansatz-ef})}
are equivalent to the following formulas:
\beac
\Delta_{I(C),J(C)} (x') =\left(\displaystyle\prod t_k\right)^{-1}\ ,
\eea
where the product is over all~$t_k$ which correspond to the following
types of crossings:

\noindent
$\bullet$~
$E$-crossings to the right of~$C$ such that $C$ lies 
between the lines intersecting
at~$t_k$;

\noindent
$\bullet$~
$F$-crossings to the left  of~$C$ such that $C$ lies 
between the lines intersecting
at~$t_k$;

\noindent
$\bullet$~
$H$-crossings to the right of~$C$ such that $C$ lies above the $E$-line
passing through~$t_k$;

\noindent
$\bullet$~
$H$-crossings to the left  of~$C$ such that $C$ lies above the $F$-line
passing through~$t_k$.
\end{theorem}

For example, in  Figure~\ref{fig:chamber-sets}, 
$\Delta_{3,1} (x') =(t_2 t_6)^{-1}\,$, 
$\Delta_{123,124} (x') =(t_1t_4t_8t_{12})^{-1}\,$,
etc.

\subsection{Applications to total positivity}
\label{sec:applications-to-TP}

In the case $G = GL_n (\CC)$, the definition of the totally 
nonnegative variety $G_{\geq 0}$ given in Section~\ref{sec:TPgeneral}
is modified as follows: $G_{\geq 0}$ is the multiplicative semigroup
generated by elementary Jacobi matrices 
(cf.~(\ref{eq:jacobi})--(\ref{eq:diag-jacobi})) 
$x_i(t), \, x_{\bar i}(t)$, and $x_\circi(t)$ with $t> 0$.
It is known \cite{loewner,whitney}
that this definition of total nonnegativity is equivalent to the
classical one: an invertible matrix $x \in G$ belongs to $G_{\geq 0}$
if and only if all minors $\Delta_{I,J}(x)$ (in particular, 
all matrix entries) are nonnegative.

For any two permutations $u$ and $v$ in $S_n\,$, the corresponding 
\emph{totally positive variety} is defined by
$$G^{u,v}_{> 0} = G_{\geq 0} \cap G^{u,v} \ .$$
Each factorization scheme of type $(u,v)$ gives rise
to a parametrization of $G^{u,v}_{> 0}$, according to the 
following analogue of Theorem~\ref{th:LusztigTP}.

\begin{theorem}
\label{th:LusztigTPGL}
For any factorization scheme $\ii= (i_1, \ldots, i_{l})$  of type $(u,v)$,
the corresponding product map 
$x_\ii: \CC_{\neq 0}^l \to G^{u,v}$ given by {\rm (\ref{eq:t(x)})} 
restricts to a bijection between $\RR_{>0}^l$ and~$G^{u,v}_{>0}\,$. 
\end{theorem}

The twist map $\zeta^{u,v}$ defined by (\ref{eq:gln-zeta}) respects 
total positivity: it sends totally nonnegative matrices in 
$G^{u,v}$ to totally nonnegative matrices in $G^{u^{-1}, v^{-1}}$
(cf.\ Theorem~\ref{th:zeta-positivity}).
Combining this fact with Theorems~\ref{th:LusztigTPGL} and 
\ref{th:main-gln} leads to a family of criteria for total positivity,
one for each factorization scheme.

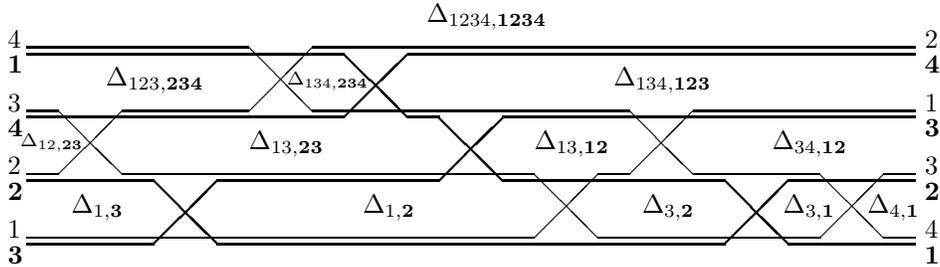
\begin{figure}[ht]
\setlength{\unitlength}{1.2pt} 
\begin{center}
\begin{picture}(280,70)(0,0)
\thicklines

  \put(0,0){\line(1,0){40}}
  \put(60,0){\line(1,0){160}}
  \put(240,0){\line(1,0){40}}
  \put(0,20){\line(1,0){40}}
  \put(60,20){\line(1,0){70}}
  \put(150,20){\line(1,0){70}}
  \put(240,20){\line(1,0){40}}
  \put(0,40){\line(1,0){100}}
  \put(120,40){\line(1,0){10}}
  \put(150,40){\line(1,0){130}}
  \put(0,60){\line(1,0){100}}
  \put(120,60){\line(1,0){160}}

  \put(40,0){\line(1,1){20}}
  \put(100,40){\line(1,1){20}}
  \put(130,20){\line(1,1){20}}
  \put(220,0){\line(1,1){20}}

  \put(40,20){\line(1,-1){20}}
  \put(100,60){\line(1,-1){20}}
  \put(130,40){\line(1,-1){20}}
  \put(220,20){\line(1,-1){20}}

  \put(283,-6){$\mathbf{1}$}
  \put(283,14){$\mathbf{2}$}
  \put(283,34){$\mathbf{3}$}
  \put(283,54){$\mathbf{4}$}
  \put(-6,-6){$\mathbf{3}$}
  \put(-6,14){$\mathbf{2}$}
  \put(-6,34){$\mathbf{4}$}
  \put(-6,54){$\mathbf{1}$}

\light{
\thinlines

  \put(0,2){\line(1,0){160}}
  \put(180,2){\line(1,0){70}}
  \put(270,2){\line(1,0){10}}
  \put(0,22){\line(1,0){10}}
  \put(30,22){\line(1,0){130}}
  \put(180,22){\line(1,0){10}}
  \put(210,22){\line(1,0){40}}
  \put(270,22){\line(1,0){10}}
  \put(0,42){\line(1,0){10}}
  \put(30,42){\line(1,0){40}}
  \put(90,42){\line(1,0){100}}
  \put(210,42){\line(1,0){70}}
  \put(0,62){\line(1,0){70}}
  \put(90,62){\line(1,0){190}}

  \put(10,22){\line(1,1){20}}
  \put(70,42){\line(1,1){20}}
  \put(160,2){\line(1,1){20}}
  \put(190,22){\line(1,1){20}}
  \put(250,2){\line(1,1){20}}

  \put(10,42){\line(1,-1){20}}
  \put(70,62){\line(1,-1){20}}
  \put(160,22){\line(1,-1){20}}
  \put(190,42){\line(1,-1){20}}
  \put(250,22){\line(1,-1){20}}

  \put(283,2){${4}$}
  \put(283,22){${3}$}
  \put(283,42){${1}$}
  \put(283,62){${2}$}
  \put(-6,2){${1}$}
  \put(-6,22){${2}$}
  \put(-6,42){${3}$}
  \put(-6,62){${4}$}
 
}


  \put(14,10){$\Delta_{\light{1},\mathbf{3}}$}
  \put(106,10){$\Delta_{\light{1},\mathbf{2}}$}
  \put(194,10){$\Delta_{\light{3},\mathbf{2}}$}
  \put(238.5,10){$\Delta_{\light{3},\mathbf{1}}$}
  \put(265,10){$\Delta_{\light{4},\mathbf{1}}$}

  \put(-2,32){$_{\Delta_{\light{12},\mathbf{23}}}$}
  \put(70,30){$\Delta_{\light{13},\mathbf{23}}$}
  \put(160,30){$\Delta_{\light{13},\mathbf{12}}$}
  \put(235,30){$\Delta_{\light{34},\mathbf{12}}$}

  \put(25,50){$\Delta_{\light{123},\mathbf{234}}$}
  \put(82,51.5){$_{\Delta_{\light{134},\mathbf{234}}}$}
  \put(185,50){$\Delta_{\light{134},\mathbf{123}}$}

  \put(126,70){$\Delta_{\light{1234},\mathbf{1234}}$}

\end{picture}
\end{center}
\caption{Minors $\Delta_{u I(C), v^{-1} J(C)}(x)$ }
\label{fig:chamber-sets-modified}
\end{figure}

For a factorization scheme $\,\ii = (i_1, \ldots, i_l)$ 
of type~$(u,v)$, let $F(\ii)$ denote the collection of $l$ minors 
$\Delta_{u I(C), v^{-1} J(C)}$, where $C$ runs over all chambers
of the arrangement~$\Arr(\ii)$, excluding the bottom chamber. 
(This notation agrees with that 
of~(\ref{eq:F(i)}).). 
We note that the pair $(u I(C), v^{-1} J(C))$ will correspond 
in the same way as above to a chamber $C$ if we relabel 
the pseudolines in $\Arr(\ii)$, 
numbering the $F$-pseudolines 1~through~$n$ bottom-up at the left end,
and 
the $E$-pseudolines bottom-up at the right end. 
See Figure~\ref{fig:chamber-sets-modified}. 

Let $F(u,v)$ denote the union of the collections $F(\ii)$
for all factorization schemes $\,\ii$ of type~$(u,v)$.
The set $F(u,v)$ can be described directly in the following way.
A subset $I \subset [1,n]$ is called
a \emph{$w$-chamber set} if, together with each element $j$ it also
contains every $i$ such that $i < j$ and $w(i) < w(j)$.
In this terminology
(originally introduced in~\cite[Section~5.3]{BFZ}), $F(u,v)$
consists of all minors $\Delta_{I, J}$ such that $I$ 
is a $u^{-1}$-chamber set while $J$ is a $v$-chamber set.

The following result specializes Theorem~\ref{th:TP-criteria}.

\begin{theorem}
\label{th:GLn-criteria}
Let $x$ be a matrix in a double Bruhat cell $G^{u,v}$,
and let $\ii$ be a factorization scheme of type~$(u,v)$.
Then the following are equivalent:

\noindent {\rm(1)} $x$ is totally nonnegative;

\noindent {\rm(2)} $\Delta_{I, J} (x) > 0$ for any $u^{-1}$-chamber set $I$ 
and any $v$-chamber set $J$;

\noindent {\rm(3)} $\Delta_{u I(C), v^{-1} J(C)}(x) > 0$ for any 
chamber $C$ of the arrangement~$\Arr(\ii)$ .
\end{theorem}

For instance, in our running example where $u$, $v$, and $\ii$ 
are given by (\ref{eq:4312-4213}) and (\ref{eq:example-fact-scheme}),
a matrix $x \in G^{u,v}$ is totally nonnegative if and only if 
the $13$ minors 
appearing in Figure~\ref{fig:chamber-sets-modified} 
are all positive if evaluated at~$x$.

Specializing Theorem~\ref{th:GLn-criteria} to the case 
$u=v=\wnot\,$, we see that the totally positive variety 
$G^{\wnot, \wnot}_{>0}$ is the classical variety 
$G_{>0}$ of the totally positive $n\times n$ matrices,
i.e., those matrices whose \emph{all} minors are (strictly) positive.
Condition (3) of Theorem~\ref{th:GLn-criteria} provides a family
of criteria for total positivity, each of which says that 
a matrix $x$ is totally positive if and only if some collection
of $n^2$ minors are positive at~$x$.

Different factorization schemes $\ii$ and $\ii'$ of the same type $(u,v)$ 
can have the same collections 
of chamber sets, thus leading to the same criteria for total positivity. 
We will say that $\ii$ and $\ii'$
(and the corresponding arrangements $\Arr(\ii)$ and $\Arr(\ii')$) 
are \emph{isotopic} if they can be obtained from each other by 
a sequence of the following ``trivial $2$-moves:''

\begin{eqnarray}
\label{eq:equiv1}
\begin{array}{cccl}
\cdots i\,j \cdots \ &\sim &\  \cdots j\,i \cdots \ ,& 
\ \ |i-j|\geq 2 \ ,\\[.1in]
\cdots \overline{i}\,\overline{j} \cdots 
   \ &\sim &\  \cdots \overline{j}\,\overline{i} \cdots \ ,& 
   \ \ |i-j|\geq 2 \ ,\\[.1in]
\cdots i \,\overline{j} \cdots 
   \ &\sim &\  \cdots \overline{j}\,i \cdots \ ,& \ \ i\neq j \ , \\[.1in]
\cdots \circled{\,i}\,\,j \cdots \ &\sim &\  \cdots j\,\circled{\,i} 
\cdots \  ,\\[.1in]
\cdots \circled{\,i}\,\,\overline{j} \cdots 
  \ &\sim &\  \cdots \overline{j}\,\circled{\,i} \cdots \ ,\\[.1in]
\cdots \circled{\,i}\,\circled{\,j} \cdots 
  \ &\sim &\  \cdots \circled{\,j}\,\circled{\,i} \cdots \ .
\end{array}
\end{eqnarray}

It is not hard to show that $\ii$ and $\ii'$ have the same collection 
of chamber sets $(I(C), J(C))$ if and only if they are isotopic.
Thus total positivity criteria in Theorem~\ref{th:GLn-criteria}
are in a bijection with ``isotopy types'' of arrangements of type~$(u,v)$.

The set of all isotopy types of arrangements of type~$(u,v)$ has a 
natural structure of a graph defined as follows.
We call two isotopy types \emph{adjacent} if the corresponding
collections of chamber sets are obtained from each other
by exchanging a single pair $(I(C), J(C))$ with another one.
The graph obtained this way is always connected, and its study is an
interesting combinatorial problem.
One can check that the adjacency relation in this graph
corresponds to the following \emph{$3$-moves}
and \emph{mixed $2$-moves} on double reduced words:  
\begin{eqnarray}
\begin{array}{rcl}
\cdots i\,j\,i \cdots &\leadsto& \cdots j\,i\,j\cdots \ ,\quad
|i-j|=1\ ,
\\[.1in]
\cdots \overline i\,\overline j \,\overline i\cdots &\leadsto&\cdots  
\overline j\,\overline i\,\overline j\cdots\ ,\quad |i-j|=1\ ,
\\[.1in]
\cdots i\,\,\overline i\,\cdots 
&\leadsto& \cdots 
\overline i\,\,i\cdots 
\end{array}
\end{eqnarray}
(cf.~Sections~\ref{sec:t-through-x} and~\ref{sec:ProofsTP}); 
the connectedness property follows 
from Proposition~\ref{pr:d-moves}. 

For $G = GL_2$ and $u=v=\wnot\,$, there are $2$ isotopy types.
The corresponding collections $F(\ii)$ 
are $\{x_{11}, \,x_{12}, \, x_{21}, \, \det (x)\}$ and 
$\{x_{22}, \,x_{12}, \, x_{21}, \, \det (x)\}$.

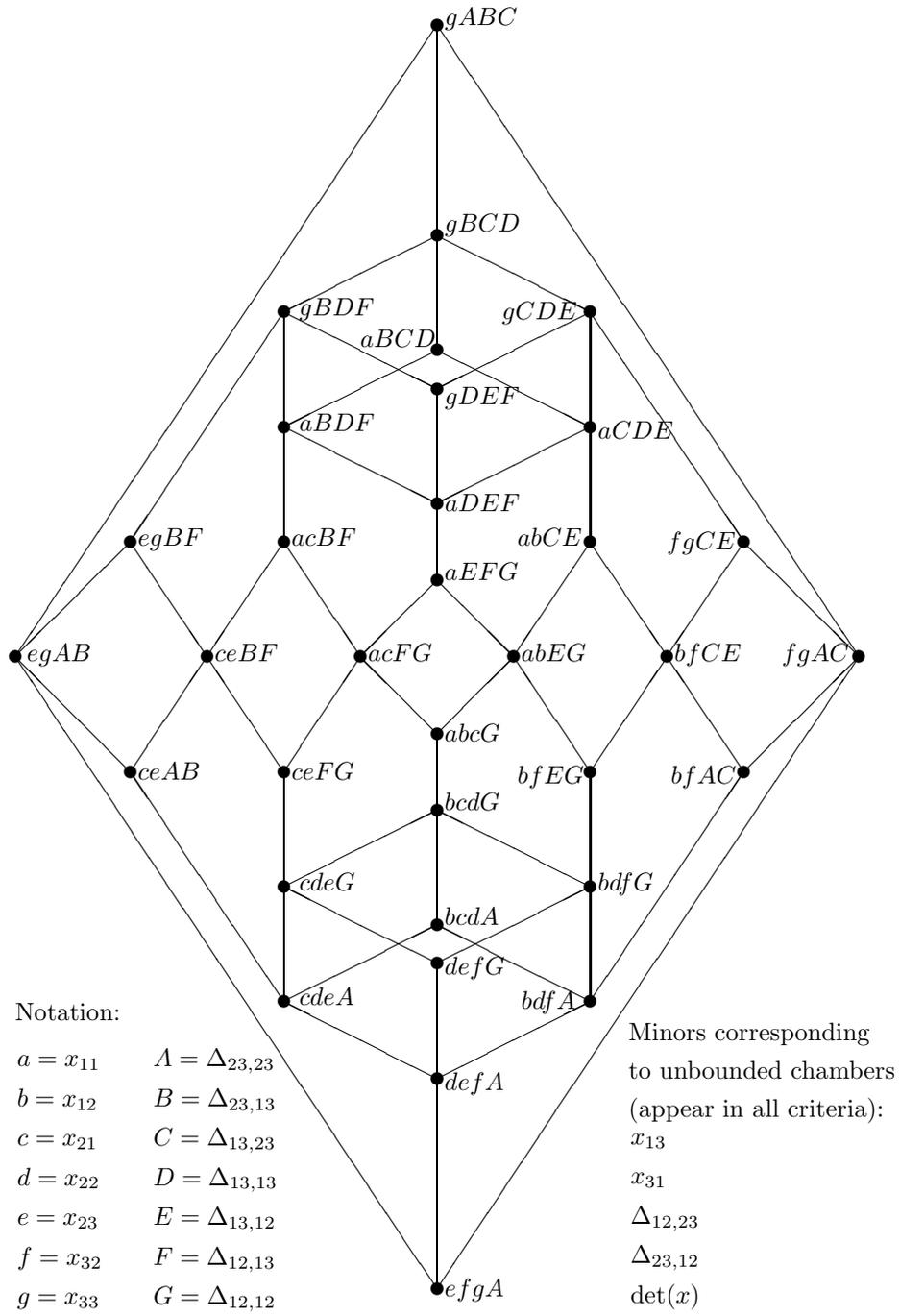
\begin{figure}[t]
\setlength{\unitlength}{1.51pt} 
\begin{center}
\begin{picture}(200,340)(-100,-170)


\put(0,-165){\line(2,3){110}}
\put(0,-165){\line(-2,3){110}}
\put(0,-165){\line(0,1){55}}

\put(0,-110){\line(0,1){30}}
\put(0,-110){\line(2,1){40}}
\put(0,-110){\line(-2,1){40}}

\put(40,-90){\line(2,3){40}}
\put(40,-90){\line(0,1){30}}
\put(40,-90){\line(-2,1){40}}

\put(0,-80){\line(-2,1){40}}
\put(0,-80){\line(2,1){40}}

\put(-40,-90){\line(2,1){40}}
\put(-40,-90){\line(0,1){30}}
\put(-40,-90){\line(-2,3){40}}
 
\put(0,-70){\line(0,1){30}}

\put(40,-60){\line(-2,1){40}}
\put(-40,-60){\line(2,1){40}}
\put(40,-60){\line(0,1){30}}
\put(-40,-60){\line(0,1){30}}

\put(0,-40){\line(0,1){20}}
 
\put(0,-20){\line(1,1){20}}
\put(0,-20){\line(-1,1){20}}
\put(40,-30){\line(2,3){20}}
\put(40,-30){\line(-2,3){20}}
\put(80,-30){\line(1,1){30}}
\put(80,-30){\line(-2,3){20}}
\put(-40,-30){\line(2,3){20}}
\put(-40,-30){\line(-2,3){20}}
\put(-80,-30){\line(2,3){20}}
\put(-80,-30){\line(-1,1){30}}


\put(0,165){\line(2,-3){110}}
\put(0,165){\line(-2,-3){110}}
\put(0,165){\line(0,-1){55}}

\put(0,110){\line(0,-1){30}}
\put(0,110){\line(2,-1){40}}
\put(0,110){\line(-2,-1){40}}

\put(40,90){\line(2,-3){40}}
\put(40,90){\line(0,-1){30}}
\put(40,90){\line(-2,-1){40}}

\put(0,80){\line(-2,-1){40}}
\put(0,80){\line(2,-1){40}}

\put(-40,90){\line(2,-1){40}}
\put(-40,90){\line(0,-1){30}}
\put(-40,90){\line(-2,-3){40}}
 
\put(0,70){\line(0,-1){30}}

\put(40,60){\line(-2,-1){40}}
\put(-40,60){\line(2,-1){40}}
\put(40,60){\line(0,-1){30}}
\put(-40,60){\line(0,-1){30}}

\put(0,40){\line(0,-1){20}}
 
\put(0,20){\line(1,-1){20}}
\put(0,20){\line(-1,-1){20}}
\put(40,30){\line(2,-3){20}}
\put(40,30){\line(-2,-3){20}}
\put(80,30){\line(1,-1){30}}
\put(80,30){\line(-2,-3){20}}
\put(-40,30){\line(2,-3){20}}
\put(-40,30){\line(-2,-3){20}}
\put(-80,30){\line(2,-3){20}}
\put(-80,30){\line(-1,-1){30}}

 
\put(2,-22){$abcG$} 
\put(-18,-1){$acFG$} 
\put(-38,-32){$ceFG$} 
\put(-36,-61){$cdeG$} 
\put(2,-40){$bcdG$} 
\put(42,-61){$bdfG$} 
\put(21,-33){$bfEG$} 
\put(22,-1){$abEG$} 
\put(2,-83){$defG$} 
\put(2,-70){$bcdA$} 
\put(-36,-91){$cdeA$} 
\put(2,-113){$defA$} 
\put(22,-92){$bdfA$} 
\put(2,-167){$efgA$} 
\put(-107,-1){$egAB$} 
\put(-78,-32){$ceAB$} 
\put(-58,-1){$ceBF$} 
\put(62,-1){$bfCE$} 
\put(61,-33){$bfAC$} 
\put(90,-1){$fgAC$} 
      
\put(2,20){$aEFG$}
 
\put(-78,29){$egBF$} 
\put(-38,29){$acBF$}
\put(21,29){$abCE$} 
\put(60,28){$fgCE$} 
 
\put(2,38){$aDEF$} 

\put(-36,59){$aBDF$} 
\put(42,57){$aCDE$}
 
\put(2,66){$gDEF$} 

\put(-20,81){$aBCD$}
 
\put(-36,89){$gBDF$} 
\put(17,88){$gCDE$} 

\put(2,111){$gBCD$} 
\put(2,165){$gABC$} 
      
 
\put(0,-20){\circle*{3}} 
\put(-20,0){\circle*{3}} 
\put(-40,-30){\circle*{3}} 
\put(-40,-60){\circle*{3}} 
\put(0,-40){\circle*{3}} 
\put(40,-60){\circle*{3}} 
\put(40,-30){\circle*{3}} 
\put(20,0){\circle*{3}} 
\put(0,-80){\circle*{3}} 
\put(0,-70){\circle*{3}} 
\put(-40,-90){\circle*{3}} 
\put(0,-110){\circle*{3}} 
\put(40,-90){\circle*{3}} 
\put(0,-165){\circle*{3}} 
\put(-110,0){\circle*{3}} 
\put(-80,-30){\circle*{3}} 
\put(-60,0){\circle*{3}} 
\put(60,0){\circle*{3}} 
\put(80,-30){\circle*{3}} 
\put(110,0){\circle*{3}} 
      
\put(0,20){\circle*{3}} 
 
\put(-80,30){\circle*{3}} 
\put(-40,30){\circle*{3}} 
\put(40,30){\circle*{3}} 
\put(80,30){\circle*{3}} 
 
\put(0,40){\circle*{3}} 

\put(-40,60){\circle*{3}} 
\put(40,60){\circle*{3}} 
 
\put(0,70){\circle*{3}} 

\put(0,80){\circle*{3}} 
 
\put(-40,90) {\circle*{3}} 
\put(40,90){\circle*{3}} 

\put(0,110){\circle*{3}} 
\put(0,165){\circle*{3}}


\put(-110,-95){Notation:}

\put(-115,-138){
$\begin{array}{ll}
a=x_{11} & \quad A=\Delta_{23,23}\\[.05in]
b=x_{12} & \quad B=\Delta_{23,13}\\[.05in] 
c=x_{21} & \quad C=\Delta_{13,23}\\[.05in] 
d=x_{22} & \quad D=\Delta_{13,13}\\[.05in] 
e=x_{23} & \quad E=\Delta_{13,12}\\[.05in] 
f=x_{32} & \quad F=\Delta_{12,13}\\[.05in] 
g=x_{33} & \quad G=\Delta_{12,12} 
\end{array}
$}
 
\put(50,-100){Minors corresponding}
\put(50,-110){to unbounded chambers}
\put(50,-120){(appear in all criteria):}

\put(45,-148){
$\begin{array}{l}
x_{13} \\[.05in] 
x_{31} \\[.05in] 
\Delta_{12,23} \\[.05in] 
\Delta_{23,12} \\[.05in] 
\det(x)
\end{array}
$}


\end{picture}
\end{center}
\caption{
Total positivity criteria for $GL_3$}
\label{fig:schemes}
\end{figure}

In the case of $GL_3\,$ and $u=v=\wnot\,$, there are $34$ isotopy
types, giving rise to 34 different total positivity criteria.
Each of these criteria involves 9 minors.
Five of them---the minors 
\[
x_{31} \,, 
x_{13}\,,
\Delta_{23,12} \,, 
\Delta_{12,23} \,, 
\det(x)
\]
---are common to all 34 criteria;
they correspond to the ``unbounded'' chambers lying on the 
periphery of each arrangement.
The other four minors that distinguish isotopy types from each other 
correspond to the bounded chambers. 
Figure~\ref{fig:schemes} shows a graph with $34$ vertices 
labeled by the quadruples of ``bounded'' minors
that appear in the corresponding total positivity criteria.

For an arbitrary $n$ (and $u=v=\wnot\,$), 
one obtains various nice (and surprising) total positivity
criteria in $GL_n$ by making particular choices of 
(the isotopy type of) a double
pseudoline arrangement in Theorem~\ref{th:GLn-criteria}. 
Let us discuss two criteria obtained in this way. 

A minor $\Delta_{I,J}$ is called \emph{solid} if both $I$ and $J$
consist of several consecutive indices.
A criterion due to Fekete~\cite{fekete} 
(see also \cite[p.~299]{gantmacher-krein}) asserts that
(strict) total positivity of a matrix is equivalent to the
positivity of all its solid minors.
Each of the two criteria described below will strengthen this result.

We will consider two factorization schemes of type
$(\wnot, \wnot)$ having the same $E$- and $F$-parts
(albeit shuffled in a different way).
For both of them, the $E$-part is 
the lexicographically minimal reduced word for $\wnot$,
i.e., the reduced word 
\[
1,2, 1, 3,2,1,\dots,n\!-\!1, n\!-\!2, \dots, 1; 
\]
the $F$-part is the same but with barred entries.
Let $\ii_1$ denote the shuffle of these parts
such that all the unbarred entries precede the barred ones.
Let $\ii_2$ denote the shuffle of the same parts such that
every unbarred entry is immediately followed by the corresponding
barred entry (so that $\ii_2$ starts with 
$1, \overline 1, 2, \overline 2, \dots$).
A direct check shows that the corresponding collections of minors
$F(\ii_1)$ and $F(\ii_2)$ are given as follows:
\begin{itemize}
\item $F(\ii_1)$ consists of solid minors 
$\Delta_{I,J}$ such that $1\in I\cup J$;
\item $F(\ii_2)$ consists of solid minors 
$\Delta_{I,J}$ such that $\min (I) + \max (J)\in \{n,n+1\}$.
\end{itemize}
Each of these two collections consists of $n^2$ minors; 
and by Theorem~\ref{th:GLn-criteria}, 
each of them provides a total positivity criterion that strengthens
the one of Fekete's:
a square matrix is totally positive if and only if all the minors
in $F(\ii_1)$ (respectively, $F(\ii_2)$) are positive. 
It should be mentioned that the first of these criteria was
(implicitly) obtained by Cryer~\cite[Theorems~1.1 and~3.1]{cryer} 
using a result of Karlin~\cite[p.~85]{karlin};
an explicit statement appears in \cite[(3.2)]{gasca-pena}. 
The second criterion seems to be new.


The equivalence of conditions $(2)$ and $(3)$ in 
Theorem~\ref{th:GLn-criteria} has the following algebraic explanation 
(cf.\ Theorem~\ref{th:TP-base}).

\begin{theorem}
\label{th:TPbaseGLn}
For any factorization scheme $\ii$ of type $(u,v)$,
the collection of minors $F(\ii)$ is a totally positive base
(cf.\ Definition~\ref{def:TP-base}) for the collection $F(u,v)$.
\end{theorem}

The most significant part of this theorem is that every minor 
from $F(u,v)$ can be written as a subtraction-free expression
in the minors from $F(\ii)$.
Such an expression can be found in a constructive way.
To do this, it will be enough to consider two arrangements 
$\Arr(\ii)$ and $\Arr(\ii')$ whose isotopy types are adjacent 
in the graph that we described above; recall that this means 
that the collection of minors $F(\ii')$ is obtained from $F(\ii)$ 
by exchanging a single minor $\Delta$ with another minor~$\Delta'$.
It suffices to show that $\Delta'$ can be written as a 
subtraction-free expression in the minors from~$F(\ii)$.
This can be done with the help 
of certain 3-term determinantal identities.
These identities are stated in the following proposition,  
which is a specialization of 
Theorems~\ref{th:minors-Plucker} and ~\ref{th:minors-Dodgson}.
We use the notation $Li$, $Lij$, etc., 
as a shorthand for $L\cup\{i\}$, $L\cup\{i,j\}$, etc.

\begin{proposition}
\label{pr:gln-3-term}
For any $i,j,k,p\in [1,n]$ and $I,L\subset [1,n]$ such that 
$i<j<k$, $|I|=|L|+1$, $L \cap \{i,j,k\} = \emptyset$, 
$I \cap \{p\} = \emptyset$, we have:
\beal
\label{eq:3-term}
\Delta_{Ip,Lik} \Delta_{I,Lj} 
= \Delta_{Ip,Lij} \Delta_{I,Lk} + \Delta_{Ip,Ljk} \Delta_{I,Li} 
\\[.1in]
\Delta_{Lik,Ip} \Delta_{Lj,I} 
= \Delta_{Lij,Ip} \Delta_{Lk,I} 
+ \Delta_{Ljk,Ip} \Delta_{Li,I} \ .
\eea
For any $i,i',j,j'\in\{1,\dots,n\}$
and $I,J\subset\{1,\dots,n\}$ such that 
$i < i'$, $j < j'$, $|I|=|J|$, 
$I \cap \{i,i'\} = J \cap \{j,j'\} = \emptyset$, we have 
\begin{equation}
\label{eq:carroll}
\Delta_{Ii,Jj} \Delta_{Ii',Jj'} = \Delta_{Ii,Jj'} \Delta_{Ii',Jj} +
\Delta_{I,J} \Delta_{Iii',Jjj'}\ .
\end{equation}
\end{proposition}

The identities~(\ref{eq:3-term})--(\ref{eq:carroll})
are well known, although their attribution
is complicated. As early as in 1819 they were proved by P.~Desnanot
(see~\cite[pp.~140-142]{muir}). 
Identities~(\ref{eq:3-term})
are special cases of the (Grassmann-)Pl\"ucker relations
(see, e.g.,~\cite[(15.53)]{fulton-harris}),
while identity~(\ref{eq:carroll}) 
plays a crucial role in C.~L.~Dodgson's condensation method,
and is because of that occasionally associated with the name of Lewis
Carroll.

It would be interesting to see which other classical determinantal identities
can be generalized to the functions $\Delta_{u \omega_i, v \omega_i}$ on any 
semisimple group. 
We conclude the paper by mentioning one challenging problem 
of this kind: find a generalization of the classical  
Binet-Cauchy formula for the minors of the product of two matrices:
$$\Delta_{I,J} (xy) = \sum_K \Delta_{I,K} (x) \Delta_{K,J} (y)  \ .$$
At present, we only know such a generalization for the \emph{minuscule} 
fundamental weights~$\omega_i\,$.


\end{document}